\documentclass[11pt,english]{article}
\usepackage[english]{babel}
\usepackage{amsfonts}
\usepackage{epsfig}
\usepackage{graphicx}
\usepackage{fancyhdr}
\usepackage{amssymb}
\usepackage{amsmath}
 \usepackage{verbatim} 
 \usepackage{appendix} 
\usepackage{xcolor}

\newtheorem{theorem}{Theorem}

\newtheorem{corollary}{Corollary}[section]

\newtheorem{prop}{Proposition}[section]
\newtheorem{lemma}[prop]{Lemma}
\newtheorem{definition}[prop]{Definition}
\newtheorem{cl}[prop]{Claim}
\newtheorem{rem}[prop]{Remark}

\def\div{{\rm div}\,}% divergence
\def\Dg {{\cal D}} %D gothique
\def\Ig {{\cal I}} %I gothique
\def\Lg {{\cal L}} %L gothique
\def\Vg {{\cal V}} %V gothique
\def\and {{\rm \; and \;}}

\def\Im {{\rm \; Im\;}}
\def\Re {{\rm \;Re\;}}
\def\exp {{\rm exp}}

\def\Dg {{\cal D}} %D gothique
\def\Ig {{\cal I}} %I gothique
\def\Lg {{\cal L}} %L gothique
\def\Vg {{\cal V}} %V gothique
\def\and {{\rm \; and \;}}

\def\Im {{\rm \; Im\;}}
\def\Re {{\rm \;Re\;}}
\def\exp {{\rm exp}}

\def\T0{T_{0,1}}

\newcommand {\R}{ \mathbb{R}}
\newcommand {\C}{ \mathbb{C}}
\newcommand {\N}{ \mathbb{N}}

\newcommand {\pa}{\partial}

\newcommand {\beqna} {\begin{eqnarray}}
\newcommand {\eeqna} {\end{eqnarray}}
\newcommand {\beqtn} {\begin{equation}}
\newcommand {\eeqtn} {\end{equation}}

\newcommand {\dsp}{\displaystyle}

\begin{document}

\title{\bf Construction of a blow-up solution for the Complex Ginzburg-Landau equation in some critical case}
\author{Nejla Nouaili,\\ {\it \small Nejla.Nouaili@dauphine.fr}\\ {\it \small CEREMADE, Universit\'e Paris Dauphine, Paris Sciences et Lettres.} \\
Hatem Zaag,\\
%\footnote{This author is supported by the ERC Advanced Grant  no. 291214, BLOWDISOL, and by ANR project no. ANR-13-BS01-0010-03, ANA\'E.},\\
{\it \small Hatem.Zaag@univ-paris13.fr}\\{\it \small Universit\'e Paris 13, Sorbonne Paris Cit\'e,}\\ {\it \small LAGA, CNRS (UMR 7539), F-93430, Villetaneuse, France.}}

\maketitle
\begin{abstract}

We construct a solution for the Complex Ginzburg-Landau (CGL) equation in some critical case, which blows up in finite time $T$ only at one blow-up point. We also give a sharp description of its profile. The proof relies on the reduction of the problem to a finite dimensional one, and the use of index theory to conclude. The interpretation of the parameters of the finite dimension problem in terms of the blow-up point and time allows to prove the stability of the constructed solution.

\end{abstract}

\section{Introduction}

We consider the following complex Ginzburg-Landau (CGL) equation
\beqtn
\begin{array}{l}
u_t=(1+i\beta)\Delta u+(1+i\delta) |u|^{p-1}u+\nu u,\\
u(.,0)=u_0\in L^\infty (\R^N,\C)
\end{array}
\mbox{     (CGL)}
\label{GL}
\eeqtn
where $\alpha,\beta,\nu\in\R$.

\medskip

In fact, in this paper we will treat only the case $\beta=0$, which is an important case in the CGL physics literature (see Aranson and Kramer \cite{AKRMP02}) and which of course can be seen as a complex valued semi-linear heat equation with no variational structure. However, in order to connect our work both to the physics and the math literature, we will consider the case $\beta\in \R$ in the introduction.

This equation, most often considered with a cubic nonlinearity ($p=3$), has a long history in physics (see Aranson and Kramer \cite{AKRMP02}). The Complex Ginzburg-Landau (CGL) equation  is one of the most-studied equations in physics, it describes a lot of phenomena including nonlinear waves, second-order phase transitions, and superconductivity. We note that the Ginzburg-Landau equation can be used to describe the evolution of amplitudes of unstable modes for any process exhibiting a Hopf bifurcation (see for example Section VI-C, page 37 and Section VII, page 40 from \cite{AKRMP02} and the references cited therein). The equation can be considered as a general normal form for a large class of bifurcations and nonlinear wave phenomena in continuous media systems. More generally, the Complex Ginzburg-Landau (CGL) equation is used to describe synchronization and collective oscillation in complex media.

\medskip

The study of blow-up, collapse or chaotic solutions of equation \eqref{GL} appears in many works; in the description of an unstable plane Poiseuille flow, see Stewartson and Stuart \cite{SSJFM71}, Hocking, Stewartson, Stuart and Brown \cite{HSSBJFM72} or in the context of binary mixtures Kolodner and \textit{al}, \cite{KBSPRL88}, \cite{KSALPNP95}, where the authors describe an extensive series of experiments on traveling-wave convection in an ethanol/water mixture, and they observe collapse solution that appear experimentally.\\
For our purpose, we consider CGL independently from any particular physical context and investigate it as a mathematical model in partial differential equations with $p>1$.
 
\medskip

We note also that the interest on the study of singular solutions in CGL comes also from the analogies with the three-dimensional Navier-Stokes. The two equations have the same scaling properties and the same energy identity (for more details see the work of Plech{\'a}{\v{c}} and {\v{S}}ver{\'a}k \cite{PSCPAM01}; the authors in this work give some evidence for the existence of a radial solution which blow up in a self-similar way). Their argument is based on matching a numerical solution in an inner region with an analytical solution in an outer region.  In the same direction we can also  cite the work of  Rottsch\"afer \cite{RPD08} and \cite{REJAM13}.

\bigskip

\medskip

The Cauchy problem for equation \eqref{GL} can be solved in a variety of spaces using the semigroup theory as in the case of the heat equation (see \cite{CNYUCIM03,GV96,GVCMP97}).\\
We say that $u(t)$ blows up or collapse in finite time $T<\infty$, if $u(t)$ exists for all $t\in [0,T)$ and 
$\lim_{t\to T}\|u(t)\|_{L^\infty}=+\infty.$
In that case, $T$ is called the blow-up time of the solution. A point $x_0\in\R^N$ is said to be a blow-up point if there is a sequence $\{(x_j,t_j)\}$, such that $x_j\to x_0$, $t_j\to T$ and $|u(x_j,t_j)|\to \infty$ as $j\to\infty$. The set of all blow-up points is called the blow-up set.\\
Let us now introduce the following definition;
\begin{definition} The exponents $(\beta,\delta)$ are said to be critical  (resp. subcritical, resp. supercritical) if 
$p-\delta^2-\beta\delta (p+1)=0$ (resp. $>0$, resp. $<0$).
\label{subcond}
\end{definition}
An extensive literature is devoted to the blow-up profiles for CGL when $\beta=\delta=0$ (which is the nonlinear heat equation),  see Vel{\'a}zquez \cite{VCPDE92, VTAMS93,VINDIANA93} and Zaag \cite{ZIHP02, ZCMP02, ZMME02} for partial results). In one space dimension, given $a$ a blow-up point, this is the situation:

\medskip

\begin{itemize}

\item either
\beqtn
\sup_{|x-a|\leq K\sqrt{(T-t)\log(T-t)}}\left|(T-t)^{\frac{1}{p-1}}u(x,t)-f\left(\frac{x-a}{\sqrt{(T-t)|\log(T-t)}}\right)|\right|\to 0,
\label{u2}
\eeqtn
\item or for some $m\in \N$, $m\ge 2$, and $C_m>0$
 \beqtn
\sup_{|x-a|<K(T-t)^{1/2m}}\left|(T-t)^{\frac{1}{p-1}}u(x,t)-f_m\left(\frac{C_m (x-a)}{(T-t)^{1/2m}}\right)\right|\to 0,
\label{um}
\eeqtn
as $t\to T$, for any $K>0$, where 
\begin{equation}
\label{definitionf}
\begin{array}{l}
f(z)=\left(p-1+b_0z^2\right)^{-\frac{1}{p-1}}
\mbox{ where }b_0=\frac{(p-1)^2}{4p},\\
f_m(z) =\left(p-1+|z|^{2m}\right)^{-\frac{1}{p-1}}.
\end{array}
\end{equation}
\end{itemize}

\medskip

If $(\beta,\delta)\not =( 0,0)$, the situation is completely understood in the  \textit{subcritical} case by Zaag \cite{ZAIHPANL98} ($\beta=0$) and Masmoudi and Zaag \cite{MZ07} ($\delta\not = 0$). More precisely, if
\[p-\delta^2-\beta\delta (p+1)>0,\] 
then, the authors construct a solution of equation \eqref{GL}, which blows up in finite time $T>0$ only at the origin such that for all $t\in [0,T)$,
\beqtn
\dsp\left\| (T-t)^{\frac{1+i\delta}{p-1}}\left |\log(T-t)\right |^{-i\mu}u(x,t)-\left(p-1+\frac{b_{sub}|x|^2}{(T-t)|\log(T-t)|}\right)^{-\frac{1+i\delta}{p-1}}
\right\|_{L^\infty}\leq \frac{C_0}{1+\sqrt{|\log(T-t)|}}
\label{profilesubc}
\eeqtn
where
\beqtn
\dsp b_{sub}=\frac{(p-1)^2}{4(p-\delta^2-\beta\delta(1+ p))}>0\mbox{ and }\mu=-\frac{2b_{sub}\beta}{(p-1)^2}(1+\delta^2).
\label{bs}
\eeqtn
Note that this result was previously obtained formally by Hocking and Stewartson \cite{HSPRSLS72} ($p=3$) and mentioned later in Popp and \textit{al} \cite{PSKKPD98}.

\medskip

When $p=3$,  many works has been devoted to the blow-up profile for the CGL in the  \textit{subcritical} case (see Definition \ref{subcond}). We cite the works of Hocking and Stewartson \cite{HSPRSLS72} and Popp and \textit{al} and the references cited therein \cite{PSKKPD98}. 

\medskip

In the  \textit{critical} and \textit{supercritical} cases, few results are known about chaotic solutions for the equation. We note that  Hocking and Stewartson \cite{HSPRSLS72} and  Popp and \textit{al} \cite{PSKKPD98} proved, formally, the existence of at least two self-similar blow-up solutions, one of them is the same solution constructed by Masmoudi and Zaag in the  \textit{subcritical} case \eqref{profilesubc}. The second kind  of blow-up is formally described in \cite{HSPRSLS72} and \cite{PSKKPD98} and approved by numerical results. 

\medskip

 Let us now give the formal result, when $p=3$, given by Popp and \textit{al} \cite{PSKKPD98} (equations (44) and (64)) in the  \textit{critical} case ($3-4\delta\beta-\delta^2=0$): the authors obtained 
\beqtn
(T-t)^{\frac{1+i\delta}{p-1}}u(x,t)\sim e^{i\psi(t)}  \left( 2+ \frac{b_p|x|^2}{(T-t)|\log(T-t)|^{\frac 12}}\right)^{-\frac{1+i\delta}{2}}.
\label{propal}
\eeqtn
where 
\beqtn
b_{p}=2\left( \sqrt{\frac {3}{2 \delta^2}(\delta^2+5)(\delta^2+1)(15-\delta^2)}\right)^{-1}.
\label{defbp}
\eeqtn
and $\psi(t)$ is given by equation (40) in \cite{PSKKPD98}. We can clearly see, that this profile exist only for $\delta^2<15$.
\begin{rem}
We will see later, in Section \ref{sectFormapp}, that we obtain the same result with another formal approach.
\end{rem}
\medskip

 In this paper, we justify the formal result of Popp and al \cite{PSKKPD98}, in the case $p=3$, and construct a solution $u(x,t)$ of \eqref{GL} in the critical case ($\beta=0$ and $\delta^2=p$) that blows up in some finite time $T$, in the sense that
\[\dsp \lim_{t\to T}\|u(.,t)\|_{L^\infty}=+\infty.\]
More precisely, this is our result:
\begin{theorem}[Blow-up profiles for equation \eqref{GL}] Consider $\beta=0$, $\delta^2=p$, then equation \eqref{GL} has a solution $u(x,t)$, which blows up in finite time $T$, only at the origin. Moreover:
(i) for all $t\in [0,T)$,
\beqtn
\left\|(T-t)^{\frac{1+i\delta}{p-1}}|\log (T-t)|^{-i\mu} u(x,t)-\varphi_0\left(\frac{x}{\sqrt{(T-t)}|\log (T-t)|^{1/4}}\right) \right\|_{L^\infty(\R^N)} \leq \frac{C_0}{1+|\log(T-t)|^{\frac 14}}
\label{Preslt}
\eeqtn
where
\beqtn
%\kappa ^{-i\delta}
\varphi_0(z)=\left ( p-1+b z^2\right )^{-\frac{1+i\delta}{p-1}},
\label{deffi0}
\eeqtn
\beqtn
b=\frac{(p-1)^2}{8\sqrt{p(p+1)}},\mu=\frac{\delta}{8p}.
\label{defbmu}
\eeqtn
(ii) For all $x\not =0$, $u(x,t) \to u^*(x)\in C^2(\R^N\backslash\{0\})$ and
\beqtn
\dsp u^{*}(x) \sim |2\log |x||^{i\mu}\left [ \frac{b|x|^2}{\sqrt{2|\log |x||}}\right]^{-\frac{1+i\delta}{p-1}}\mbox { as }x\to 0.
\label{finalprofile}
\eeqtn 
\label{thm1}
\end{theorem}

\medskip

\begin{rem} We will consider CGL, given by \eqref{GL}, only when $\nu=0$. The case $\nu\not =0$ can be done as in \cite{EZSMJ11}. In fact, when $\nu\not =0$, exponentially small terms will be added to our estimate in self-similar variable (see \eqref{chauto} below), and that will be absorbed in our error terms, since our trap $\Vg_A(s)$ defined in Definition \ref{defthess} is given in polynomial scales.
\label{nu0}
\end{rem}

\begin{rem} 

The derivation of the blow-up profile \eqref{deffi0} can be understood through a formal analysis, using  the matching asymptotic expansions (see Section \ref{sectFormapp} below). 
This method was used by  Galaktionov, Herrero and Vel\'azquez \cite{VGH91} to derive all the possible behaviors of the blow-up solution given by (\ref{u2},\ref{um}) in the heat equation ($\beta=\delta=0$). This formal method was used recently by Tayachi an Zaag \cite{TZ17} in the case of a nonlinear heat equation with a critical power nonlinear gradient term.\\
However, we would like to emphasize on the fact that our formal analysis is far from being a simple adaptation of this work.  We will see in Section \ref{sectFormapp} that we need much more effort to obtain the profile, this is due to the criticality of the problem ($p=\delta^2$).

\end{rem}
\begin{rem}
The exhibited profile  \eqref{deffi0} is new in two respects: 
\begin{itemize}
\item The scaling law in the critical case is $\sqrt{(T-t) |\log (T-t)|^{\frac 12}}$ instead of the laws of subcritical case, $\sqrt{(T-t)| \log (T-t)|}$. 
%or  \textcolor{violet}{$(T-t)^{\frac{1}{2m}}$}, where $m\geq 2$ is an integer.
\item The profile function: $\varphi_0(z)=(p-1+b |z|^2)^{-\frac{1+i\delta}{p-1}}$ is different from the profile of the subcritical case , namely $f(z)=(p-1+b_{sub} |z|^2)^{-\frac{1+i\delta}{p-1}}$, in the sense that  $b\not = b_{sub}$ (see \eqref{bs} and \eqref{defbmu}).
\end{itemize}
\end{rem}
\begin{rem}
In the subcritical case $p-\delta^2-\beta\delta (p+1)>0 $ ($p-\delta^2>0$, when $\beta=0$), the final profile of the CGL is given by 

\[u^*(x)\sim |2\log |x||^{i\mu}\left [ \frac b 2 \frac{|x|^2}{|\log|x||}\right]^{-\frac{1+i\delta}{p-1}}\mbox{ as $x\to 0$  }\]
with
\[b=\frac{(p-1)^2}{4(p-\delta^2-\beta\delta(p+1))}\mbox{ and }\mu=-\frac{2b\beta}{(p-1)^2}(1+\delta^2).\]
In the critical case $\beta=0,\delta^2=p$, the final profile is given by \eqref{finalprofile}.
\end{rem}

\begin{rem}We strongly believe that our construction will not work for all $(\beta,\delta)$ satisfying the critical condition . This is due to the formal study given by Popp and \textit{al} \cite{PSKKPD98} in the cubic case. %From the formal analysis done in \cite{PSKKPD98}, the authors obtain the profile \beqtn{propal}. and 
We can see from the definition \eqref{defbp} that $b_p$ can be defined only if $\delta^2<15$.
\end{rem}

As a consequence of our techniques, we show the stability of the constructed solution with respect to perturbations in initial data. More precisely, we have the following result.

\begin{theorem}[Stability of the solution constructed in Theorem \ref{thm1} ($\beta=0$)] Let us denote by $\hat u(x,t)$ the solution constructed in Theorem 1 and by $\hat T$ its blow-up time. Then, there exists a neighborhood $\Vg_0$ of $\hat u(x,0)$ in $L^\infty$ such that for any $u_0\in\Vg_0$, equation \eqref{GL} has a unique solution $u(x,t)$ with initial data $u_0$, and $u(x,t)$ blows up in finite time $T(u_0)$ at one single blow-up point $a(u_0)$. Moreover estimate \eqref{Preslt} is satisfied by $u(x-a,t)$ and
\[T(u_0)\to \hat T,\;\; a(u_0)\to 0\mbox{ as }u_0\to \hat u_0\mbox{ in }L^\infty(\R^N,\C).\]
\label{thstablity}
\end{theorem}

\begin{rem}
We will not give the proof of Theorem \ref{thstablity}
because the stability result follows from the reduction to a finite dimensional case as in \cite{MZDuke97}  (see Theorem 2 and its proof in Section 4) 
and \cite{MZ07} (see Theorem 2 and its proof in Section 6) with the same argument. Hence, we only prove the existence result (Theorem \ref{thm1}) and kindly refer the reader to \cite{MZDuke97} and \cite{MZ07} for the proof of the stability.
\end{rem}

\medskip

Let us give an idea of the method used to prove the results. We construct the blow-up solution with the profile in Theorem \ref{thm1}, by following the method of \cite{MZDuke97}, \cite{BKN94}, though we are far from a simple adaptation, since we are studying the critical problem, which make the technical details harder to elaborate. This kind of methods has been applied for various nonlinear evolution equations. For hyperbolic equations, it has been successfully used for the construction of multi-solitons for semilinear wave equation in one space dimension (see \cite{CZCPAM13}). For parabolic equations, it has been used in \cite{MZ07} and \cite{ZCPAM01} for the Complex Ginzburg Landau (CGL) equation with no gradient structure, the critical harmonic heat flow in  \cite{RSCPAM13}, the two dimensional Keller-Segel equation in \cite{RSMA14} and the nonlinear heat equation involving nonlinear gradient term in \cite{EZSMJ11} and \cite{TZ17}. Recently, this method has been applied for a non variational parabolic system in  \cite{NZCPDE15} for a logarithmically perturbed nonlinear equation in  \cite{NVTZ15}.

\medskip

Unlike in subcritical case \cite{MZ07} and \cite{ZCPAM01}. the criticality of the problem induces substantial changes in the blow-up profile as pointed out in the comments following Theorem \ref{thm1}. Accordingly, its control requires special arguments .So, working in the framework of \cite{MZDuke97} and \cite{MZ07}, some crucial modifications are needed. In particular, we have overcome the following challenges:
\begin{itemize}
\item The prescribed profile is not known and not obvious to obtain. See Section \ref{sectFormapp} for a formal approach to justify such a profile.
\item The profile is different from the profile in \cite{MZDuke97} and \cite{MZ07}, therefore new estimates are needed.
\item In order to handle the new scaling, we introduce a new shrinking set to trap the solution. See Definition \ref{defthess}. Finding such set is not trivial, in particular in the critical case, where we need much more details in the expansions of the rest term see Appendix \ref{Detailscal}.
\end{itemize}
Then, following \cite{MZDuke97}, the proof is divided in two steps. First, we reduce the problem to a finite dimensional case. Second, we solve the finite time dimensional problem and conclude by contradiction using index theory. 
%CORRECTION RQ 8, p1
%Here we follow the method developed by Merle and Zaag \cite{MZDuke97} for the construction of a stable blowup solution. 
The proof is performed in the framework of the similarity variables defined below in \eqref{chauto}. 
We linearize the self-similar solution around the profile $\varphi_0$ and we obtain $q$ (see \eqref{definitionq} below). Our goal is to guaratee that $q(s)$ belongs to some set $\Vg_A(s)$ (introduced in Definition \ref{defthess}), which shrinks to $0$ as $s\to +\infty$. The proof relies on two arguments:

\begin{itemize}
\item The linearized equation gives two positives mode; $q_0$ and $q_1$, one zero modes ($q_2$) and an infinite dimensional negative part. The negative part is easily controlled by the effect of the heat kernel. The control of the zero mode is quite delicate (see Part 2: Proof of Proposition \ref{propode}, page \pageref{delicatq}). Consequently, the control of $q$ is reduced to the control of its positive modes.

\item The control of the positive modes $q_0$ and $q_1$ is handled thanks to a topological argument based on index theory (see the argument at page \pageref{indextheory}). 
\end{itemize}
%CORRECTION

\medskip

The organization of the rest of this paper is as follows. In Section \ref{sectFormapp}, we explain formally how we obtain the profile. In Section \ref{formpb}, we give a formulation of the problem in order to justify the formal argument. Section \ref{existence} is divided in two subsections. In Subsection \ref{pwtr} we give the proof of the existence of the profile assuming technical details. In particular, we construct a shrinking set and give an example of initial data giving rise to prescribed blow-up profile. Subsection \ref{tecnicsect} is devoted to the proof of technical results which are needed in the proof of existence. Finally, in Section \ref{pth1}, we give the proof of Theorem \ref{thm1}.

\medskip

\textbf{Acknowledgment:} The authors would like to thank the referees for their valuable suggestions which (we hope) made our paper much clearer and reader friendly.

\section{Formal approach }\label{sectFormapp}
We consider CGL, given by \eqref{GL}, when $\nu=0$, as we mentioned before in Remark \ref{nu0}.\\
The aim of this section is to explain formally how we derive the behavior given in Theorem \ref{thm1}. In particular, how to obtain the profile $\varphi_0$ in \eqref{deffi0}, the parameter $b$ in \eqref{defbmu}. Consider an arbitrary $T$ and the self-similar transformation of \eqref{GL}
\beqtn
w(y,s)=(T-t)^{\frac{1+i\delta}{p-1}}u(x,t),\;\;y=\frac{x}{\sqrt{T-t}},\;\; s=-\log(T-t).
\label{chauto}
\eeqtn
% Let us introduce $v(y,s)$
% \[w(y,s)=\kappa+v(y,s).\]

If $u(x,t)$ satisfies \eqref{GL} for all $(x,t)\in \R^N\times [0,T)$, then $w(y,s)$ satisfies for all $(x,t)\in \R^N\times [-\log T,+\infty)$ the following equation
\beqtn
\frac{\pa w}{\pa s}=\Delta w-\frac 1 2 y.\nabla w-\frac{1+i\delta}{p-1}w+(1+i\delta)|w|^{p-1}w,
\label{GLauto}
\eeqtn
for all $(y,s)\in \R^N\times [-\log T,+\infty)$. Thus constructing a solution $u(x,t)$ for the equation \eqref{GL} that blows up at $T<\infty$ like $(T-t)^{-\frac{1}{p-1}}$ reduces to constructing a global solution $w(y,s)$ for equation \eqref{GLauto} such that
\beqtn
0<\varepsilon\leq \lim_{s\to\infty} \|w(s)\|_{L^\infty(\R^N)}\leq \frac 1 \varepsilon.
\eeqtn
A first idea to construct a blow-up solution for \eqref{GL}, would be to find a stationary solution of \eqref{GLauto}, yielding a self-similar solution for \eqref{GL}. It happens, that in the subcritical case, the second author together with Masmoudi were able in \cite{MZ07} to construct such a solution. 
%Correction RQ 9, p1
Of course, the construction of \cite{MZ07} is in the same spirit as the approach used for the heat equation in \cite{MZDuke97} and \cite{BKN94}, in the sense that the authors in \cite{MZ07} cut in space the solution on similarity variables. However, the profile they get is different from the case of the standard heat equation, in the sense that they obtain blow-up for the mode and also for the phase, with a different constant in the profile (see $b_{sub}$ \eqref{bs}).   

%Correction RQ 9, p1
In the critical case ($\beta=0,p=\delta^2$), there is no self-similar solution apart from the trivial constant solution $w\equiv \kappa$ of \eqref{GLauto} (see Proposition 2.1 in \cite{NZTAMS08} or simply multiply equation \eqref{GLauto} by $w\rho$, where
\beqtn
\rho(y)=\frac{e^{-\frac{|y|^2}{4}}}{(4\pi)^{N/2}},\; y\in \R^N
\label{defrho}
\eeqtn
and integrate in space).
 
%Correction RQ 9, p1

%Correction RQ 9, p1
\subsection{Inner expansion.} Following the approach of Bricmont and Kupiainen \cite{BKN94} and  Masmoudi and Zaag \cite{MZ07}, we may look for a solution $w$ to \eqref{chauto} such that $v\to 0$ as $s\to\infty$, if $v$ is defined by the following ansatz;
\beqtn w(y,s)=e^{i\mu \log s}(v(y,s)+\kappa)\mbox{ where }\kappa=(p-1)^{-\frac{1}{p-1}}.
\label{eqwmu}
\eeqtn
Using \eqref{chauto}, we see that $v$  should satisfy the equation

\beqtn
\frac{\pa v}{\pa s}=\tilde\Lg v+F(v)-i\frac \mu s \kappa-i\frac \mu s v
\label{eqv}
\eeqtn
with
\beqtn
\tilde \Lg v=\Delta v-\frac 1 2  y.\nabla v+ (1+i\delta)v_1\mbox{, where }v_1=\Re (v),
\eeqtn
\beqtn
F(v)=(1+i\delta)\left( |v+\kappa|^{p-1}(v+\kappa)-\kappa^p-\frac{v}{p-1}-v_1\right).
\eeqtn

Let us recall some properties of $\tilde \Lg$. The operator $\tilde \Lg$ is a $\R-$linear operator defined on $L^{2}_{\rho}(\R^N,\C)$ where

\[\dsp L^{2}_{\rho}(\R^N,\C)=\left\{f\in L^{2}_{loc}(\R^N,\C)|\; \int_{\R^N}|f(y)|^2 \rho(y)dy <\infty    \right\}\] 
and $\rho(y)$ is defined by \eqref{defrho}.

The spectrum of $\tilde \Lg$ is explicitly given by 
\[\{1-\frac m2|\; m\in \N\}.\]
For $N=1$, the eigenfunctions are given by 
\[\{(1+i\delta) h_m(y), ih_m(y)|m\in\N\},\]
where $h_m$ are rescaled versions of Hermite polynomials given by:
\beqtn
\dsp h_m(y)=\sum_{n=0}^{[\frac m 2]}\frac{m!}{n!(m-2n)!}(-1)^ny^{m-2n}.
\label{eigfun}
\eeqtn
Note that we have
\beqtn
\left\{
\begin{array}{lll}
\tilde\Lg ((1+i\delta)h_m)&=&\dsp \left(1-\frac m 2\right)(1+i\delta)h_m,\\
\tilde\Lg (ih_m)&=&-\frac m2 ih_m.
\end{array}
\right .
\label{spectLtilde}
\eeqtn

Since the family $h_m$ is a basis of the Hilbert space $L^2_\rho(R^N, R)$, we see that each $r \in L^{2}_{\rho} $ can be uniquely decomposed as
\beqtn
\label{decompr}
r(y)=(1+i\delta) R_1 (y)+i R_2(y)=\dsp (1+i\delta) \left(\sum_{m=0}^{+\infty}R_{1m}h_m\right)+i \left(\sum_{m=0}^{+\infty}R_{2m}h_m\right),
\eeqtn
where
\beqtn
\begin{array}{l}
\dsp R_1(y)=\Re (r(y))\mbox{,          }R_2(y)=\Im(r(y))-\delta \Re (r(y)),\\
\dsp \mbox{and  for }i={1,2},\;\; R_{im}(y)=\dsp\int _{\R} R_i(y)\frac{h_m(y)}{\|h_m(y)\|_{L^{2}_{\rho}}^{2}}\rho(y) dy.
\end{array}
\eeqtn
Note that $\tilde \Lg$ is not $\C -$linear. In addition it is not self-adjoint, but fortunately it can be written as follows in a diagonal form. Indeed, using the notation \eqref{decompr}, we see that
\[
\tilde \Lg r=(1+i\delta) \Lg R_1+i (\Lg -\Ig)R_2, 
\]
where $\Lg h=\Delta h-\frac 12 y\cdot \nabla h+h$.\\

\medskip

\noindent In compliance with the spectral properties of $\tilde \Lg$, we may look for a solution expanded as follows:
\[v(y,s)=\displaystyle(1+i\delta)\sum_{m\in \N} \bar v_{m}(s) h_{m}(y)+ i \sum_{m\in \N} \hat v_{m}(s) h_{m}(y).\]
Since the eigenfunctions for $m\geq 3$ correspond to negative eigenvalues of $\tilde \Lg$, assuming $v$ is even, we may  look for a solution $v(y,s)$ for equation \eqref{eqv} in the form:
\beqtn
v(y,s)=(1+i\delta)\left(\bar v_0 h_0+\bar v_2 h_2\right)+ i \hat v_0 h_0(y).
\label{formaldecomp}
\eeqtn
Then projecting equation \eqref{eqv}, we derive the following ODE system
%Correction Rq 4 p3
\beqtn
\left \{
\begin{array}{lll}
\bar v_{0}^{'}&=& \displaystyle\bar v_0+\frac{\mu\delta}{s} \bar v_0+\frac{\mu}{s} \hat v_0+\frac{1}{2\kappa}\hat v_{0}^{2}-\frac{(p+1)p}{3\kappa^2} \bar v_{0}^{3}-\frac{(p+1)\delta}{\kappa^2} \bar v_{0}^{2}\hat v_0\displaystyle-\frac{(p+1)}{\kappa^2} \bar v_{0}\hat v_{0}^{2}-8\frac{(p+1)}{\kappa^2} \bar v_{0} \bar v_{2}^{2}\\

&&\displaystyle-\frac{\delta }{2\kappa^2}\hat v_{0}^{3}-8\frac{(p+1)\delta}{\kappa^2}\hat v_0 \bar v_{2}^{2}-\frac{64}{3}\frac{(p+1)p}{\kappa^2} \bar v_{2}^{3}+R_1,\\

\bar v_{2}^{'}&=& \frac{\mu\delta}{s} \bar v_2-40\displaystyle\frac{(p+1)p}{\kappa^2} \bar v_{2}^{3}-8\frac{(p+1)p}{\kappa^2} \bar v_{2}^{2}\bar v_0-8\frac{(p+1)\delta}{\kappa^2} \bar v_{2}^{2}\hat v_0\\
&&\displaystyle -\frac{(p+1)p}{\kappa^2} \bar v_{2}\bar v_{0}^{2}-\frac{(p+1)}{\kappa^2} \bar v_{2}\hat v_{0}^{2}-2\frac{(p+1)\delta}{\kappa^2} \bar v_{2}\bar v_0\hat v_{0}+R_1\\
%\displaystyle\frac{p-\delta^2}{2\kappa} \bar v_{0}\bar v_2+4\frac{p-\delta^2}{\kappa} \bar v_{2}^{2}+R,\\
\hat v_{0}^{'}&=& -\frac{\mu\kappa}{s}-\frac{\mu (1+p)}{s}\bar v_0 -\frac{\mu\delta}{s}\hat v_0+
\displaystyle\frac{1+p}{\kappa} \hat v_{0}\bar v_0+\frac{(1+p)\delta}{\kappa} \bar v_{0}^{2}+\frac{8(1+p)\delta}{\kappa} \bar v_{2}^2\\
&&\displaystyle -\frac{\delta (p+1)^2}{\kappa^2} \bar v_{0}^{3}+\frac{(2p-1)(p+1)}{\kappa^2}\bar v_{0}^{2}\hat v_0+\frac{3\delta(p+1)}{2\kappa^2}\bar v_0\hat v_{0}^{2}+\frac{(p+1)}{2\kappa^2}\hat v_{0}^{3},\\
&&\displaystyle -24\frac{\delta (p+1)^2}{\kappa^2}\bar v_0\bar v_{2}^{2}+8\frac{(2p-1)(p+1)}{\kappa^2}\hat v_0 \bar v_{2}^{2}-64 \frac{\delta (p+1)^2}{\kappa^2}\bar v_{2}^{3}+R_1,
\end{array}
\right .
\label{sysdn}
\eeqtn

where $R_1= \displaystyle O(|\bar v_0|^4+|\hat v_0|^4+|\bar v_2|^4)$.\\%Correction Rq 4 p3
Looking for a solution satisfying the ansatz
\beqtn
\bar v_2\sim\dsp\frac{\alpha}{\sqrt{s}},\;\; \hat v_0 \sim \dsp\frac{\beta}{\sqrt{s}}\mbox{ and }\bar v_0 =O(\frac 1s)\mbox{, where }\alpha,\beta\in \R,
 \label{hypsys}
 \eeqtn
we obtain by the last equation of \eqref{sysdn}

\beqtn
\mu =\frac{8\delta (p+1)}{\kappa^2}\alpha^2.
\label{hypsys1}
\eeqtn
%NEW OCT 05
Note that the ansatz \eqref{hypsys}  seems to be compatible with our ODE \eqref{sysdn}. Of course, in order to determine $\alpha$ and $\beta$, we need to refine \eqref{hypsys} adding smaller terms, and the process never stops. For that reason, we abandon that procedure, and proceed in a different way in order to determine $\alpha$.\\

Indeed, in order to determine the value of $\alpha$, one may use the formal approach of the physicists in \cite{PSKKPD98}, Section 3.3, page 93, and see that when we have
\[p=3,\;\;\;\alpha=-\frac {\sqrt 2}{16 \sqrt 3}.\]
If $p\not = 3$, we are able to successfully carry on that method and find the value of $\alpha$, and also $\mu$ by \eqref{hypsys1}:

\beqtn
\dsp  \alpha=-\frac{\kappa}{8\sqrt{p(p+1)}},\;\; \mu =\frac{\delta}{8 p}.
\label{hypsys2}
\eeqtn
In fact, this calculation for $p\not =3$ is easier than in \cite{PSKKPD98} since we already know the good scaling from \eqref{hypsys}. However, in order to keep the paper in a reasonable length, we choose to not include that calculation here and we refer the interested reader to the physicists paper \cite{PSKKPD98}.

\medskip

In fact, we should keep in mind that we have a better justification of the value of $\alpha$ in \eqref{hypsys2}: that value is precisely the one that makes our rigorous proof work in the next section! Indeed, we see from Corollary \ref{NprojR*} below that this particular value of $\alpha$, makes some projections small and allows us to conclude.

\medskip

Therefore, since we took $ w(y,s)=e^{i\mu \log s}(v(y,s)+\kappa)$ from \eqref{eqwmu}, it follows from \eqref{formaldecomp}, \eqref{hypsys} and \eqref{hypsys2} that
\beqtn
\begin{array}{lll}
\dsp w(y,s)&=&e^{i\mu \log s}\left[\kappa-\frac{\kappa}{8\sqrt{p(p+1)}}\frac{y^2-2}{\sqrt s}(1+i\delta)+i\frac{\beta}{\sqrt s}+o\left (\frac{1}{\sqrt s}\right)\right ],\\
      &=&e^{i\mu \log s+i\frac{\beta}{\kappa \sqrt s}}\left[\kappa-\frac{\kappa}{8\sqrt{p(p+1)}}\frac{y^2-2}{\sqrt s}(1+i\delta)+o\left (\frac{1}{\sqrt s}\right)\right ],

\end{array}
\label{calfom2}
\eeqtn
in $L^{2}_{\rho}(\R^N)$, and also uniformly on compact sets by standard parabolic regularity.

\medskip

We would like to note that the term $i\beta/(\kappa \sqrt s)$ in the phase is compatible with our rigorous proof. Indeed,
the modulation parameter $\theta(s)$ we introduce in our formulation of the problem below in \eqref{formulav} will be shown to satisfy $|\theta(s)| \le C/\sqrt s$ below in \eqref{estitetas} (note that we may take $\theta_0=0$ in \eqref{estitetas} from the rotation invariance of equation \eqref{GLauto}). Accordingly, we could have chosen to make a modulation technique already in this formal step, but
as one may easily check, there is a one-to-one relation between the modulation parameter and the coordinate $\hat v_0(s)$
and he result and well as the calculations will be completely the same.

\subsection{Outer expansion} From \eqref{calfom2}, we see that the variable
\[z=\frac{y}{s^{1/4}},\]
as given by \eqref{Preslt}, is perhaps the relevant variable for blow-up. Unfortunately,  \eqref{calfom2} provides no shape, since it is valid only on compact sets (note that $z\to 0$ as $s\to \infty$ in this case). In order to see some shape, we may need to go further in space, to the ``outer region``, namely when $z \not = 0$. In view of  \eqref{calfom2}, we may try to find an expression of $w$ of the form
\beqtn
w(y,s)=e^{i\mu \log s+i\frac{\beta}{\kappa \sqrt s}}\left [\varphi_0(z)+(1+i\delta)\frac{\kappa}{4\sqrt{p(p+1)}}\frac{1}{\sqrt s}+o(\frac{1}{\sqrt s})\right]\mbox{  with  }z=\frac{y}{s^{1/4}}.
\label{oe1}
\eeqtn
Plugging this ansatz in equation \eqref{eqv}, keeping only the main order, we end-up with the following equation on $\varphi_0$: 
\beqtn
-\frac 12 z\cdot \nabla \varphi_0(z)-\frac{1+i\delta}{p-1}\varphi_0(z)+(1+i\delta)|\varphi_0(z)|^{p-1}\varphi_0(z)=0,\;\; z=\frac{y}{s^{1/4}}.
\label{oe2}
\eeqtn
Recalling that our aim is to find $w$ a solution of \eqref{eqv} such that $w\to \kappa$ as $s\to\infty$ (in $L^{2}_{\rho}$, hence uniformly on every compact set), we derive from \eqref{oe1} ( with $y=z=0$) the natural condition 
\[ \varphi_0 (0) =\kappa.\]
Therefore, integrating equation \eqref{oe2}, we see that 

\[\varphi_0(z)= \left(p-1+bz^2\right)^{-\frac{1+i\delta}{p-1}},\]
for some $b\in \R$. Recalling also that we want a solution $w\in L^\infty(\R^N)$, we see that $b\geq 0$ and for a nontrivial solution, we should have
\[b>0.\]-
Thus, we have just obtained from \eqref{oe1} that
\beqtn
w(y,s)=e^{i\mu \log s+i\frac{\beta}{\kappa \sqrt s}}\left [\left(p-1+b z^2\right)^{-\frac{1+i\delta}{p-1}}
+(1+i\delta)\frac{\kappa}{4\sqrt{p(p+1)}}\frac{1}{\sqrt s}+o(\frac{1}{\sqrt s})\right]\mbox{  with  }z=\frac{y}{s^{1/4}}.
\label{oe3}
\eeqtn
We should understand this expansion to be valid at least on compact sets in $z$, that is for $|y|<R s^{\frac 14}$, for any $R>0$.

\subsection{Matching asymptotics}
Since \eqref{oe3} holds for $|y|<R s^{\frac 14}$, for any $R>0$, it holds also uniformly on compact sets, leading to the following expansion for $y$ bounded:
\[\dsp w(y,s)=e^{i\mu \log s+i\frac{\beta}{\kappa \sqrt s}}\left [\kappa-(1+i\delta)\frac{\kappa b}{(p-1)^2}\frac {y^2}{s^{\frac12}} +(1+i\delta)\frac{a}{s^{\frac 12}} +o\left(\frac{1}{\sqrt s}\right)   \right ]\]
Comparing with \eqref{calfom2}, we find the following values of $a$ and $b$:
\[b=\frac{(p-1)^2}{8\sqrt{p(p+1)}},\;\;\; a=\frac{\kappa}{4\sqrt{p(p+1)}}.\]
In conclusion, we see that we have just derived the following profile for $w(y,s)$:
\[w(y,s)\sim e^{i\mu \log s+i\frac{\beta}{\kappa \sqrt s}} \varphi(y,s),\]
with 
\beqtn
\begin{array}{l}
\dsp\varphi\left (y,s\right )=\varphi_0\left (\frac{y}{s^{1/4}}\right )+(1+i\delta)\frac{a}{s^{1/2}}\equiv\kappa ^{-i\delta}\left ( p-1+b\frac{|y|^2}{s^{1/2}}\right )^{-\frac{1+i\delta}{p-1}}+(1+i\delta)\frac{a}{s^{1/2}}\\

\displaystyle b=\frac{(p-1)^2}{8\sqrt{p(p+1)}},\;\;\; a=\frac{\kappa}{4\sqrt{p(p+1)}},\;\;\mu= \frac{\delta}{8 p} \mbox{ and } \kappa=(p-1)^{-\frac{1}{p-1}}.
\end{array}
\label{defifi}
\eeqtn

\section{Formulation of the problem}\label{formpb}
We recall that we consider CGL, given by \eqref{GL}, when $\nu=0$, as we mentioned before in Remark \ref{nu0}.\\
The preceding calculation is purely formal. However, the formal expansion provides us with the profile of the function ($w(y,s)=e^{i\mu \log s}\left (\varphi_0(\frac{y}{s^{1/4}})+...\right ))$. Our idea is to linearize equation \eqref{GLauto} around that profile and prove that the linearized equation as well as the nonlinear equation have a solution that goes to zero as $s\to \infty$. Let us introduce $q(y,s)$ and $\theta (s)$ such that

\beqtn
\displaystyle w(y,s)= e^{i\left (\mu \log s+\theta (s)\right )}\left (\varphi(y,s)+q(y,s) \right ),
\label{formulav}
\eeqtn

\[\displaystyle\mbox{where }\varphi\left (y,s\right )=\varphi_0\left (\frac{y}{s^{1/4}}\right )+(1+i\delta)\frac{a}{s^{1/2}}\equiv\kappa ^{-i\delta}\left ( p-1+b\frac{|y|^2}{s^{1/2}}\right )^{-\frac{1+i\delta}{p-1}}+(1+i\delta)\frac{a}{s^{1/2}},\]
\beqtn
\displaystyle b=\frac{(p-1)^2}{8\sqrt{p(p+1)}},\;\;\; a=\frac{\kappa}{4\sqrt{p(p+1)}},\;\;\mu= \frac{\delta}{8 p} \mbox{ and } \kappa=(p-1)^{-\frac{1}{p-1}}.
\label{definitionq}
\eeqtn

In order to guarantee the uniqueness of the couple $(q,\theta)$ an additional constraint is needed, see \eqref{eqmod} below; we will choose $\theta(s)$ such that we kill one the neutral modes of the linearized operator.\\

Note that $\varphi_0(z)=w_0(z)$ has been exhibited in the formal approach and satisfies the following equation 
\beqtn
-\frac 12 z\nabla w_0 -\frac{1+i\delta}{p-1}w_0+(1+i\delta)|w_0|^{p-1}w_0=0,
\label{eqfi0}
\eeqtn
 which makes $\varphi (y,s)$ an approximate solution of  \eqref{GLauto}. If $w$ satisfies equation \eqref{GLauto}, then $q$ satisfies the following equation
\beqtn
\label{eqq1}
\frac{\pa q}{\pa s}=\Lg_0 q-\frac{(1+i\delta)}{p-1} q +L(q,\theta ', y, s)+R^*(\theta ',y,s)
%i\left (\frac \mu s+\theta'(s)\right ) q +V_1 q+ V_2 \bar q+B(q,y,s)+R^*(\theta',y,s),
\eeqtn
where
\beqtn
\begin{array}{lll}
%\tilde\Lg q&=&\Lg_0+(1+i\delta )\Re q \mbox{ where }
\Lg_0 q&=&\Delta q -\frac 12 y\cdot\nabla q,\\
L(q,\theta',y,s)&=&(1+i\delta)\left\{|\varphi+q|^{p-1}(\varphi+q)-|\varphi|^{p-1}\varphi -i\left(\frac \mu s+\theta'(s)\right)q\right\}\\
R^*(\theta',y,s)&=&R(y,s)-i\left (\frac \mu s+\theta '(s)\right )\varphi,\\
R(y,s)&=&-\frac{\pa \varphi}{\pa s}+\Delta \varphi-\frac 12 y\cdot \nabla \varphi-\frac{(1+i\delta)}{p-1}\varphi+(1+i\delta)|\varphi|^{p-1}\varphi
\end{array}
\label{eqqd1}
\eeqtn
 
\medskip

Our aim is to find a $\theta\in C^1([-\log T,\infty),\R$ such that equation \eqref{eqq} has a solution $q(y,s)$ defined for all $(y,s)\in \R^N\times [-\log T,\infty)$ such that
\[\|q(s)\|_{L^\infty}\to 0 \mbox{ as }s\to \infty.\]

From \eqref{eqfi0}, one sees that the variable $z=\frac{y}{s^{1/4}}$ plays a fundamental role. Thus, we will consider the dynamics for $|z|>K$and $|z|<2K$ separately for some $K>0$ to be fixed large.
\subsection{The outer region where $|y|>Ks^{1/4}$}
Let us consider a non-increasing cut-off function $\chi_0\in C^{\infty}(\R^+,[0,1])$ such that $\chi_0(\xi)=1$ for $\xi<1$ and $\chi_0(\xi)=0$ for $\xi>2$ and introduce
\beqtn
\chi(y,s)=\dsp\chi_0\left(\frac{|y|}{Ks^{1/4}}\right),
\label{defchi14}
\eeqtn
where $K$ will be fixed large. Let us define
\beqtn
q_e(y,s)=e^{\frac{i\delta}{p-1}s}q(y,s)\left( 1-\chi(y,s)\right)
\label{defiqe}
\eeqtn
$q_e$ is the part of $q(y,s)$ for $|y|>Ks^{1/4}$. As we will explain in subsection \eqref{outreg}, the linear operator of the equation satisfied by $q_e$ is negative, which makes it easy to control $\|q_e(s)\|_{L^\infty}$. This is not the case for the part of $q(y,s)$ for $|y|<2Ks^{1/4}$, where the linear operator has two positive eigenvalues, a zero eigenvalue in addition to infinitely many negative ones. Therefore, we have to expand $q$ with respect to these eigenvalues in order to control $\|q(s)\|_{L^\infty(|y|<2 K s^{1/4})}$. This requires more work than for $q_e$. The following subsection is dedicated to that purpose. From now on, $K$ will be fixed constant which is chosen such that $\|\varphi(s')\|_{L^\infty(|y|> K {s'}^{1/4})}$ is small enough, namely $\|\varphi_0(z)\|_{L^\infty(|z|>K)}^{p-1}\leq \frac{1}{C(p-1)}$ (see subsection  \eqref{outreg} below, for more details). 

\subsection{The inner region where $|y|< 2K s^{1/4}$}
If we linearize the term $L(q,\theta ', y, s)$ in equation \eqref{eqq1}, then we can write \eqref{eqq1} as

\beqtn
\label{eqq}
\frac{\pa q}{\pa s}=\tilde\Lg q-i\left (\frac \mu s+\theta'(s)\right ) q +V_1 q+ V_2 \bar q+B(q,y,s)+R^*(\theta',y,s),
\eeqtn
where
\beqtn
\begin{array}{lll}
\tilde\Lg q&=&\Delta q -\frac 12 y\cdot\nabla q+(1+i\delta )\Re q,\\
V_1(y,s)&=&(1+i\delta)\frac{p+1}{2}\left(|\varphi|^{p-1}-\frac{1}{p-1}\right),\;\;V_2(y,s)=(1+i\delta)\frac{p-1}{2}\left(|\varphi|^{p-3}\varphi^2-\frac{1}{p-1}\right),\\
B(q,y,s)&=&(1+i\delta)\left ( |\varphi+q|^{p-1}(\varphi+q)-|\varphi|^{p-1}\varphi-|\varphi|^{p-1} q-\frac{p-1}{2}|\varphi|^{p-3}\varphi(\varphi \bar q+\bar\varphi q)\right),\\ 
R^*(\theta',y,s)&=&R(y,s)-i\left (\frac \mu s+\theta '(s)\right )\varphi,\\
R(y,s)&=&-\frac{\pa \varphi}{\pa s}+\Delta \varphi-\frac 12 y\cdot \nabla \varphi-\frac{(1+i\delta)}{p-1}\varphi+(1+i\delta)|\varphi|^{p-1}\varphi
\end{array}
\label{eqqd}
\eeqtn
Note that the term $B(q,y,s)$ is built to be quadratic in the inner region $\dsp |y|\leq K s^{1/4}$. Indeed, we have for all $K\geq 1$ and $s\geq 1$, 
\beqtn
\sup_{|y|\leq 2 K s^{1/4}}|B(q,y,s)|\leq C(K)|q|^2
\label{estiquadinn}
\eeqtn
Note also that $R(y,s)$ measures the defect of $\varphi(y,s)$ from being an exact solution of \eqref{GLauto}. However, since $\varphi(y,s)$ is an approximate solution of \eqref{GLauto}, one easily derives the fact that
\beqtn
\|R(s)\|_{L^\infty}\leq \frac{C}{\sqrt{s}}.
\label{estR*}
\eeqtn
Therefore, if $\theta'(s)$ goes to zero as $s\to \infty$, we expect the term $R^*(\theta',y,s)$ to be small, since \eqref{eqq} and \eqref{estR*} yield
\beqtn
|R^*(\theta',y,s)|\leq \frac{C}{\sqrt{s}}+|\theta'(s)|.
\label{Rexpect}
\eeqtn
Therefore, since we would like to make $q$ go to zero as $s\to \infty$, the dynamics of equation \eqref{eqq} are influenced by the asymptotic limit of its linear term,
\[\tilde \Lg+V_1 q+V_2\bar q,\]
as $s\to\infty$. In the sense of distribution (see the definition of $V_1$ and $V_2$ \eqref{eqq} and  $\varphi$ \eqref{defifi}) this limit is $\tilde \Lg$.

\subsection{Decomposition of $q$}\label{sectidecompq}
For the sake of controlling $q$ in the region $|y|<2K s^{1/4}$, by the spectral properties of $\tilde \Lg$ \eqref{spectLtilde},we will expand the unknown function $q$ with respect to the family $h_n$ and then, with respect to the families 
\beqtn
\tilde h_n=(1+i\delta)h_n\mbox{ and }\hat h_n=i h_n.
\label{defthhh} 
\eeqtn
We start by writing
\beqtn
q(y,s)=\sum_{n\leq M} Q_n(s) h_n(y)+q_{-}(y,s),
\label{decomp1}
\eeqtn
where $h_n$ is the eigenfunctions of $\Lg$ defined in \eqref{eigfun}, $Q_n(s)\in \C$, $q_{-}$ satisfy

\beqtn
Q_n(s)=\dsp \frac{\int q h_n\rho}{\int h_{n}^{2}\rho},
\label{eqQn}
\eeqtn

\[\int q_{-}(y,s)h_n(y)\rho(y)dy=0\mbox{ for all }n\leq M,\]
and $M$ is a fixed even integer satisfying
\beqtn
\dsp M\geq 4\left( \sqrt{1+\delta^2}+1+2 \max_{i=1,2,y\in \R,s\geq 1}|V_i(y,s)|\right),
\label{boundM}
\eeqtn
with $V_{i=1,2}$ defined in \eqref{eqqd}.\\
The function $q_{-}(y,s)$ can be seen as the projection of $q(y,s)$ onto the spectrum of $\tilde \Lg$, which is smaller than $(1-M)/2$. We will call it the infinite dimensional part of $q$ and we will denote it $q_-=P_{-,M}(q)$. We also introduce $P_{+,M}=Id -P_{-,M}$. Notice that $P_{-,M}$ and $P_{+,M}$ are projections. In the sequel, we will denote $P_-=P_{-,M}$ and $P_+=P_{+,M}$.\\
The complementary part $q_+=q-q_-$ will be called the finite dimensional part of $q$. We will expand it as follows
\beqtn
\dsp q_+(y,s)=\sum_{n\leq M}Q_n(s)h_n(y)=\sum_{n\leq M}\tilde q_n(s)\tilde h_n(y)+ \hat q_n(s) \hat h_n(y),
\label{decomp2}
\eeqtn
where $\tilde q_n,\hat q_n\in \R$, $\tilde h_n=(1+i\delta)h_n$ and $\hat h_n=i h_n$. Finally, we notice that for all $s$, we have
\[\dsp\int q_-(y,s)q_+(y,s)\rho(y)dy=0.\] 
Our purpose is to project \eqref{eqq} in order to write an equation for $\tilde q_n$ and $\hat q_n$. Note that
\beqtn
\tilde P_n(q)=\tilde q_n(s)= \Re Q_n(s),\;\hat  P_n(q)=\hat q_n(s)= \Im Q_n(s)-\delta \Re Q_n(s). 
\label{decomp3}
\eeqtn
We conclude from \eqref{decomp1} and \eqref{decomp2}, that
\beqtn
q(y,s)=\left(\sum_{n\leq M}\tilde q_n(s) \tilde h_n(y)+\hat q_{n}(s) \hat h_n(y)\right)+q_-(y,s),
\label{decompq}
\eeqtn
where $\tilde h_n$ and $\hat h_n$ are given by \eqref{defthhh}. 
%\[\tilde h_n=(1+i\delta)h_n,\;\;\hat h_n =i h_n\]
we should keep in mind that this decomposition is unique.
%NEWWWWWWw
\section{Existence}\label{existence}
In this section, we prove the existence of a solution $q(s), \theta(s)$ of problem \eqref{eqq1}-\eqref{eqmod}  such that
\beqtn
\lim_{s\to \infty}\|q\|_{L^\infty}=0,\mbox{ and }|\theta '(s)|\leq C\frac{A^{5}}{s^{3/2}}\mbox{ for all }s\in [-\log T,+\infty).
\label{limitev}
\eeqtn
%where $\theta_0\in\R$.\\
Hereafter, we denote by $C$ a generic positive constant, depending only on $p$ and $K$ introduced in \eqref{defchi14}, itself depending on $p$. In particular, $C$ does not depend on $A$ and $s_0$, the constants that will appear shortly and throughout the paper and need to be adjusted for the proof.\\
We proceed in two subsections. In the first, we give the proof assuming the technicals details. In the second subsection we give the proof of the technicals details.

\subsection{Proof of the existence assuming technical results}\label{pwtr} Our construction is built on a careful choice of the initial data of $q$ at a time $s_0$. We will choose it in the following form:
\begin{definition}[Choice of initial data] Let us define, for $A\geq 1$, $s_0=-\log T>1$ and $d_0,d_1\in \R$,  the function
\beqtn\begin{array}{l}
 \psi_{s_0,d_0,d_1}(y)=\dsp\left [\frac{A}{s_{0}^{3/2}} (1+i\delta)\left(  d_0h_0(y) + d_1 h_1(y) \right)+id_2 \right]\chi (2y,s_0)\\
 \mbox{ where }s_0=-\log T,
 \end{array}
%d_3(s_0,d_0,d_1,d_2)=&-\dsp\frac{A}{s_{0}^{3/2}}\frac{d_0\hat P_0\left ((1+i\delta)\chi (2y,s_0)\right )+d_1\hat P_0\left((1+i\delta)y\chi (2y,s_0) \right)}{\hat P_{0}\left (i\chi (2y,s_0)\right)} 
\label{definitdata1}
\eeqtn
where $h_i$, $i=0,1,2$ are defined by \eqref{eigfun}, $\chi$ is defined by \eqref{defchi14} and $d_2=d_2(d_0,d_1)$ will be fixed later in (i) of Proposition \ref{propinitialdata}. 
\end{definition}
\begin{rem}
Let us recall that we will modulate the parameter $\theta$ to kill one of the neutral modes, see equation \eqref{eqmod} below. It is natural that this condition must be satisfied for the initial data at $s=s_0$. Thus, it is  necessary that we choose $d_2$ to satisfy condition \eqref{eqmod}, see  \eqref{d2} below. 
\end{rem}

The solution of equation \eqref{eqq} will be denoted by $q_{s_0,d_0,d_1}$ or $q$ when there is no ambiguity. We will show that if $A$ is fixed large enough, then, $s_0$ is fixed large enough depending on $A$, we can fix the parameters $(d_0,d_1)\in [-2,2]^2$, so that the solution $v_{s_0,d_0,d_1}\to 0$ as $s\to \infty$ in $L^{\infty}(\R)$, that is \eqref{limitev} holds. Owing to the decomposition given in \eqref{decompq}, it is enough to control the solution in a shrinking set defined as follows
\begin{definition}[A set shrinking to zero]
For all $K>1$, $A\geq 1$ and $s\geq e$, we define $\Vg_A(s)$ as the set of all $r\in L^\infty(\R)$ such that
\beqtn
\begin{array}{ll}
\|r_e\|_{L^\infty(\R)}\leq \frac{A^{M+2}}{s^{\frac 14}}    ,                               &\|\frac{r_-(y)}{1+|y|^{M+1}}\|_{L^\infty(\R)}\leq  \frac{A^{M+1}}{s^{\frac{M+2}{4}}}, \\
|\hat r_j|,\;|\tilde r_j|\leq \frac{A^j}{s^{\frac{j+1}{4}}}\mbox{ for all }3\leq j\leq M,  & |\tilde r_0|,\;|\tilde r_1|\leq \frac{A}{s^{\frac 32}},\\
 |\hat r_0| \leq \frac{1}{s^{\frac 32}},\;\;|\hat r_1| \leq \frac{A^4}{s^{\frac 32}},&|\tilde r_2|\leq \frac{A^5}{s}, \;\;   |\hat r_2|\leq \frac{A^3}{s},                                                                        
\end{array}
\label{thess}
\eeqtn
\label{defthess}
\end{definition}
Since $A\geq 1$, the  sets $\Vg_A(s)$ are increasing (for fixed $s$) with respect to $A$ in the sense of inclusions.\\
We also show the following property of elements of $\Vg_A(s)$:\\
For all $A\geq 1$, there exists $s_{01}(A)\geq 1$, such that for all $s\geq s_{01}$ and $r\in\Vg(A)$, we have
\beqtn
\|r\|_{L^\infty}(\R)\leq C(K)\frac{A^{M+2}}{s^{\frac 14}}
\label{estilinfty}
\eeqtn
where $C$ is a positive constant (see Claim \ref{propshrinset} below for the proof).\\
By \eqref{estilinfty}, if a solution $q$ stays in $\Vg(A)$ for $s\geq s_{01}$, then it converges to $0$ in $L^\infty(\R)$.

\medskip
%Correction Rq 9 p3

\noindent So far, the phase $\theta(s)$ introduced in \eqref{formulav} is arbitrary, in fact as we will show below in Proposition \ref{enopropmod}. We can use a modulation technique to choose $\theta(s)$ in such a way that we impose the condition 

\beqtn
\hat P_{0} (q(s))=0,
\label{eqmod}
\eeqtn
which allows us to kill the neutral direction of the operator $\tilde \Lg$ defined in \eqref{eqq}.
%Correction Rq 9 p3
 %\medskip
Reasonably, our aim is then reduced to the following proposition:

\begin{prop}[Existence of a solution trapped in $\Vg_A(s)$] There exists $A_2\geq 1$ such that for $A\geq A_2$ there exists $s_{02}(A)$ such that for all $s_0\geq s_{02}(A)$, there exists $(d_0,d_1)$ such that if $q$ is the solution of \eqref{eqq}-\eqref{eqmod}, with initial data given by \eqref{definitdata1} and \eqref{d2}, then $v\in \Vg_A(s)$, for all $s\geq s_0$.
\label{propsop1}
\end{prop}
This proposition gives the stronger convergence to $0$ in $L^\infty(\R)$ thanks to \eqref{estilinfty}.\\
Let us first be sure that we can choose the initial data such that it starts in $\Vg_A(s_0)$. In other words, we will define a set where where will be selected the good parameters $(d_0,d_1)$ that will give the conclusion of Proposition \ref{propsop1}. More precisely, we have the following result:

\begin{prop}[Properties of initial data] For each $A\geq 1$, there exists $s_{03}(A)>1$ such that for all $s_0\geq s_{03}$:\\
(i) $\hat P_{0}\left (i\chi (2y,s_0)\right)\not = 0$ and the  parameter $d_2(s_0,d_0,d_1)$ given by 
\beqtn\label{d2}
d_2(s_0,d_0,d_1)=-\dsp\frac{A}{s_{0}^{3/2}}\frac{d_0\hat P_0\left ((1+i\delta)\chi (2y,s_0)\right )+d_1\hat P_0\left((1+i\delta)y\chi (2y,s_0) \right)}{\hat P_{0}\left (i\chi (2y,s_0)\right)} 
\eeqtn

is well defined, where $\chi$ defined in \eqref{defchi14}.\\(ii) If $\psi$ is given by \eqref{definitdata1} and \eqref{d2} with $d_2$ defined by  \eqref{d2}.Then, there exists a quadrilateral $\Dg_{s_0}\subset [-2,2]^2$ such that the mapping $(d_0,d_1)\to (\tilde \psi_0,\tilde \psi_1)$ (where $\psi$ stands for $\psi_{s_0,d_0,d_1}$) is linear, one to one from $\Dg_{s_0}$ onto $[-\frac{A}{s_{0}^{3/2}},\frac{A}{s_{0}^{3/2}}]^2$. Moreover it is of degree $1$ on the boundary.\\
(iii) For all $(d_0,d_1)\in \Dg_{s_0}$, $\psi_e\equiv 0$, $\hat \psi_0=0$, $ |\tilde\psi_i|+|\hat\psi_i|\leq C A e^{-\gamma s_0}$ for some $\gamma>0$, for some $\gamma>0$ and for all $3\leq i\leq M$ and $1\leq j \leq M$. Moreover , $\|\frac{\psi_-(y)}{(1+|y|)^{M+1}}\|_{L^\infty(\R)}\leq C\frac{A}{s_{0}^{\frac M4+1}}$.\\
(iv) For all $(d_0,d_1)\in \Dg_{s_0}$, $\psi_{s_0,d_0,d_1}\in \Vg_A(s_0)$ with strict inequalities except for $(\tilde \psi_0,\tilde\psi_1)$.

\label{propinitialdata}
%\lable{propsop2}
\end{prop}
The proof of previous proposition is postponed to subsection \ref{tecnicsect}. 

\medskip

In the following, we find a local in time solution for equation \eqref{eqq} coupled with the condition \eqref{eqmod}.

\begin{prop}(\textbf{Local in time solution and modulation for problem \eqref{eqq}-\eqref{eqmod} with initial data  \eqref{definitdata1}-\eqref{d2}})\label{enopropmod}
For all $A\geq 1$, there exists $T_3(A)\in (0,1/e)$ such that for all $T\leq T_3$, the following holds:\\
For all $(d_0,d_1)\in D_T$, there exists $s_{max}>s_0=-\log T$ such that problem \eqref{eqq}-\eqref{eqmod} with initial data at $s=s_0$,
\[(q(s_0),\theta(s_0))=(\psi_{s_0,d_0,d_1},0),\]
where $\psi_{s_0,d_0,d_1}$ is given by \eqref{definitdata1} and \eqref{d2}, has a unique solution $q(s),\theta(s)$ satisfying $q(s)\in V_{A+1}(s)$ for all $s\in[s_0,s_{max})$.
\end{prop}
The proof of this proposition will be given later in page \pageref{prolocalintime}.

\medskip

Let us now give the proof of Proposition \ref{propsop1}.\\
\textit{Proof of Proposition \ref{propsop1}}: Let us consider $A\geq 1$, $s_0\geq s_{03}$, $(d_0,d_1)\in \Dg_{s_0}$, where $s_{03}$ is given by Proposition \ref{propinitialdata}. From the existence theory (which follows from the Cauchy problem for equation \eqref{GL}), starting in $\Vg_A(s_0)$ which is in $\Vg_{A+1}(s_0)$, the solution stays in $\Vg_A(s)$ until some maximal time $s_*=s_*(d_0,d_1)$. Then, either:\\
$\bullet$ $s_*(d_0,d_1)=\infty$ for some $(d_0,d_1)\in \Dg_{s_0}$, then the proof is complete.\\
$\bullet$ $s_*(d_0,d_1)<\infty$, for any  $(d_0,d_1)\in \Dg_{s_0}$, then we argue by contradiction. By continuity and the definition of $s_*$, the solution on $s_*$ is in the boundary of $\Vg_A(s_*)$. Then, by definition of $\Vg_A(s_*)$, one at least of the inequalities in that definition is an equality. Owing to the following proposition, this can happen only for the first two components $\tilde q_0,\tilde q_1$. Precisely we have the following result

\begin{prop}[Control of $q(s)$ by $(q_0(s),q_1(s))$ in $\Vg_A(s)$]. There exists $A_4\geq 1$ such that for each $A\geq A_4$, there exists $s_{04}\in \R$ such that for all $s_0\geq s_{04}$. The following holds:\\
If $q$ is a solution of \eqref{eqq} with initial data at $s=s_0$ given by \eqref{definitdata1} and \eqref{d2} with $(d_0,d_1)\in \Dg_{s_0}$, and $q(s)\in \Vg(A)(s)$ for all $s\in [s_0,s_1]$, with $q(s_1)\in \pa \Vg_A(s_1)$ for some $s_1\geq s_0$, then:\\
(i)\textbf{(Smallness of the modulation parameter $\theta$ defined in \eqref{eqwmu})} For all $s\in [s_0,s_1]$, 
\[|\theta '(s)|\leq \frac{CA^5}{s^{3/2}}.\]
(ii) \textbf{(Reduction to a finite dimensional problem)} We have: 
\[(\tilde q_0(s_1),\tilde q_1(s_1))\in \pa\left( \left [-\frac{A}{s_{1}^{\frac 32}},\frac{A}{s_{1}^{\frac 32}} \right ]^2\right).\]
(iii)\textbf{(Transverse crossing)} There exists  $m\in \{0,1\}$ and $\omega\in \{-1,1\}$ such that
\[\omega \tilde q_m (s_1)=\frac{A}{s_{1}^{\frac 32}}\mbox{ and }\omega \frac{d\tilde q_m}{ds} (s_1)>0.\]
\label{propcontrol}
\end{prop}
\label{indextheory}
Assume the result of the previous proposition, for which the proof is given below in page \pageref{proofiniti}, and continue the proof of Proposition \ref{propsop1}. Let $A\geq A_4$ and $s_0\geq s_{04}(A)$. It follows from Proposition \ref{propcontrol}, part (ii) that $\left (q_0(s_*),q_1(s_*)\right)\in\pa\left( \left [-\frac{A}{s_{*}^{\frac 32}},\frac{A}{s_{*}^{\frac 32}} \right ]^2\right)$, and the following function
\[
\begin{array}{ll}
\phi&:\Dg_{s_0}\to \pa ([-1,1]^2)\\
&(d_0,d_1)\to \frac{s_{*}^{3/2}}{A}\left( q_0, q_1\right)_{(d_0,d_1)}(s_*)\mbox{, with }s_*=s_*(d_0,d_1),
\end{array}
\]
is well defined. Then, it follows from Proposition \ref{propcontrol}, part (iii) that $\phi$ is continuous. On the other hand, using Proposition \ref{propinitialdata} (ii)-(iv)
together with the fact that $q(s_0)=\psi_{s_0,d_0,d_1}$, we see that when $(d_0,d_1)$ is in the boundary of the rectangle $\Dg_{s_0}$, we have strict inequalities for the other components.\\
Applying the transverse crossing property given by (iii) of Proposition \ref{propcontrol}, we see that $q(s)$ leaves $\Vg_A(s)$ at $s=s_0$, hence $s_*(d_0,d_1)=s_0$. Using Proposition \ref{propinitialdata}, part (ii), we see that the restriction of $\phi$ to the boundary is of degree 1. A contradiction, then follows from the index theory. Thus there exists a value $(d_0,d_1)\in \Dg_{s_0}$ such that for all $s\geq s_{0}$, $q_{s_0,d_0,d_1}(s)\in \Vg_A(s)$. This concludes the proof of Proposition \ref{propsop1}.\\
Using (i) of Proposition \ref{propcontrol}, we get the bound on $\theta'(s)$. This concludes the proof of \eqref{limitev}.

\subsection{Proof of the technical results result}\label{tecnicsect}
This section is devoted to the proof of the existence result given by Theorem \ref{thm1}. We proceed in 4 steps, each of them making a separate 
subsection.
\begin{itemize}
\item In the first subsection, we give some properties of the shrinking set $\Vg_A(s)$ defined by \eqref{thess} and translate our goal of making $q(s)$ go to $0$ in $L^\infty(\R)$ in terms of belonging to $\Vg_A(s)$. We also give the proof of Proposition \ref{propinitialdata}.
%exhibit a two parameter initial data family for equation \eqref{eqqd1} whose coordinates are very small (with respect to the requirement of $\Vg_A(s)$), except the two first $\tilde q_0$ and $\tilde q_1$.

\item In second subsection, we solve the local in time Cauchy problem for equation \eqref{eqq} coupled with some orthogonality condition.

\item In the third subsection using the spectral properties of equation \eqref{eqq}, we reduce our goal from the control of $q(s)$ (an infinite dimensional variable) in $\Vg_A(s)$ to control its two first components ($\tilde q_0$,$\tilde q_1$) a two variables in $[-\frac{A}{s^{\frac 32}},\frac{A}{s^{\frac32} }]^2$.

\item  In the fourth subsection, we solve the finite dimensional problem using the index theory and conclude the proof of Theorem \ref{thm1}
.

\end{itemize}

\bigskip

\subsubsection{Properties of the shrinking set $\Vg_A(s)$ and preparation of initial data}
In this subsection, we give some properties of the shrinking set defined by \eqref{thess}. Let us first introduce the following claim:
\begin{cl}[Properties of the shrinking set defined by  \eqref{thess}]
 %For all $K>1$, $A\geq 1$ and $s\geq e$, we define $\Vg_A(s)$ as the set of all $r\in L^\infty(\R)$ such that
%\beqtn
%\begin{array}{ll}
%\|r_e\|_{L^\infty(\R)}\leq \frac{A^{M+2}}{s^{\frac 14}}    ,                               &\|\frac{r_-(y)}{1+|y|^{M+1}}\|_{L^\infty(\R)}\leq  \frac{A^{M+1}}{s^{\frac{M+2}{4}}}, \\
%|\hat r_j|,\;|\tilde r_j|\leq \frac{A^j}{s^{\frac{j+1}{4}}}\mbox{ for all }3\leq j\leq M,  & |\tilde r_0|,\;|\tilde r_1|\leq \frac{A}{s^{\frac 32}},\\
%|\tilde r_2|\leq \frac{A^5}{s},                                                                         & |\hat r_1| \leq \frac{A^4}{s^{\frac 32}},\\
%|\hat r_2|\leq \frac{A^5}{s^{\frac 32}},                                                                         & |\hat r_0| \leq \frac{1}{s^{\frac 32}},
%\end{array}
%\label{thess}
%\eeqtn
For all $r\in \Vg_A(s)$,

\medskip

\noindent (i) $\|r\|_{L^\infty(|y|<2K s^{\frac 14})}\leq C(K) \frac{A^{M+1}}{s^{\frac 14}}$ and $\|r\|_{L^\infty}(\R)\leq C(K)\frac{A^{M+2}}{s^{\frac 14}}$. \\
\noindent (ii) for all $y\in \R$, $|r(y)|\leq C \frac{A^{M+1}}{s}(1+|y|^{M+1})$. 
\label{propshrinset}
\end{cl}
\textit{Proof: }Take $r\in \Vg_A(s)$ and $y\in \R$.\\
(i) If $|y|\geq 2 K s^{\frac 14}$, then we have from the definition of  $r_e$  \eqref{defiqe}, $|r(y)|=|r_e(y)|\leq \frac{A^{M+2}}{s^{\frac 14}} $.\\
Now, if $|y| < 2 K s^{\frac 14}$, since we have for all $0\leq j\leq M$, $|\tilde r_j|+|\hat r_j|\leq C\frac{A^j}{s^{\frac{j+1}{4}}}$ from \eqref{thess} (use the fact that $M\geq 4$), we write from \eqref{decompq}
\beqtn
\begin{array}{ll}
|r(y)|&\leq \dsp\left(\sum_{j\leq M}    |\tilde h_j||\tilde h_j|+ |\hat h_j||\hat h_j|\right )+|r_-(y)|,\\
&\leq C\dsp\sum_{j\leq M}  \frac{A^{M+1}}{s^{\frac{j+1}{4}}}(1+|y|)^j+ \frac{A^{M+1}}{s^{\frac{M+2}{4}}}(1+|y|)^{M+1},\\
& \leq C\dsp\sum_{j\leq M}  \frac{A^{M+1}}{s^{\frac{j+1}{4}}}(1+K s^{\frac 14})^j+ \frac{A^{M+1}}{s^{\frac{M+2}{4}}}(1+K s^{\frac 14})^{M+1} \leq C \frac{(KA)^{M+1}}{s^{\frac 14}},
\end{array}
\label{estipss}
\eeqtn
which gives (i).\\
(ii) Just use \eqref{estipss} together with the fact that for all $0\leq j\leq M$, $|\tilde r_j|+|\hat r_j|\leq C\frac{A^{M+1}}{s}$ from  \eqref{thess}. This ends the proof of Claim \ref{propshrinset}. $\blacksquare$

\medskip

Let us now give the proof of Proposition \ref{propinitialdata}.\\
%\textbf{Remark: }In some sense, $\psi_{s_0,d_1,d_2,d_3}$ is reduced to its components on $\tilde h_{i,i=0,1,2}$.\\
\textit{Proof of Proposition \ref{propinitialdata}} For simplicity, we write $\psi$ instead of  $\psi_{s_0,d_0,d_1}$. We note that, from Claim \ref{propshrinset}, (iv) follows from (ii) and (iii) by taking $s_0=-\log T$ large enough (that is $T$ is small enough). Thus, we only prove  (i), (ii) and (iii). Consider $K\geq 1$, $A\geq 1$ and $T\leq 1/e$. Note that $s_0=-\log T\geq 1$.\\
The proof of (i) is a direct consequence of (iii) of the following claim
 \begin{cl}
There exists $\gamma=\frac{1}{64}>0$ and $T_2<1/e$ such that for all $K\geq 1$ and $T\leq T_2$, if $g$ is given by $(1+i\delta)\chi(2y,s_0)$, $(1+i\delta)y\chi(2y,s_0)$, $(1+i\delta)h_2(y)\chi(2y,s_0)$ or $i\chi(2y,s_0)$, then $\left\|\frac{g_-(y)}{1+|y|^{M+1}}\right\|_{L^\infty}\leq \frac{C}{s_{0}^{\frac M4}}$ and all $\hat g_i$, $\tilde g_i$ for $0\leq i\leq M$ are less than $Ce^{-\gamma s_0}$. expect:\\
i) $ |\tilde g_0-1|\leq C e^{-\gamma s_0}$ when $g=(1+i\delta)\chi(2y,s_0)$.\\
ii)  $|\tilde g_1-1|\leq C e^{-\gamma s_0}$ when $g=(1+i\delta)y\chi(2y,s_0)$.\\
%iii)  $|\tilde g_2-1|\leq C e^{-\gamma s_0}$ when $g=(1+i\delta)h_2(y)\chi(2y,s_0)$.\\
iii)  $|\hat g_0-1|\leq C e^{-\gamma s_0}$ when $g=i\chi(2y,s_0)$.
\label{cl}
\end{cl}
\textit{Proof: } In all cases, we write 
\beqtn
g(y)=p(y)+r(y)\mbox{ where }p(y)=(1+i\delta)\mbox{ or }(1+i\delta)y \mbox{ or }i \mbox{ and }r(y)=p(y)(\chi(2y,s_0)-1).
\label{pcl1}
\eeqtn
From the uniqueness of the decomposition \eqref{decompq}, we see that $p_-\equiv 0$ and al $\hat p_i$, $\tilde p_i$ are zero except\\
$\tilde p_0=1$ (when $p(y)=(1+i\delta)$), $\tilde p_1=1$ (when $p(y)=(1+i\delta)y$) and $\hat p_0=1$ (when $p(y)=i$).\\
Concerning the cases $2|y|<K s^{\frac 14}$ and $2|y|>K s^{\frac 14}$, we have the definition of $\chi$  \eqref{defchi14},
\[1-\chi(2y,s)\leq \left(\frac{2|y|}{K s_{0}^{\frac 14}}\right)^{M-1},\]
\[\rho(y)(1-\chi(2y,s))\leq \sqrt{\rho(y)}\sqrt{\rho\left(\frac K2s^{\frac 14}\right)}\leq Ce^{-\frac{K^2 s_0}{64}}\sqrt{\rho(y)}.\]
Therefore, from \eqref{decompq} and \eqref{pcl1}, we see that
\beqtn
\begin{array}{l}
|r(y)|\leq C(1+|y|^2)\left(\frac{2|y|}{Ks_{0}^{\frac 14}}\right)^{M-1}\leq C\frac{(1+|y|^{M+1})}{s_{0}^{\frac M4}},\\
|\hat r_j|+|\tilde r_j|\leq Ce^{-\frac{K^2\sqrt {s_0}}{64}}\mbox{ for all }j\leq M.
\end{array}
\label{pcl2}
\eeqtn
Hence, using  \eqref{pcl2} and \eqref{decomp1} and the fact that $|r_j(y)|\leq C (1+|y|)^M$, for all $j\leq M$, we get also 
\[|r_-(y)|\leq C\frac{(1+|y|)^M}{s_{0}^{\frac M4}}.\]
Using \eqref{pcl1} and the estimates for $p(y)$ stated below, we conclude the proof of Claim \ref{cl} and (i) of Proposition \ref{propinitialdata}. \\
(ii) of Proposition \ref{propinitialdata}: From \eqref{definitdata1} and \eqref{d2}, we see that
\beqtn
\left (
\begin{array}{l}
\tilde \psi_0\\
\tilde \psi_1
%\tilde \psi_2
\end{array}\right )
=G
\left (\begin{array}{l}
d_0\\d_1
\end{array}
\right )
\mbox{ where }G=(g_{i,j})_{0\leq i,j\leq 1}.
\eeqtn
Using Claim \ref{cl}, we see from \eqref{definitdata1} and \eqref{d2} that
\beqtn
|d_2|\leq C(|d_0|+|d_1|)e^{-\gamma s_0}
\label{initiad3}
\eeqtn
for $T$ small enough. Using again Claim \ref{cl}. We see that $\frac{s_{0}^{\frac 32}}{A} G\to Id$ and

as $s_0\to\infty$ (for fixed $K$ and $A$), which concludes the proof of (ii) of Proposition \ref{propinitialdata}.\\
(iii) of Proposition \ref{propinitialdata}: Since $supp (\psi)\subset B(0,K s_{0}^{\frac 14})$ by  \eqref{definitdata1} and \eqref{d2}, we see that $\psi_e\equiv 0$ and that
\[
\begin{array}{ll}
\hat \psi_0 =\hat P_{0}(\psi)&=\frac{A}{s_{0}^{\frac 32}}\left(d_0 \hat P_{0}((1+i\delta)\chi(2y,s_0))+d_1 \hat P_{0}((1+i\delta)y\chi(2y,s_0))\right)\\
&+d_2 \hat P_{0}(i\chi(2y,s_0)) ,
\end{array}
\]
which is zero from the definition of $d_2$  \eqref{definitdata1} and \eqref{d2}. Using the fact that $|d_{i,i=0,1}|\leq 2$ and the bound on $d_2$ by \eqref{initiad3}, we see that the estimates on $\hat \psi_j$ and $\tilde \psi_j$ and $\psi_-$ in (iii) follows from \eqref{definitdata1} and \eqref{d2} and Claim \ref{cl}. This concludes the proof of Proposition \ref{propinitialdata}. $\blacksquare$

\medskip

In the following we give the proof of Local in time solution for problem \eqref{eqq}-\eqref{eqmod}. In fact, we impose some orthogonality condition given by \eqref{eqmod}, killing the one of the zero eigenfunctions of the linearized operator of equation \eqref{eqq}.

\medskip

\textit{Proof of Proposition \ref{enopropmod}:}\label{prolocalintime} From solution of the local in time Cauchy problem for equation \eqref{GL} in $L^\infty(\R)$, there exists $s_1>s_0$ such that equation \eqref{GLauto} with initial data (at $s=s_0$) $\varphi(y,s_0)+ \psi_{s_0,d_0d_1}(y)$, where $\varphi(y,s)$ is given by \eqref{defifi} has a unique solution $w(s)\in C([s_0,s_1),L^\infty(\R))$. Now, we have to find a unique $(q(s), \theta(s))$ such that
\beqtn
w(y,s)=e^{i(\mu\log s+\theta(s))}\left(\varphi(y,s)+q(y,s)\right)
\label{promod1}
\eeqtn
and \eqref{eqmod} is satisfied. Using \eqref{decomp3}, we can write  \eqref{eqmod} as follows
\[\hat P_0 (q)=\Im \left(\int q(y,s)\rho(y)dy\right )-\delta \Re \left(\int q(y,s)\rho(y)dy\right )=\Im \left((1-i\delta)\int q(y,s)\rho(y)dy\right)=0,\]
or using \eqref{promod1}
\[F(s,\theta)\equiv \Im\left((1-i\delta)\int \left(e^{-i(\mu\log s+\theta(s))}  w(y,s)-\varphi(y,s)\right)\rho(y)dy \right)=0.\]
Note that
\[ \frac{\pa F}{\pa\theta} (s,\theta)=-\Re\left((1-i\delta)\int e^{-i(\mu\log s+\theta(s))}  w(y,s) \rho(y)dy\right).\]
From (iii) in Proposition \ref{propinitialdata}, $F(s_0,0)=P_{0,M}(\psi_{s_0,d_0,d_1})=0$ and 
\[\frac{\pa F}{\pa\theta} (s_0,0)=-\Re \left((1-i\delta)\int (\varphi(y,s_0)+\psi_{s_0,d_0d_1}(y)  )\rho(y)dy \right)=-\kappa +O\left(\frac{1}{s_{0}^{1/4}}\right )\mbox{ as }s_0\to \infty,\]
for fixed $K$ and $A$.\\
Therefore, if $T$ is small enough in terms of $A$, then $\frac{\pa F}{\pa \theta}(s_0,0)\not =0$, and from the implicit function Theorem, there exists $s_2\in(s_0,s_1)$ and $\theta\in C^1([s_0,s_2),\R)$ such that $F(s,\theta(s))=0$ for all $s\in [s_0,s_2)$.Defining $q(s)$ by \eqref{promod1} gives a unique solution of the problem  \eqref{eqqd}-\eqref{eqmod} for all $s\in[s_0,s_2)$. Now, since we have from (iv) of Proposition \ref{propinitialdata}, $q(s_0)\in V_{A}(s_0)  \begin{array}{l}\subset\\ \not =\end{array} V_{A+1}(s_0)  $, there exists $s_3\in (s_0,s_2)$ such that for all $s\in [s_0,s_3)$, $q(s)\in V_{A+1}(s)$. This concludes the proof of Proposition \ref{enopropmod}. $\blacksquare$

\subsubsection{Reduction to a finite dimensional problem}
In the following we give the proof of Proposition  \ref{propcontrol}:\\
%\textit{Proof of Proposition  \ref{propcontrol}}:
\label{proofiniti}\\
The idea of the proof is to project equation  \eqref{eqq} on the different components of the decomposition \eqref{decompq}. More precisely, we claim that Proposition \ref{propcontrol}  is a consequence of the following

\begin{prop} There exists $A_5\geq 1$ such that for all $A\geq A_5$, there exists $s_5(A)$ such that the following holds for all $s_0\geq s_5$:\\
Assuming that for all $s\in [\tau,s_1]$ for some $s_1\geq \tau\geq s_0$, $q(s)\in \Vg_A(s)$ and $\hat q_0(s)=0$, then  the following holds for all $s\in [\tau, s_1]$:\\

\noindent (i) (\textbf{Smallness of the modulation parameter}):
\[|\theta '(s)|\leq C\frac{A^5}{s^{\frac 32}}\]
(ii) (\textbf{ODE satisfied by the expanding mode}): For $m=0$ and $1$, we have
\[|\tilde q_{m}^{'}-\left(1-\frac m2\right)\tilde q_m|\leq \frac{C}{s^{\frac 32}}.\]
(iii) (\textbf{ODE satisfied by the null mode}):
\[|\tilde q_2' +\frac 2s \tilde{q}_2 |\leq \frac{C A^3}{s^2}.\]
 (iv) (\textbf{Control of null and negative modes}):

\[|\hat q_1(s)|\leq e^{-\frac{(s-\tau)}{2}}|\hat q_1(\tau)|+\frac{CA^5}{s^{\frac 32}},\]
\[|\hat q_2(s)|\leq e^{-(s-\tau)}|\hat q_2(\tau)|+\frac{C}{s},\]
\[|\hat q_j(s)|\leq e^{-j\frac{(s-\tau)}{2}}|\hat q_j(\tau)|+\frac{CA^{j-1}}{s^{\frac{j+1}{4}}},\mbox{ for all }3\leq j\leq M,\]
\[|\tilde q_j(s)|\leq e^{-(j-2)\frac{(s-\tau)}{2}}|\tilde q_j(\tau)|+\frac{CA^{j-1}}{s^{\frac{j+1}{4}}},\mbox{ for all }3\leq j\leq M,\]
\[ \left\|\frac{q_-(y,s)}{1+|y|^{M+1}} \right\|_{L^\infty} \leq e^{-\frac{M+1}{4}(s-\tau)}\left\|\frac{q_-(\tau)}{1+|y|^{M+1}}  \right\|_{L^\infty} +C\frac{A^M}{s^{\frac{M+2}{4}}}, \]
\[\|q_e(y,s)\|_{L^\infty}\leq e^{-\frac{(s-\tau)}{2(p-1)}}\|q_e(\tau)\|_{L^\infty}  +\frac{C A^{M+1}}{\tau^{\frac 14}}(1+s-\tau)\]

\medskip
\label{propode}
\end{prop}
The idea of the proof of Proposition \ref{propode} is to project equations \eqref{eqq1} and \eqref{eqq} according to the decomposition \eqref{decompq}. However because of the number of parameters and coordinates in \eqref{decompq}, the computation become too long. That is why Subsection \ref{proofpropode} is devoted to the proof of Proposition \ref{propode}. 

\begin{rem} The coefficient in front of $\frac{\tilde q_2}{s}$ in (iii) of Proposition \ref{propode} is $`2`$. In our proof, see page \pageref{2generique} below, that coefficient is the sum of four contributions, which depend on $p$ in a non trivial way. Thus, it may appear miraculous to see the sum of such contributions equal to $`2`$. The same phenomena occur in the subcritical range of parameters, see \cite{MZ07} and also the heat equation, with a critical gradient term (see \cite{TZ17}). In fact, adopting the approach of Pierre Rapha{\"e}l and co-authors, one may see that the coefficient $`2`$ appears in a natural way due to scaling considerations (From a personal communication of Pierre Rapha{\"e}l).
\end{rem}

\medskip

Let us now derive Proposition \ref{propcontrol} from Proposition \ref{propode}.

\medskip

\textit{Proof of Proposition \ref{propcontrol} assuming Proposition \ref{propode}:}\\
We will take $A_4\geq A_5$. Hence, we can use the conclusion of Proposition \ref{propode}.\\
(i) The proof follows from (i) of Proposition \ref{propode}. Indeed by choosing $T_4$ small enough, we can make $s_0=-\log T$ bigger than $s_5(A)$.\\

(ii) We notice that from Claim \ref{propshrinset} and the fact that $\hat q_0(s)=0$, it is enough to prove that for all $s\in [s_0,s_1]$,

\beqtn\label{goalnullmode}
|\tilde q_2(s)|<\frac{A^5}{s}.
\eeqtn

\beqtn
\begin{array}{ll}

\|q_e\|_{L^\infty(\R)}\leq \frac{A^{M+2}}{2s^{\frac 14}}, &\left\|\frac{q_-(y,s)}{1+|y|^{M+1}} \right\|_{L^\infty} \leq \frac{A^{M+1}}{2 s^{\frac{M+2}{4}}},\\
|\hat q_j|, |\tilde q_j|\leq \frac{A^j}{2s^{\frac{j+1}{4}}}\mbox{ for all }3\leq j\leq M,& |\hat q_1|\leq \frac{A^4}{2s^{\frac 32}},\;\;\;\;\;\; |\hat q_2|\leq \frac{A^3}{2s}.
\end{array}
\label{object}
\eeqtn
Let us first prove \eqref{goalnullmode}. Arguing by contradiction, we assume that 
\[\tilde q_2(s_*)=\omega \frac{A^5}{s_*}\mbox{ and for all }s\in [s_0,s^*[,\;\; |\tilde q_2(s)|<\frac{A^5}{s}.\]
Of course, we can reduce to the case $\omega=1$. Note by (iv) of Proposition \ref{propinitialdata} that $|\tilde q_2(s_0)|<\frac{A^5}{s_0}$, hence $s_*>s_0$, and the interval $[s_0,s_*]$ is not empty.\\
 
 By minimality, it follows that $\tilde q_{2}'(s)\geq \frac{\pa}{\pa s}\left( \frac{A^5}{s}\right)_{|s=s_*}$,
\beqtn\label{pq1}
\tilde q_2'(s)\geq -\frac{A^5}{s_{*}^{2}},
\eeqtn
in the one hand, recalling, from (iii) of Proposition \ref{propode},  that
\[|\tilde q_2' +\frac 2s \tilde{q}_2 |\leq \frac{CA^3}{s^2}.\]
We write 
\beqtn\label{pq2}
\tilde q_2'(s_*)\leq-2\frac{\tilde q_2(s_*)}{s_*}+\frac{CA^3}{s_{*}^{2}}=\frac{-2A^5+CA^3}{s_{*}^{2}}, 
\eeqtn
for $A$ large enough a contradiction follows from \eqref{pq1} and \eqref{pq2}. Thus \eqref{goalnullmode} holds.

\medskip

Now, let us deal with \eqref{object}. Define $\sigma =\log A$ and take $s_0\geq \sigma$ (that is $T\leq e^{-\sigma}=1/A$) so that for all $\tau \geq s_0$ and $s\in [\tau, \tau +\sigma]$, we have 
\beqtn
\tau \leq s\leq \tau+\sigma\leq \tau +s_0\leq 2\tau\mbox{ hence }\frac{1}{2\tau} \leq \frac 1s\leq \frac 1\tau\leq \frac{2}{s}. 
\label{taus}
\eeqtn
We consider two cases in the proof.\\
\textbf{ Case 1:  $s\leq s_0+\sigma$.}\\
Note that \eqref{taus} holds with $\tau=s_0$. Using (iv) of Proposition \ref{propode} and estimate (iii) of Proposition \ref{propinitialdata} on the initial data $q(.,s_0)$ (where we use \eqref{taus} with $\tau =s_0$), we write
\beqtn
\begin{array}{l}
|\hat q_1(s)|\leq C A e^{-\gamma \frac s2}+ \frac{CA^5}{s^{3/2}},\\
|\hat q_2(s)|\leq  C A e^{-\gamma \frac s2}+ \frac{C}{s},\\
|\tilde q_j(s)|\leq CA e^{-\gamma \frac s2}+ \frac{C A^{j-1}}{s^{\frac{j+1}{4}}}\mbox{ for all }3\leq j\leq M,\\
|\hat q_j(s)|\leq CA e^{-\gamma \frac s2}+ \frac{C A^{j-1}}{s^{\frac{j+1}{4}}}\mbox{ for all }3\leq j\leq M,\\
\left\|\frac{q_-(s)}{1+|y|^{M+1}}\right\|_{L^\infty}\leq C\frac{A}{\left(\frac s2\right)^{\frac M4+2}}+C\frac{A^M}{s^{\frac{M+2}{4}}},\\
\|q_e(s)\|_{L^\infty}\leq \frac{CA^{M+1}}{\left(\frac s2\right)^{\frac 14}}(1+\log A).
\end{array}
\eeqtn
Thus, if $A\geq A_6$ and $s_0\geq s_6(A)$ (that is $T\leq e^{-s_6(A)}$) for some positive $A_6$ and $s_6(A)$, we see that \eqref{object} holds.\\
\textbf{ Case 2:  $s> s_0+\sigma$.}\\
Let $\tau=s-\sigma>s_0$. Applying (iv) of Proposition \ref{propode} and using the fact that $q(\tau)\in \Vg_A(\tau)$, we write (we use \eqref{taus} to bound any function of $\tau$ by a function of $s$)
\beqtn
\begin{array}{l}
|\hat q_1(s)|\leq  e^{-\frac\sigma 2}\frac{A^6}{\left(\frac s2\right)^{3/2}}+ \frac{CA^5}{s^{3/2}},\\
|\hat q_2(s)|\leq e^{-\sigma}\frac{A^2}{\left(\frac s2\right)}+ \frac{C}{s},\\
|\tilde q_j(s)|\leq e^{-\frac{(j-2)\sigma}{2}}\frac{A^j}{\left(\frac s2\right )^{\frac{j+1}{4}}}+ \frac{C A^{j-1}}{s^{\frac{j+1}{4}}}\mbox{ for all }3\leq j\leq M,\\
|\hat q_j(s)|\leq e^{-\frac{j\sigma}{2}}\frac{A^j}{\left(\frac s2\right )^{\frac{j+1}{4}}}+ \frac{C A^{j-1}}{s^{\frac{j+1}{4}}}\mbox{ for all }3\leq j\leq M,\\
\left\|\frac{q_-(s)}{1+|y|^{M+1}}\right\|_{L^\infty}\leq e^{-\frac{M+1}{4}\sigma}\frac{A^{M+1}}{\left (\frac s 2\right)^{\frac{M+2}{4}}}+C\frac{A^M}{s^{\frac{M+2}{4}}},\\
\|q_e(s)\|_{L^\infty}\leq e^{-\frac{\sigma}{2(p-1)}}\frac{A^{M+2}}{\left (\frac s2\right)^{\frac 14}}+ \frac{CA^{M+1}}{\left(\frac s2\right)^{\frac 14}}(1+\sigma).
\end{array}
\eeqtn
For all the coordinates, it is clear that if $A\geq A_7$ and $s_0\geq s_7(A)$ for some positive $A_7$ and $s_7(A)$, then \eqref{goalnullmode} and \eqref{object} is satisfied (remember that $\sigma=\log A$).

\medskip

Conclusion of (ii): If $A\geq \max (A_6,A_7,A_8)$ and $s_0\geq \max (s_6(A), s_7(A),s_8(A))$, then \eqref{object} is satisfied. Since we know that $q(s_1)\in \pa V_A(s_1)$, we see from the definition of $V_A(s)$ that $(\tilde q_0(s_1),\tilde q_1(s_1)) \in \pa [-\frac{A}{s_{1}^{3/2}},\frac{A}{s_{1}^{3/2}}]^2$. This concludes the proof of (ii) of Proposition \ref{propcontrol}.

\medskip

(iii) From (ii), there is $m=0$, $1$  and $\omega=\pm 1$ such that $\tilde q_m(s_1)=\omega \frac{A}{s_{1}^{3/2}}$. \\

Using (ii) of Proposition \ref{propcontrol}, we see that for $m=0$ or $1$
\[\omega \tilde q^{'}_{m}(s_1)\geq (1-\frac m2)\omega \tilde q_m(s_1)-\frac{C}{s_{1}^{3/2}}.\]
Taking $A$ large enough gives $\omega \tilde q_{m}^{'}(s_1)>0$, for $m=0,1$ and concludes the proof of Proposition \ref{propcontrol}. $\blacksquare$

\subsection{Proof of Proposition \ref{propode}}\label{proofpropode}
In this section, we prove Proposition \ref{propode}. We just have to project equations \eqref{eqq1} and \eqref{eqq} to get equations satisfied by the different coordinates of the decomposition \eqref{decompq}.  We proceed as Section 5 in \cite{MZ07}, taking into account the new scaling law $\frac{y}{s^{1/4}}$. We note that the projections of $V_1 q+V_2\bar q$, $B$ and $R^*$ in  \eqref{eqqd}, will need much more effort and this is due to the fact that we are dealing with the critical case.

\medskip

More precisely, the proof will be carried out in 3 subsections
\begin{itemize}
\item In the first subsection, we deal with equation \eqref{eqq} to write equations satisfied by $\tilde q_j$ and $\hat q_j$. Then, we prove (i), (ii), (iii) and (iv) (expect the two last identities) of Proposition \ref{propode}.
\item In the second subsection, we first derive from equation \eqref{eqq} an equation satisfied by $q_-$ and prove the last but one identity in (iv) of Proposition \ref{propode}.
\item In the third subsection, we project equation \eqref{eqq1} (which is simpler than \eqref{eqq}) to write an equation satisfied by $q_e$ and prove the last identity in (iv) of Proposition \ref{propode}.
\end{itemize}
\subsubsection{The finite dimensional part $q_+$}
We proceed in 2 parts:
\begin{itemize}
\item In Part 1, we project equation \eqref{eqq} to get equations satisfied by $\tilde q_j$ and $\hat q_j$. 
\item In Part 2, we prove (i), (ii) and (iii) of Proposition \ref{propode}, together with the estimates concerning $\tilde q_j$ and $\hat q_j$ in (iv).

\end{itemize}
\textbf{Part 1: The projection of equation \eqref{eqq} on the eigenfunctions of the operator $\tilde \Lg$ }
In the following, we will find the main contribution in the projections $\tilde P_{n,M}$ and $\hat P_{n,M}$ of the six terms appearing in equation \eqref{eqq}: $\pa_s q$, $\tilde \Lg q$, $-i(\frac \mu s+\theta'(s))q$, $V_1 q+V_2 \bar q$, $B(q,y,s)$ and $R^*(\theta',y,s)$. Most of the time, we give two estimates of error terms, depending on whether we use or not the fact that $q(s)\in  \Vg_A(s)$.\\
\textbf{First term: $\dsp\frac{\pa q}{\pa_s}$.}\\ 
From \eqref{decomp3}, its projection on $\tilde h_n$ and $\hat h_n$ is $\tilde q_{n}'$ and $\hat q_{n}'$ respectively:
\beqtn
\tilde P_{n}\left(\frac{\pa q}{\pa s}\right)=\tilde q_n'\mbox{ and }\hat P_{n}\left(\frac{\pa q}{\pa s}\right)=\hat q_n'.
\eeqtn
\textbf{Second term: $\tilde \Lg q$.}\\ 
We can easily see that
\beqtn
\begin{array}{ll}
\tilde P_{n}\left(\tilde \Lg q\right)&=\dsp (1-\frac n 2)\tilde q_n,\\
\hat P_{n}\left(\tilde \Lg q\right)&=\dsp -\frac n2 \hat q_n.
\end{array}
\eeqtn
\textbf{Third term: $-i(\frac \mu s+\theta ')q$.} 
It is enough to project $i q$, from \eqref{decomp3}, we have
\beqtn
\begin{array}{ll}
\tilde P_{n}\left(-i(\frac \mu s+\theta ')q\right)&=-\dsp\left(\frac \mu s+\theta'(s)\right)(-\hat q_n-\delta \tilde q_n),\\
\hat P_{n}\left(-i(\frac \mu s+\theta ')q\right)&=-\dsp\left(\frac \mu s+\theta'(s)\right)(\delta \hat q_n+(1+\delta^2)\tilde q_n).
\end{array}
\eeqtn
If in addition $q(s)\in \Vg_A(s)$, then the error estimates can be bounded from Definition \ref{defthess} as follows:
\begin{corollary} For all $A\geq 1$, there exists $s_{10}(A)\geq 1$ such that for all $s \geq s_{10}(A) $, if $q\in \Vg_A(s)$ and $|\theta'(s)|\leq \frac{C A^5}{s^{3/2}}$, then:\\
a) for all $1\leq n\leq M$, we have 
\[\left| \hat P_{n}\left(-i(\frac \mu s+\theta ')q\right)  \right|\leq C\frac{A^n}{s^{\frac{n+5}{4}}}.,\]
b) for $1\leq n\leq M$, we have
\[\left| \tilde P_{n}\left(-i(\frac \mu s+\theta ')q\right)  \right|\leq C\frac{A^n}{s^{\frac{n+5}{4}}},\]
c) for $n=0$, $\left| \tilde P_{0}\left(-i(\frac \mu s+\theta ')q\right) \right|  +\left|  \hat P_{0}\left(-i(\frac \mu s+\theta ')q\right) \right|\leq \frac{C}{s^{3/2}}$
\end{corollary}
\textbf{Fourth term: $V_1 q+V_2\bar q$.}

We claim the following
\begin{lemma}[Projection of $V_1 q$ and $V_2 \bar q$ ] 
(i) It holds that
\beqtn
|V_i(y,s)|\leq C\frac{(1+|y|^2)}{s^{1/2}}, \mbox{ for all $y\in \R$ and $s\geq 1$,}
\label{bdVi1}
\eeqtn
and for all $k\in\N^*$
\beqtn
V_i(y,s)=\dsp \sum_{j=1}^{k}\frac{1}{s^{j/2}}W_{i,j}(y)+\tilde W_{i,k}(y,s),
\label{dcVi}
\eeqtn
where $W_{i,j}$ is an even polynomial of degree $2j$ and $\tilde W_{i,k}(y,s)$ satifies
\beqtn
\mbox{for all $s\geq 1$ and $|y|\leq s^{1/4}$, } \left| \tilde W_{i,k}(y,s)\right|\leq  C\frac{(1+|y|^{2k+2})}{s^{\frac{k+1}{2}}}.
\eeqtn
(ii) The projection of $V_1 q$ and $V_2 \bar q$ on $(1+i\delta)h_n$ and $i h_n$, and we have
\beqtn
\begin{array}{l}
\dsp|\tilde P_n(V_1 q)|+|\hat P_n(V_1 q)|\\
\leq \dsp\frac{C}{s^{1/2}} \dsp \sum_{j=n-2}^{M}(|\tilde q_j|+|\hat q_j|)
+\sum_{j=0}^{n-3}\frac{C}{s^{\frac{n-j}{4}}}(|\tilde q_j|+|\hat q_j|)
+\frac{C}{s^{1/2}}\left\| \frac{q_-}{1+|y|^{M+1}}\right\|_{L^\infty},

\end{array}
\label{bdVi}
\eeqtn
and the same holds for $V_2\bar q$

\label{lembdVi}
\end{lemma}
\begin{rem} If $n\leq 2$, the first sum in \eqref{bdVi} runs for $j=0$ to $M$ and the second sum doesn't exist.
\end{rem}
If in addition $q(s)\in \Vg_A(s)$, then the error estimates can be bounded from Definition \ref{defthess} as follows:
\begin{corollary} For all $A\geq 1$, there exists $s_{11}(A)\geq 1$ such that for all $s \geq s_{11}(A) $, if $q\in \Vg_A(s)$, then:\\
a) for $3\leq n\leq M$, we have 
\[\dsp |\tilde P_n(V_1 q)|+|\hat P_n(V_1 q)|\leq C\frac{A^{n-2}}{s^{\frac{n+1}{4}}},\]
b) for $n=0,1$ or $2$, we have
\[\dsp |\tilde P_n(V_1 q)|+|\hat P_n(V_1 q)|\leq \frac{CA^5}{s^{\frac{3}{2}}},\]
\label{cbdVi}
\end{corollary}

\textit{Proof of Lemma \ref{lembdVi}: }\\
(i)  The estimates of $V_1q$ and $V_2 \bar q$ are the same, so we only deal with $V_1q$. Let $F(u)=\frac{(p+1)}{2}(1+i\delta)\left [|u|^{p-1}-\frac{1}{p-1}\right]$, where $u\in\C$ and consider $z=\frac{y}{s^{1/4}}$. Note that from \eqref{eqq} and \eqref{eqfi0}, we have 
\[V_1(y,s)=F(\varphi(y,s))\mbox{, where } \varphi(y,s)=\varphi_0(\frac{y}{s^{1/4}})+\frac{a}{s^{1/2}}(1+i\delta).\]
Note that there exist positive constant $c_0$ and $s_0$ such that $\varphi_0(z)|$ and $|\varphi(y,s)|=|\varphi_0(\frac{y}{s^{1/4}})+\frac{a}{s^{1/2}}(1+i\delta)|$
are both larger than $\frac{1}{c_0}$ and smaller than $c_0$, uniformly in $|z|<1$ and for $s\geq s_0$. Since $F(u)$ is $C^\infty$ for $\frac{1}{c_0}\leq |u|\leq c_0$, we expand it around $u=\varphi_0(z)$ as follows: for all $s\geq s_0$ and $|z|<1$,
\[
\begin{array}{lll}
\left|F\left(\varphi_0(z)+\frac{a}{s^{1/2}}(1+i\delta)\right)-F\left(\varphi_0(z)\right)\right|&\leq &\dsp\frac{C}{s^{1/2}},\\
\left|F\left(\varphi_0(z)+\frac{a}{s^{1/2}}(1+i\delta)\right)-F\left(\varphi_0(z)\right)-\dsp \sum_{j=1}^{n}\frac{1}{s^{j/2}}F_j(\varphi_0(z))\right|&\leq &\dsp\frac{C}{s^{\frac{n+1}{2}}},
\end{array}
\]
where $F_j(u)$ are $C^\infty$. Hence, we can expand $F(u)$ and $F_j(u)$ around $u=\varphi_0(0)$ and write for all $s\geq s_0$ and $|z|<1$,

\[
\begin{array}{l}
\left|F\left(\varphi_0(z)+\frac{a}{s^{1/2}}(1+i\delta)\right)-F\left(\varphi_0(0)\right)\right|\leq C z^2+\frac{C}{s^{1/2}},\\
\left|F\left(\varphi_0(z)+\frac{a}{s^{1/2}}(1+i\delta)\right)-F\left(\varphi_0(0)\right)-\dsp \sum_{l=1}^{n} c_{0,l}z^{2l}-\sum_{j=1}^{n}\sum_{l=0}^{n-j}\frac{c_{j,l}}{s^{j/2}} z^{2l}\right|\\
\leq C |z|^{2n+2}+\dsp\sum_{j=1}^{n}\frac{C}{j^{1/2}}|z| ^{2(n-j)+2}           \frac{C}{s^{\frac{n+1}{2}}}.
\end{array}
\]
Since $F(\varphi_0(0))=F(\kappa)=0$ and $z=\frac{y}{s^{1/4}}$, this gives us estimates in (i),  when $s\geq s_0$ and $|y|<s^{1/4}$. Since $V_1$ is bounded, the inequalities still valid when $|y|\geq s^{1/4}$ and then when $s\geq 1$.

\medskip

(ii) Note first that it is enough to prove the bound \eqref{decomp2} for the projection of $V_i q$ onto $h_n$ to get the same bound for $\tilde P_n(V_i q)$ and $\hat P_n(V_i q)$. Since in addition, the proof for $V_2 \bar q$ is the same as for $V_1 q$, we only prove \eqref{bdVi} for the projection of $V_1 q$ onto $h_n$. Using \eqref{decompq} and the fact that $\tilde h_n=(1+i\delta )h_n$ and $\hat h_n=ih_n$, we see that the projection is given by
\beqtn
\begin{array}{ll}
\dsp \int h_n V_1 q\rho&=\dsp\int h_n V_1 q_-\rho+\sum_{j=0}^{M}\tilde q_j \int h_n\tilde h_j V_1\rho+\sum_{j=0}^{M}\hat q_j \int h_n\hat h_j V_1\rho.\\
&\dsp=\int h_n V_1 q_-\rho+(1+i\delta)\sum_{j=0}^{M}\tilde q_j \int h_n h_j V_1\rho+i\sum_{j=0}^{M}\hat q_j \int h_n h_jV_1\rho
\label{bd1}
\end{array}
\eeqtn
The first term can be bounded by
\beqtn
\dsp\int h_n V_1 \left(\frac{1+|y|^2}{s^{1/2}}\right)|q_-|\rho\leq  \frac{C}{s^{1/2}}\left\| \frac{q_-}{1+|y|^{M+1}}\right\|_{L^\infty}.
\label{bd2}
\eeqtn
Now we deal with the second term. We only focus on the terms involving $h_j$ .\\
If $j\geq n-2$, we use \eqref{bdVi1} to write $|\int h_n h_j V_1\rho|\leq \frac{C}{s^{1/2}}$.\\
If $j\leq n-3$, then we claim that
\beqtn
\dsp\left |\int h_n h_j V_1\rho\right|\leq \frac{C}{s^{\frac{n-j}{4}}},
\label{bd3}
\eeqtn
(this actually vanishes if $j$ and $n$ have different parities). It is clear that \eqref{bdVi} follows from \eqref{bd1}, eqref{bd2} and \eqref{bd3}.\\
Let us prove \eqref{bd3}. Note that $k\equiv\left [ \frac{n-j-1}{2}\right ]$ (which is in $\N^*$ since $j\leq n-3$) is the largest integer such that $j+2k <n$. We use \eqref{dcVi} to write
\beqtn
\begin{array}{l}
\dsp \int h_nh_j V_1 \rho=\dsp\int_{|y|< s^{1/4}}h_nh_j V_1 \rho+\int_{|y|>s^{1/4}}h_nh_j V_1 \rho,\\
= \dsp\sum_{l=1}^{k}\frac{1}{s^{l/2}}\int_{|y|< s^{1/4}}h_n h_jW_{1,l}\rho +O\left(\frac{1}{s^{\frac{ \left [ \frac{n-j-1}{2}\right ]+1  }{2}}}\int(1+|y|^{n-j+1})|h_n||h_j|\rho dy\right)+\int_{|y|>s^{1/4}}h_nh_j V_1 \rho,\\
= \dsp\sum_{l=1}^{k}\frac{1}{s^{l/2}}\int_{\R^N}h_n h_jW_{1,l}\rho +O\left(\frac{1}{s^{\frac{ \left [ \frac{n-j-1}{2}\right ]+1  }{2}}}\right)
-\sum_{l=1}^{k}\frac{1}{s^l}\int_{|y|>s^{1/4}} h_n h_j W_{1,l}\rho +\int_{|y|>s^{1/4}}h_nh_j V_1 \rho,
\end{array}
\eeqtn
since $deg(h_j W_{1,l})=j+2l\leq j+2k<n=deg(h_n)$, $h_n$ is orthogonal to $h_j W_{1,l}$ and
\[\int_{\R^N}h_n h_j W_{1,l}\rho=0.\]
Since $|\rho(y)|\leq \dsp C e^{-cs^{1/2}}$ when $|y|>s^{1/4}$, the integrals over the domain ${|y|>s^{1/4}}$ can be bounded by
\[\dsp Ce^{-cs^{1/2}} \int_{|y|>s^{1/4}} |h_n||h_j|(1+|y|^{2k})\sqrt{\rho} \leq Ce^{-cs^{1/2}}. \]
Using that $\left [ \frac{n-j-1}{2}\right ]+1\geq \frac{n-j}{2}$, we deduce that \eqref{bd3} holds. Hence, we have proved \eqref{bdVi} and this concludes the proof of Lemma \ref{lembdVi}. $\blacksquare$.\\

We need further refinements when $n=0,2$ for the terms $\tilde P_{2,M}(V_1 q)$, $\tilde P_{2,M}(V_2\bar q)$, $\hat P_{0,M}(V_1 q)$ and $\hat P_{0,M}(V_2 \bar q)$. More precisely

\begin{lemma}{Projection of $V_1 q$ and $V_2\bar q$ on $\hat h_0$ and $\tilde h_2$}

(i) It holds that for $i=1,2$
\beqtn
\forall s\geq 1 \mbox{ and } |y|< s^{1/4},\;\; \left| V_i(y,s)-\frac{1}{s^{1/2}} W_{i,1}(y)-\frac{1}{s} W_{i,2}(y) \right|\leq \frac{C}{s^{3/2}}(1+|y|^6),
\label{estiVii}
\eeqtn
where
\beqtn
\begin{array}{lll}
W_{1,1}&=&-(1+i\delta)\frac{b(p+1)}{2(p-1)^2}h_2(y)\\
W_{1,2}&=&(1+i\delta)\frac{b^2(p+1)}{2(p-1)^3}h_{2}^{2}(y),\\
W_{2,1}&=&-(1+i\delta)\frac{p-1}{2}\frac{b}{(p-1)^3}\left(p-1+2i\delta\right)h_2(y)\\
W_{2,2}&=&(1+i\delta)\frac{b^2}{2(p-1)^3}\Big[(p^2-4p+1)h_{2}^{2}(y)+i\delta\left (8(p-2)(1-y^2)+y^43(p-1)\right)\Big].
\end{array}
\label{defWii}
\eeqtn

(ii) The projection of $V_1 q$ and $V_2 \bar q$ on $\tilde h_2$ satisfy
\beqtn
\begin{array}{l}
\dsp\Big | \tilde P_2(V_1 q)+\frac{1}{s^{1/2}} \left [\sum_{j\geq 0} \tilde A_{j}^{1} \tilde q_j+  \sum_{j\geq 0} \hat B_{j}^{1} \hat q_j   \right ] +\frac{\tilde q_2 }{s} 60b^2(p+1)\frac{1}{(p-1)^2}\\

\dsp+ \tilde P_2(V_2 \bar q)+\frac{1}{s^{1/2}} \left [\sum_{j\geq 0} \tilde A_{j}^{2} \tilde q_j+  \sum_{j\geq 0} \hat B_{j}^{2} \hat q_j   \right ]   -\frac{\tilde q_2 }{s} 60b^2(p+1)\frac{p^2-4p+1}{2(p-1)^3}\Big |\\
\leq \dsp\frac{C}{s}\sum_{j=0,j\not =2,}^{M} |\tilde q_j| + \frac{C}{s} \sum_{j=0}^{M} |\hat q_j|+\frac{C}{s} \left\|\frac{q_-}{1+|y|^M}\right \|_{L^\infty}
+\frac{C}{s^{3/2}}|\tilde q_2(s)|.
%+\frac{C}{s}\left\| \frac{q}{1+|y|^M}\right\|_{L^\infty}.
\end{array}
\label{bdVtilde2}
\eeqtn
where, for all $j\geq 0$ 

\[
\begin{array}{ll}
\tilde A_{j}^{1}=\tilde P_2(W_{1,1}\tilde h_j)  &     \hat B_{j}^{1}=\tilde P_2(W_{1,1}\hat h_j)\\
\tilde A_{j}^{2}=\tilde P_2(W_{2,1}\bar{\tilde h}_j)    &    \hat B_{j}^{2}=\tilde P_2(W_{2,1}\bar{\hat h}_j),
\end{array}
\]

and
\[
\tilde A_{j}^{1}=-\tilde A_{j}^{2},  \;\;\; \hat B_{j}^{1}=-\hat B_{j}^{2} \]

ii) The projection of $V_1 q$ and $V_2 \bar q$ on $\hat h_0$ satisfy
\beqtn
\begin{array}{l}
\left| \hat P_0(V_1 q)+\frac{\tilde q_2}{s^{1/2}} 4\delta b(p+1)\frac{(p+1)}{(p-1)^2}\right |\\
+\left| \hat P_0(V_2 \bar q))-\frac{\tilde q_2}{s^{1/2}} 4\delta b(p+1)\frac{(p-3)}{(p-1)^2} \right |\\
\leq\dsp\frac{C}{s^{1/2}}\sum_{j=0,j\not =2}^{M} |\tilde q_j| + \frac{C}{s^{1/2}}\sum_{j=0}^{M} |\hat q_j|+\frac{C}{s^{1/2}} \left\|\frac{q_-}{1+|y|^M}\right \|_{L^\infty}
+\frac{C}{s^{3/2}}|\tilde q_2(s)|.
%+\frac{C}{s}\left\| \frac{q}{1+|y|^M}\right\|_{L^\infty}.
\end{array}
\label{bdVhat0}
\eeqtn
\label{projPtls}
\end{lemma}
In addition , if $q(s)\in \Vg_A(s)$, then the error estimates can be bounded from \eqref{thess} as follows;
\begin{corollary} For all $A\geq 1$, there exists $s_{12}(A)\geq 1$ such that for all $s\geq s_{12}(A)$, if $q(s)\in \Vg_A(s)$, then 

\[\left| \tilde P_2(V_1 q) +\frac{\tilde q_2 }{s} 60b^2(p+1)\frac{1}{(p-1)^2}+ \tilde P_2(V_2 \bar q)-\frac{\tilde q_2 }{s} 60b^2(p+1)\frac{p^2-4p+1}{2(p-1)^3}\right |\leq C\frac{A^3}{s^2}
\]
\[\left| \hat P_0(V_1 q)+\frac{\tilde q_2}{s^{1/2}} 4\delta b(p+1)\frac{(p+1)}{(p-1)^2}\right |\leq C\frac{A^3}{s^{3/2}}\]
\[\left| \hat P_0(V_2 \bar q))-\frac{\tilde q_2}{s^{1/2}} 4\delta b(p+1)\frac{(p-3)}{(p-1)^2} \right |\leq C\frac{A^3}{s^{3/2}}\]
\label{corprojPtls}
\end{corollary}
\textit{Proof of Lemma \ref{projPtls}: }(i) This is a simple, but lengthy computation that we omit. For more details see Appendix \ref{detailsV1V2}.\\
(ii) Using \eqref{estiVii} and \eqref{decompq}, we see that 
\beqtn
\begin{array}{ll}
V_1 q&\dsp =\frac{1}{s^{1/2}}W_{1,1}q+\frac{1}{s} W_{1,2} q+O\left (\frac{q(1+|y|^6)}{s^{3/2}}\right )\\
&\dsp=\frac{1}{s^{1/2}} W_{1,1}\left(\sum_{j=0}^{M} \tilde q_j\tilde h_j +\sum_{j=0}^{M} \hat q_j\hat h_j \right)\\
&\dsp+\frac{1}{s} W_{1,2}\left(\sum_{j=0}^{M} \tilde q_j\tilde h_j +\sum_{j=0}^{M} \hat q_j\hat h_j \right)+O\left(\frac{q(1+|y|^6)}{s^{3/2}}\right),
\end{array}
\label{projV1}
\eeqtn
where $O$ is uniform with respect to $|y|< s^{1/4}$. 
By the definition of the components of $q$:
\beqtn
\dsp |q(y,s)|\leq \left [ \sum_{i=1}^{M}(|\tilde q_i|+|\hat q_i|)+\frac{\|q(s)\|_{L^\infty}}{1+|y|^M}\right ](1+|y|^M),
\label{compoq}
\eeqtn
and we may replace the occurrence of $q$ by its components.

\medskip

When projecting \eqref{projV1} on $\tilde h_2$ (use \eqref{decomp3} for the definition of that projection), we write using the definition \eqref{defWii} of $W_{1,1}$ and $W_{1,2}$ and 
\eqref{compoq}
\beqtn
\begin{array}{l}
\dsp\Big |  \tilde P_2(V_1 q) -\frac{1}{s^{1/2}} \left [\sum_{j\geq 0} \tilde P_2(W_{1,1}\tilde h_j) \tilde q_j+  \sum_{j\geq 0} \tilde P_2(W_{1,1}\hat h_j) \hat q_j   \right ]  -\frac{\tilde q_2(s)}{s} \frac{b^2(p+1)}{2(p-1)^3}\tilde P_2((\tilde h_2)^2 h_2)\\

\dsp+ \tilde P_2(V_2 \bar q)-\frac{1}{s^{1/2}} \left [\sum_{j\geq 0} \tilde P_2(W_{2,1}\tilde h_j) \tilde q_j+  \sum_{j\geq 0} \tilde P_2(W_{2,1}\hat h_j) \hat q_j   \right ] 

-\frac{\tilde q_2(s)}{s} \frac{b^2(p+1)(p^2-4p+1)}{2(p-1)^3}\tilde P_2 (h_{2}^{3})\Big |\\
\leq 
\dsp\frac{C}{s^{1/2}}\sum_{j=0,\not =2}^{M} |\tilde q_j| + \frac{C}{s^{1/2}} \sum_{j=0}^{M}|\hat q_j|+\frac{C}{s^{1/2}} \left\|\frac{q_-}{1+|y|^M}\right \|_{L^\infty}+ \frac{C}{s^{3/2}}|\tilde q_2(s)|

\end{array}
\eeqtn
We note by \eqref{defWii}

\[ \begin{array}{lll}

W_{1,1}\tilde h_j&=&\left [ (p-1)(1+i\delta)-i\delta(p+1)\right ] \frac{b (p+1)}{2(p-1)^2}h_j h_2 \\
 W_{2,1}\tilde h_j&=&\left [ -(p-1)(1+i\delta)+i\delta (p-3)\right ] \frac{b (p+1)}{2(p-1)^2}h_j h_2 \\
 W_{1,1}\hat h_j&=&\left [ \delta(1+i\delta)-i(p+1)\right ] \frac{b (p+1)}{2(p-1)^2}h_j h_2 \\
 W_{2,1}\hat h_j&=&\left [ -\delta(1+i\delta)+i(p-1)\right ] \frac{b (p+1)}{2(p-1)^2}h_j h_2 
\end{array}
 \]
and
we deduce that for all $k\geq 0$
\beqtn
\begin{array}{lll}
\tilde P_k(W_{1,1}\tilde h_j)  &=&-  \tilde P_k(W_{2,1}\bar{\tilde h}_j) \\
 \tilde P_k(W_{1,1}\hat h_j)  & =&-   \tilde P_k(W_{2,1}\bar{\hat h}_j),
\end{array}
\label{order12}
\eeqtn

In particular when $k=2$, terms involving $\frac{\tilde q_j}{s^{1/2}} $ and $\frac{\hat q_j}{s^{1/2}} $ disappear, it happens that its occurrence coming from $\tilde P_2(V_1 q)$ and $\tilde P_2(V_2\bar q)$ do cancels.\\

Therefore, the problem is reduced to the projection of  $(\tilde h_2)^2 h_2$  and $h_{2}^{3}$ on $\tilde h_2$

\[\tilde P_2((\tilde h_2)^2h_2)=120 (1-p),\;\; \tilde P_2(\tilde h_{2}^{3})=120.\]

The other bounds on\eqref{bdVtilde2} and \eqref{bdVhat0} are similar, thus we skip it.$\blacksquare$
\begin{rem}\label{bettrestiVi}
From equation \eqref{order12}, we can see that  terms involving $\frac{\tilde q_j}{s^{1/2}} $ and $\frac{\hat q_j}{s^{1/2}} $ disappear on $\tilde P_k(V_1q)+\tilde P_k(V_2\bar q)$. Then, arguing as in Lemma \ref{lembdVi} and Corollary \ref{cbdVi} we can obtain a better estimate on $|\tilde P_k(V_1q)+\tilde P_k(V_2\bar q)|$, for all $k\geq 0$:

\beqtn
|\tilde P_k(V_1q)+\tilde P_k(V_2\bar q)|\leq \frac{C}{s^{3/2}},
\label{NewbdVi}
\eeqtn
\end{rem}
\medskip

\textbf{Fifth term: $B(q,y,s)$}
Let us recall from \eqref{eqqd} that:
\[B(q,y,s)=(1+i\delta)\left ( |\varphi+q|^{p-1}(\varphi+q)-|\varphi|^{p-1}\varphi-|\varphi|^{p-1} q-\frac{p-1}{2}|\varphi|^{p-3}\varphi(\varphi \bar q+\bar\varphi q)\right).\]
We have the following
\begin{lemma}\label{lemfifB} The function $B=B(q,y,s)$ can be decomposed for all $s\geq 1$ and $|q|\leq 1$ as
\beqtn
\dsp \sup_{|y|<s^{1/4}}
\left|B-\sum_{l=0}^{M}  \sum_{\begin{array}{l}  0\leq j,k\leq M+1\\2\leq j+k\leq M+1   \end{array}  }\frac{1}{s^{l/2}} \left [   B_{j,k}^{l}(\frac{y}{s^{1/4}}) q^j\bar q^k+\tilde B_{j,k}^{l}(y,s) q^j \bar q^k \right ]\right|
\leq C|q|^{M+2}+\frac{C}{s^{\frac{M+1}{2}}}
\eeqtn
where $ B_{j,k}^{l}(\frac{y}{s^{1/4}})$ is an even polynomial of degree less or equal to $M$ and the rest $\tilde  B_{j,k}^{l}(y,s)$ satisfies 
\[\forall s\geq 1 \mbox{ and } |y|< s^{1/4},\; \left|   \tilde  B_{j,k}^{l}(y,s)\right |\leq C\frac{1+|y|^{M+1}}{s^{\frac{M+1}{2}}}.\]
Moreover,
\[\forall s\geq 1  \mbox{ and } |y|< s^{1/4},\; \left|   B_{j,k}^{l}(\frac{y}{s^{1/4}})+\tilde  B_{j,k}^{l}(y,s)\right |\leq C.\]
On the other hand, in the region $|y|\geq s^{1/4}$, we have
\beqtn
|B(q,y,s)|\leq C|q|^{\bar p},
\eeqtn
for some constant $C$ where $\bar p=\min (p,2)$.
\end{lemma}
\textit{Proof:} See the proof of Lemma 5.9, page 1646 in \cite{MZ07}. $\blacksquare$

\begin{lemma}[The quadratic term $B(q,y,s)$] For all $A\geq 1$, there exists $s_{13}$ such that for all $s\geq s_{13}$, if $q(s)\in \Vg_A(s)$, then:\\
a) the projection of $B(q,y,s)$ on $(1+i\delta) h_n$ and on $i h_n$, for $n \geq 3$ satisfies
\beqtn
|\tilde P_n (B(q,y,s))|+|\hat P_n (B(q,y,s)) |\leq C\frac{A^n}{s^{\frac{n+2}{4}}}.
\label{lembin1}
\eeqtn
b) For $n=0,1$, we have
\beqtn
|\tilde P_n (B(q,y,s))|+|\hat P_n (B(q,y,s)) |\leq \frac{C}{s^{\frac{5}{2}}},
\label{lembin2}
\eeqtn
and 
\beqtn
|\tilde P_2 (B(q,y,s))|+|\hat P_2 (B(q,y,s)) |\leq \frac{C}{s^{\frac 52}},
\label{lembin3}
\eeqtn
\label{lembin}
\end{lemma}
\textit{Proof :} We only prove estimate \eqref{lembin1}, since \eqref{lembin2} can be proved in the same way.  The estimation \eqref{lembin3}, is also proved in the same way and using the Taylor expansion of $B$, established in Appendix \ref{TEB}.

\medskip

It is enough to prove estimate \eqref{lembin1} for the projection on $h_n$ since it implies the same estimate on $\tilde P_n$ and $\hat P_n$ through \eqref{decomp3}. We have
\[ \dsp\int h_n B(q,y,s)\rho dy =\int_{|y|<s^{1/4}}h_n B(q,y,s)\rho dy +\int_{|y|>s^{1/4}}h_n B(q,y,s)\rho dy.\]
Using Lemma \ref{lemfifB}, we deduce that
\[\begin{array}{l}
\left| 
%\Big
\dsp\int_{|y|<s^{1/4}}h_n B(q,y,s)\rho dy-\int_{|y|<s^{1/4}}h_n \rho\sum_{l=0}^{M}  \sum_{\begin{array}{l}  0\leq j,k\leq M+1\\2\leq j+k\leq M+1   \end{array}  }\frac{1}{s^{l/2}} \left [   B_{j,k}^{l}(\frac{y}{s^{1/4}}) q^j\bar q^k+\tilde B_{j,k}^{l}(y,s) q^j \bar q^k \right ]
% \Big |
\right|\\
\leq C\dsp\int_{|y|<s^{1/4}}|h_n|| \rho |(|q|^{M+2}+\frac{1}{s^{\frac{M+1}{2}}}).
\end{array}
\]
Let us write
\[\dsp B_{j,k}^{l}(\frac{y}{s^{1/4}}) =\sum_{i=0}^{M/2} b_{j,k}^{l,i}(\frac{y}{s^{1/4}}) \left(\frac{y}{s^{1/4}}\right)^{2i},\]

\[\dsp q^j=\left( \sum_{m=0}^{M}\tilde q_m \tilde h_m+ \hat q_m \hat h_m+q_-\right)^j,\;\;q^k=\left( \sum_{m=0}^{M}\tilde q_m \tilde h_m+ \hat q_m \hat h_m+q_-\right)^k,\]
where $b_{j,k}^{l,i}$ are the coefficients of the polynomials $B_{j,k}^{l}$. Using the fact that $\|q(s)\|_{L^\infty}\leq 1$ (which holds for $s$ large enough, from the fact that $q(s)\in \Vg_A(s)$ and (i) of Proposition \ref{thess}). We deduce that
\[|q^j-q_{+}^{j}|\leq C(|q_-|^j+|q_-|).\]
Using that $q(s)\in \Vg_A(s)$ and the fact that $\sqrt s\geq 2 A^2$, we deduce that in the region $|y|\leq s^{1/4}$, we have $ |q_-|\leq \frac{1}{s^{1/4}}(\frac{A}{s^{1/4}})^{M+1}(1+|y|)^{M+1}$ and that
\[\dsp |q^j-(\sum_{m=0}^{M} \tilde q_m \tilde h_m+\hat q_m \hat h_m)^j| \leq C  \left(\frac{A}{s^{1/4}}\right)^{M+1} \frac{1}{s^{1/4}}(1+|y|)^{jM+j}.\]
In the same way, we have
 \[\dsp |q^k-(\sum_{m=0}^{M} \tilde q_m \tilde h_m+\hat q_m \hat h_m)^k| \leq C  \left(\frac{A}{s^{1/4}}\right)^{M+1} \frac{1}{s^{1/4}}(1+|y|)^{kM+k},\]
hence, the contribution coming from $q_-$ is controlled by the right-hand side of  \eqref{lembin1}. Moreover for all $j,k$ and $l$, we have
\beqtn
\dsp\left | \int_{|y|< s^{1/4}} h_n\rho  B_{j,k}^{l}(\frac{y}{s^{1/4}}) q_{+}^{j}\bar q^{k}_{+}-\int h_n\rho  B_{j,k}^{l}(\frac{y}{s^{1/4}}) q_{+}^{j}\bar q^{k}_{+} \right |\leq C e^{-C\sqrt s}.
\label{lembin4}
\eeqtn
To compute the second term on the left had side of \eqref{lembin4}, we notice that $B_{j,k}^{l}(\frac{y}{s^{1/4}}) q_{+}^{j}\bar q^{k}_{+}$  is a polynomial in $y$ and that the coefficient of the term of degree $n$ is controlled by the right had side of \eqref{lembin1} since $q\in \Vg_A$.\\
Moreover, using that $\sqrt s\geq 2 A^2$, we infer that $|q|\leq \frac{1}{s^{1/4}}(1+|y|)^{M+1}$ in the region $|y|\leq s^{1/4}$
 and hence for all $j,k$ and $l$, we have
 \[\dsp \left| \int_{|y|< s^{1/4}}h_n \rho \frac{1}{s^{l/2}} \tilde B_{j,k}^{l}(y,s) q^j \bar q^k \right| \leq C\frac{1}{s^{\frac l 2+\frac{M+1+j+k}{4}}} \]
and
\[\dsp \left|\int _{|y|< s^{1/4}} h_n \rho(|q|^{M+2}+\frac{1}{s^{\frac{M+1}{2}}}) \right|\leq C \frac{1}{s^{\frac{M+2}{4}}}.\]
The terms appearing in these two inequalities are controlled by the right hand side of \eqref{lembin1}.\\
Using the fact that $\|q(s)\|_{L^\infty}\leq 1$ and \eqref{eqqd}, we remark that $|B(q,y,s)|\leq C$. Since $|\rho(y)|\leq C e^{-s\sqrt s} $ for all $|y|>s^{1/4}$, it holds that
\[\left | \int_{|y|>s^{1/4}}h_n B(q,y,s)\rho dy \right|\leq Ce^{-C\sqrt s}.\]
This concludes the proof of Lemma \ref{lembin}. $\blacksquare$

\medskip

\textbf{Sixth term: $R^*(\theta',y,s)$}\\
In the following, we expand $R^*$ as a power series of $\frac 1 s$ as $s\to\infty$, uniformly for $|y|\leq s^{1/4}$.
\begin{lemma}[Power series of $R^*$ as $s\to\infty$] For all $n\in\N$,
\beqtn
R^*(\theta',y,s)=\Pi_n(\theta',y,s)+\tilde\Pi_n(\theta',y,s),
\eeqtn
where,
\beqtn
\Pi_n(\theta',y,s)=\sum_{k=0}^{n-1}\frac{1}{s^{\frac{k+1}{2}}} P_k(y)-i\theta'(s)\left(\frac{a}{s^{1/2}}(1+i\delta)+\sum_{k=0}^{n-1}e_k\frac{y^{2k}}{s^{k/2}}\right),
\eeqtn
and
\beqtn
\forall |y|<s^{1/4},\;\;\left |\tilde\Pi_n(\theta',y,s)\right|\leq C(1+s|\theta'(s)|)\frac{(1+|y|^{2n})}{s^{\frac{n+1}{2}}},
\eeqtn
where $P_k$ is a polynomial of order $2k$ for all $k\geq 1$ and $e_k\in \R$.\\
In particular,
\beqtn
\begin{array}{l}
\dsp \sup_{|y|\leq s^{1/4}}\left | R^*(\theta',y,s)-\sum_{k=0}^{1} \frac{1}{s^{\frac{k+1}{2}}}P_k(y)+i\theta'\left[\kappa+\frac{(1+i\delta)}{s^{1/2}}\left( a-\frac{b\kappa y^2}{(p-1)^2}\right) \right]\right |\\
\leq C\left (\frac{1+|y|^4}{s^{3/2}}+C|\theta '|\frac{y^4}{s^2}\right).
\end{array}
\eeqtn
\label{decompR}
\end{lemma}
\textit{Proof: }
Using the definition of $\varphi$ \eqref{defifi}, the fact that $\varphi_0$ satisfies \eqref{eqfi0} and \eqref{eqq}, we see that $R^*$ is in fact a function of $\theta'$, $z=\frac{y}{s^{1/4}}$ and $s$ that can be written as 
\beqtn
\begin{array}{l}
\displaystyle R^*(\theta',y,s)=\frac{1}{4}\frac{z}{s}\nabla_z \varphi_0(z)+\frac{1}{2}\frac{a}{s^{3/2}}(1+i\delta)+ \frac{1}{s^{1/2}}\Delta_z \varphi_0(z)\\
\displaystyle -\frac{a(1+i\delta)^2}{(p-1)s^{1/2}}+\left( F\left(\varphi_0(z)+\frac{a}{s^{1/2}}(1+i\delta)\right)-F(\varphi_0)\right)\\
\displaystyle -i\frac{\mu}{s}\left( \varphi_0(z)+\frac{a}{s^{1/2}} (1+i\delta)\right)-i\theta'(s)\left( \varphi_0(z)+\frac{a}{s^{1/2}} (1+i\delta)\right)\mbox{ with, }F(u)=(1+i\delta)|u|^{p-1}u
\end{array}
\label{eqreste}
\eeqtn
Since $|z|<1$, there exists positive $c_0$ and $s_0$ such that $|\varphi_0(z)|$ and $|\varphi_0(z)+\frac{a}{s^{1/2}}(1+i\delta)|$ are both larger that $\frac{1}{c_0}$ ans smaller than $c_0$, uniformly in $|z|<1$ and $s>s_0$. Since $F(u)$ is $C^\infty$ for $\frac{1}{c_0}\leq |u|\leq c_0$, we expand it around $u=\varphi_0(z)$ as follows
\[\left | F\left(\varphi_0(z)+\frac{a}{s^{1/2}}(1+i\delta)\right)-F\left(\varphi_0(z)\right)-  \sum_{j=1}^{n}\frac{1}{s^{j/2}}F_j\left(\varphi_0(z)\right)\right |\leq C\frac{1}{s^{\frac{n+1}{2}}},\]
where $F_j(u)$ are $C^{\infty}$. Hence, we can expand $F_j(u)$ around $u=\varphi_0(0)$ and write 
\[\left | F\left(\varphi_0(z)+\frac{a}{s^{1/2}}(1+i\delta)\right)-F\left(\varphi_0(z)\right)-  \sum_{j=1}^{n}\sum_{l=0}^{n-j}  \frac{c_{j,l}}{s^{j/2}} z^{2l}\right |\leq\sum_{j=1}^{n}\frac{C}{s^{\frac{j}{2}}}|z|^{2(n-j)+2}+ \frac{C}{s^{\frac{n+1}{2}}},\]
Similarly, we have the following
\[\left|\frac{z}{s}\nabla_z\varphi_0(z)-\frac{|z|^2}{s}\sum_{j=0}^{n-2}d_j z^{2j} \right |\leq \frac{C}{s}|z|^{2n}, \]
\[\left |\frac{1}{s^{1/2}}\Delta_z\varphi_0(z)-\frac{1}{s^{1/2}}\sum_{j=0}^{n-1}b_j z^{2j}\right |\leq \frac{C}{s^{1/2}}|z|^{2n}\mbox{  and  }\left|\varphi_0(z)-\sum_{j=0}^{n-1}e_j z^{2j}\right|\leq C|z|^{2n}.\]
Recalling that $z=\frac{y}{s^{1/4}}$, we get the conclusion of the Lemma. $\blacksquare$

\medskip

In the following, we introduce $F_j(R^*)(\theta,s)$ as the projection of the rest term $R^*(\theta',y,s)$ on the standard Hermite polynomial, introduced in \eqref{eigfun}.
\begin{lemma}[Projection of $R^*$ on the eigenfunctions of $\Lg$]
It holds that $F_j(R^*)(\theta',s)\equiv 0$ when $j$ is odd, and $|F_j(R^*)(\theta',s)|\leq C\frac{1+s|\theta'(s)|}{s^{\frac j 4+\frac 12}} $, when $j$ is even and $j\geq 4$.\\
\[
\begin{array}{lll}
F_0(R^*)(\theta',s)&=&\dsp -i\theta'(s)\left(\kappa+ O(\frac{1}{s^{1/2}})\right)+(1+i\delta)\left\{-2\frac{\kappa b}{(p-1)^2}+a\right\}\frac{1}{s^{1/2}}
+\frac{a(p-\delta^2)}{\kappa}\frac 1s
\\&&\dsp +i \kappa\delta\left \{\frac{a^2}{\kappa^2}(1+p)-\frac \mu\delta+12\frac{ b^2(p+1)}{(p-1)^4}-\frac{4ab (p+1)}{\kappa(p-1)^2}+O(\frac{1}{s^{1/2}})
\right\}\frac{1}{s}+O(\frac{1}{s^{3/2}})
%\\&=&\dsp -i\theta'(s)\left(\kappa+ O(\frac{1}{s^{1/2}})\right)+O(\frac{1}{s^{3/2}})
\end{array}
\]
with the fact that $\delta^2=p$, we obtain

\[\begin{array}{lll}
 F_0(R^*)(\theta',s)&=&\dsp -i\theta '(s)\left(\kappa+ O(\frac{1}{s^{1/2}})\right)+(1+i\delta)\left\{-2\frac{\kappa b}{(p-1)^2}+a\right\}\frac{1}{s^{1/2}}
%+\frac{a(p-\delta^2)}{\kappa}\frac 1s
\\&&\dsp +i \kappa\delta\left \{\frac{a^2}{\kappa^2}(1+p)-\frac \mu\delta+12\frac{ b^2(p+1)}{(p-1)^4}-\frac{4ab (p+1)}{\kappa(p-1)^2}+O(\frac{1}{s^{1/2}})\right\}\frac{1}{s}+O(\frac{1}{s^{3/2}})
\end{array}
\]
%and the chose of $a$, $b$ and $\mu$ satisfying
%\beqtn
%\dsp a=\frac{2\kappa b}{(p-1)^2},
%\label{cond1}
%\eeqtn
%and 
%\beqtn
%\dsp\mu=8\delta (p+1)\frac{b^2}{(p-1)^4}.
%\label{cond2}
%\eeqtn

If $j=2$, then
\[
\begin{array}{lll}
F_2(R^*)(\theta',s))&=&\dsp i\theta'(s)\left(\frac{\kappa b}{(p-1)^2}(1+i\delta)\frac{1}{s^{1/2}}-6\kappa \frac{b^2 (p-\delta^2)}{(p-1)^4}\frac 1s-i6\kappa \frac{b^2 \delta (p+1)}{(p-1)^4}\frac 1s +O(\frac{1}{s^{3/2}})\right)\\
&&\dsp \left(6\kappa\frac{ b^2(p-\delta^2)}{(p-1)^4}-\frac{2ab (p-\delta^2)}{(p-1)^2}\right)\frac{1}{s}+i\left(6\kappa\frac{\delta b^2(p+1)}{(p-1)^4}-\frac{2ab\delta (p+1)}{(p-1)^2}\right)\frac{1}{s}\\
&&+\frac{\kappa b}{(p-1)^2}\left\{-(\frac 12+\mu\delta)+p(p+1)\left(\frac{a^2}{\kappa^2}+60\frac{b^2}{(p-1)^4}
-12 \frac{ab^2}{(p-1)^4}\right)

\right\}\frac{1}{s^{3/2}}+iO(\frac{1}{s^{3/2}})+O(\frac{1}{s^{2}}).\\

\end{array}
\]
Then, by the fact that $\delta^2=p$,
\[
\begin{array}{lll}
F_2(R^*)(\theta',s))&=&\dsp i\theta'(s)\left(\frac{\kappa b}{(p-1)^2}(1+i\delta)\frac{1}{s^{1/2}}-i6\kappa \frac{b^2 \delta (p+1)}{(p-1)^4}\frac 1s +O(\frac{1}{s^{3/2}})\right)\\
%&&\dsp \left(6\kappa\frac{ b^2(p-\delta^2)}{(p-1)^4}-\frac{2ab (p-\delta^2)}{(p-1)^2}\right)\frac{1}{s}
&&\dsp+i\left(6\kappa\frac{\delta b^2(p+1)}{(p-1)^4}-\frac{2ab\delta (p+1)}{(p-1)^2}\right)\frac{1}{s}+iO(\frac{1}{s^{3/2}})\\
&&\dsp+\frac{\kappa b}{(p-1)^2}\left\{-(\frac 12+\mu\delta)+p(p+1)\left(\frac{a^2}{\kappa^2}+60\frac{b^2}{(p-1)^4}
-12 \frac{ab^2}{(p-1)^4}\right)

\right\}\frac{1}{s^{3/2}}+O(\frac{1}{s^{2}}).\\

\end{array}
\]

\end{lemma}
\textit{Proof: } Since $R^*$ is even in the $y$ variables and $f_j$ is odd when $j$ is odd, $F_j(R^*)(\theta',s)\equiv 0$, when $j$ is odd.\\
Now, when $j$ is even, we apply Lemma \ref{decompR} with $n=[\frac j2]$ and write
\[R^*(\theta',y,s)=\Pi _{\frac j2}(\theta',y,s)+O\left(\frac{1+s|\theta'(s)|+|y|^j}{s{\frac j4+\frac 12}}\right),\]
where $\Pi_{\frac j2 }$ is a polynomials in $y$ of degree less than $j-1$. Using the definition of $F_j(R^*)$ (projection on the $h_j$ of $R^*$), we write
\beqtn
\begin{array}{l}
\dsp\int_{\R^N}R^* h_j \rho=\int_{|y|<s^{1/4}} R^* h_j \rho dy+\int_{|y|>s^{1/4}} R^* h_j \rho dy\\

=\displaystyle\int_{|y|< s^{1/4}} \Pi _{\frac j2} h_j \rho dy+O\left( \int_{|y|<s^{1/4}}  \frac{1+s|\theta'(s)|+|y|^j}{s{\frac j4+\frac 12}}h_j \rho dy\right)+\int_{|y|>s^{1/4}} R^* h_j \rho dy\\

\dsp=\displaystyle\int_{\R^N} \Pi _{\frac j2} h_j \rho dy+O\left(  \frac{1+s|\theta'(s)|}{s{\frac j4+\frac 12}}\right)+\int_{|y|>s^{1/4}} R^* h_j \rho dy+\int_{|y|>s^{1/4}} \Pi_{\frac j2} h_j \rho dy

\end{array}
\eeqtn
We can see that $\int_{\R^N} \Pi _{\frac j2} h_j \rho dy=0$ because $h_j$ is orthogonal to all polynomials of degree less than $j-1$. Then, note that both integrals over the domain $\{|y|>s^{1/4}\}$ are controlled by
\[\dsp\int _{s^{1/4}}\left(|R^*(\theta',y,s)|+1+|y|^j\right)(1+|y|^j)\rho dy.\]
Using the fact that $R(y,s)$ measures the defect of $\varphi(y,s)$ from being an exact solution of \eqref{GLauto}. However, since $\varphi$ is an approximate solution of \eqref{GLauto}, one easily derive the fact that 

\beqtn
\begin{array}{l}
\|R(s)\|_{L^\infty}\leq \frac{ C}{\sqrt{s}}\mbox{, and}\\
|R^*(\theta', y,s)|\leq \frac{C}{\sqrt{s}} +|\theta'(s)|. 
\end{array}
\eeqtn 
Using the fact that $|\rho(y)|\leq C e^{\sqrt{-s}}$, for $|y|>s^{1/4}$, we can bound our integral by
\[\dsp C(1+|\theta'(s)|)\int_{\R^N}(1+|y|^j)^2 e^{c\sqrt{-s}}\sqrt{\rho}dy=C(j)(1+|\theta'(s)|)e^{c\sqrt{-s}}.\]
This inequality gives us the result for $j\geq 4$.\\ 
If $j=0$ or $j=2$, one has to refine Lemma \ref{decompR} in straightforward but long way and do as we did for general $j$. The details are given in Appendix \ref{Detailscal}. This concludes the proof of the lemma. $\blacksquare$

\medskip

\begin{corollary}{ Projection of $R^*$ on the eigenfunctions of $\tilde \Lg$}
\newline

If $j$ is even and $j\geq 4$, then $\tilde P_j(R^*)(\theta',s)$ and $\hat {P}_j(R^*)(\theta',s)$ are $O\left(\frac{1+s|\theta'|}{s^{\frac j4+\frac 12}}\right)$.\\
If $j$ is odd, then  $\tilde P_j(R^*)(\theta',s)=\hat {P}_j(R^*)(\theta',s)=0$.\\
If $j=0$, then,

 \[\begin{array}{lll}
 \dsp\hat P_0(R^*)(\theta',s)&=&-\theta'(s)\left(\kappa+O\left(\frac{1}{s^{1/2}}\right)\right)\\
&&+\dsp\kappa\delta\left \{\frac{a^2}{\kappa^2}\delta(1+p)-\mu\kappa+12\frac{ b^2(p+1)}{(p-1)^4}-\frac{4ab (p+1)}{\kappa(p-1)^2} \right\}\frac 1s+O\left(\frac{1}{s^{3/2}}\right)
\end{array}
\]
and 
\[\dsp \tilde P_0(R^*)(\theta',s)=O\left(\frac{\theta'(s)}{s^{1/2}}\right)+\left\{-2\frac{\kappa b}{(p-1)^2}+a\right\}\frac{1}{s^{1/2}}
%+\frac{a(p-\delta^2)}{\kappa}\frac 1s
+O\left(\frac{1}{s^{3/2}}\right).\]
If $j=2$, then $\dsp \hat P_2(R^*)(\theta',s)=O\left(\frac 1s\right)+O\left(\frac{\theta'(s)}{s^{1/2}}\right)$ and

%and by condition \eqref{cond3}, we obtain
\[
\begin{array}{lll}
\dsp\tilde P_2(R^*)&=&\dsp\theta'(s)\left(-\frac{\kappa b\delta}{(p-1)^2}\frac{1}{s^{1/2}}+O\left(\frac{1}{s}\right)\right)\\
%&&\dsp +(p-\delta^2)\left(6\frac{\kappa b^2}{(p-1)^4}-2\frac{ab}{(p-1)^2}\right)\frac{1}{s}+\\

&&\dsp+\frac{\kappa b}{(p-1)^2}\left\{-\frac 12-\mu\delta+
p(p+1)\left(\frac{a^2}{\kappa^2}+60\frac{b^2}{(p-1)^4}-12 \frac{ab^2}{(p-1)^4}\right )
\right\}\frac{1}{s^{3/2}}+O\left(\frac{1}{s^{2}}\right).

%&=&\dsp\theta'(s)\left(-\frac{\kappa\delta}{8\sqrt{p(p+1)}}\frac{1}{s^{1/2}}+O\left(\frac{1}{s}\right)\right)+O\left(\frac{1}{s^{2}}\right).
\end{array}
\]
\label{projR*}
\end{corollary}
In the following, we give a new version of Corollary \ref{projR*}.
\begin{corollary} If we choose $a$, $b$ and $\mu$ as follows
\beqtn
\begin{array}{l}
\dsp a=\frac{2\kappa b}{(p-1)^2},\\

\dsp\mu=8\delta (p+1)\frac{b^2}{(p-1)^4},\\

(\frac 12+\mu\delta)=p(p+1)\left(\frac{a^2}{\kappa^2}+60\frac{b^2}{(p-1)^4}-12 \frac{ab^2}{(p-1)^4}\right ),
\end{array}
\label{cond}
\eeqtn
If $j$ is even and $j\geq 4$, then $\tilde P_j(R^*)(\theta',s)$ and $\hat {P}_j(R^*)(\theta',s)$ are $O\left(\frac{1+s|\theta'|}{s^{\frac j4+\frac 12}}\right)$.\\
If $j$ is odd, then  $\tilde P_j(R^*)(\theta',s)=\hat {P}_j(R^*)(\theta',s)=0$.\\
If $j=0$, then,

 \[\dsp\hat P_0(R^*)(\theta',s)=-\theta'(s)\left(\kappa+O\left(\frac{1}{s^{1/2}}\right)\right)
%\kappa\delta\left \{\frac{a^2}{\kappa^2}\delta(1+p)-\mu\kappa+12\frac{ b^2(p+1)}{(p-1)^4}-\frac{4ab (p+1)}{\kappa(p-1)^2} \right\}\frac %1s
+O\left(\frac{1}{s^{3/2}}\right)\]
and 
\[\dsp \tilde P_0(R^*)(\theta',s)=O\left(\frac{\theta'(s)}{s^{1/2}}\right)+O\left(\frac{1}{s^{3/2}}\right).\]
If $j=2$, then $\dsp \hat P_2(R^*)(\theta',s)=O\left(\frac 1s\right)+O\left(\frac{\theta'(s)}{s^{1/2}}\right)$ and

\[
\dsp\tilde P_2(R^*)=\dsp\theta'(s)\left(-\frac{\kappa b\delta}{(p-1)^2}\frac{1}{s^{1/2}}+O\left(\frac{1}{s}\right)\right)+O(\frac{1}{s^2})
\]
\label{NprojR*}
\end{corollary}
\begin{rem} It is very important to note that $a$, $b$ and $\mu$ chosen by \eqref{cond} are the same given by our formal approach (See Section \ref{sectFormapp}).
\end{rem}
\medskip

\textbf{Part 2: Proof of Proposition \ref{propode}}
\label{delicatq}
\medskip

In this part, we consider  $A\geq 1$ and take $s$ large enough so that Part 1 applies.\\
(i) We control $\theta '(s)$, from the projection of \eqref{eqq} on $\hat h_0=i h_0$, taking on consideration the modulation, we obtain
\beqtn
 \hat q_0'=-(\frac \mu s+\theta')\left [(1+\delta^2)\tilde q_0+\delta \hat q_0\right ]+\hat P_0(V_1 q)+\hat P_0(V_2 \bar q)+\hat{P}_0(B)+\hat P_0(R^*) 
 \label{Projtild0}
 \eeqtn
By the definition of the shrinking set $\Vg_A$ Definition \ref{defthess}, Corollary \ref{corprojPtls}, Lemma \ref{lembin} and Corollary \ref{NprojR*}, since $q(s)\in \Vg_A(s)$ and $\hat q_0=0$ for all $s\in [\tau,s_1]$, this yields
\[|\theta'(s)|\leq \frac{CA^5}{s^{\frac3 2}}.\]
%\[\hat P_0(R^*)(\theta',s)=-\theta'(s)\kappa+O\left(\frac{\theta'}{s^{1/2}}\right)+O\left(\frac{1}{s^{3/2}}\right)\]
(ii) From the projection of \eqref{eqq}  respectively on $\tilde h_0=(1+i\delta) h_0$ and $\tilde h_1=(1+i\delta) h_1$, we obtain;
\[
\begin{array}{lll}
 \tilde q_0'&=&\tilde q_0+(\frac \mu s+\theta')\left [\delta\tilde q_0+\hat q_0\right ]+\tilde P_0(V_1 q)+\tilde P_0(V_2 \bar q)+\tilde P_0(B)+\tilde P_0(R^*),\\
 \tilde q_1'&=&\frac 12 \tilde q_1+(\frac \mu s+\theta')\left [\delta\tilde q_1+\hat q_1\right ]+\tilde P_1(V_1 q)+\tilde P_1(V_2 \bar q)+\tilde P_1(B)+\tilde P_1(R^*). 
 \end{array}
 \]

The two inequalities in (ii) of Proposition \ref{propode} are  a direct consequence of Definition \ref{defthess}, Remark \ref{bettrestiVi}, Lemma \ref{lembin} and Corollary \ref{NprojR*}, provided that $s_0$ is large enough.

\medskip

(iii) Estimate of $\tilde q_2$. By Corollary \ref{corprojPtls},  Corollary \ref{projR*}, equation \eqref{Projtild0} and the fact that $q(s)\in \Vg_A(s)$, we obtain
%\leq C\frac{A^5}{s^{3/2}}
%\[\left| \hat P_0(V_1 q)+\frac{\tilde q_2}{s^{1/2}} 4\delta b(p+1)\frac{(p+1)}{(p-1)^2}\right |\leq\frac{C(A)}{s^{3/2}}??,
%;\;\left| \hat P_0(V_2 \bar q))-\frac{\tilde q_2}{s^{1/2}} 4\delta b(p+1)\frac{(p-3)}{(p-1)^2} \right|\leq \frac{C(A)}{s^{3/2}}?? \]
%Then, using the fact $\hat q_0=0$, we obtain
\beqtn
\left|\theta' \kappa-16\delta \frac{b(p+1)}{(p-1)^2} \frac{\tilde{q}_2}{s^{1/2}}  \right|\leq\frac{C A^3}{s^{\frac 32}}
\label{estitetaR}
\eeqtn

Let us project the different terms of \eqref{eqq} on $\tilde h_2$. We use Part 1, (Lemma \ref{projPtls} and Corollary \ref{projR*}) and the fact that $q(s)\in \Vg_A(s)$ and we obtain:
\[\begin{array}{l}
\tilde P_2 \left(\frac{\pa q}{\pa s}\right)=\tilde q_2',\\
|\tilde P_2 (-i(\frac \mu s+\theta'(s)q)-\mu\delta \frac{\tilde q_2}{s}|\leq C\frac{A^3}{s^2},\\
\left| \tilde P_2(V_1 q) +\frac{\tilde q_2 }{s} 60b^2(p+1)\frac{1}{(p-1)^2}+ \tilde P_2(V_2 \bar q)-\frac{\tilde q_2 }{s} 60b^2(p+1)\frac{p^2-4p+1}{2(p-1)^3}\right |\leq C\frac{A^3}{s^2}\\
%|\tilde P_2 (V_1 q)- 32 \frac{b(p+1)}{(p-1)} \frac{\tilde q_2}{s^{1/2}}   + 60b^2(p+1)\frac{1}{(p-1)^2}\frac{\tilde q_2 }{s} |\leq \frac{CA^6}{s^2},\\
%|\tilde P_2 (V_2 \bar q)+ 32 \frac{b(p+1)}{(p-1)} \frac{\tilde q_2}{s^{1/2}}   - 60b^2(p+1)\frac{p^2-4p+1}{2(p-1)^3} \frac{\tilde q_2 }{s}|\leq \frac{C A^6}{s^2},\\

%|\tilde P_2(B(q,y,s))|\leq \frac{C}{s^2}
|\tilde P_2 (B(q,y,s))|+|\hat P_n (B(q,y,s)) |\leq \frac{C }{s^2},
\end{array}
\]
and 
\[
| \tilde P_2(R^*(q,y,s))+16\frac{b^2}{(p-1)^4}p(p+1)\frac{\tilde q_2}{s}|\leq \frac{CA^3}{s^2},
\]
here we used \eqref{estitetaR}, the estimate on $\theta'$ and Corollary \ref{NprojR*}.

\label{2generique}

 Adding all these contributions gives $-2$ as the coefficient of $\dsp\frac{\tilde q_2(s)}{s}$ in the following ODE
\[|\tilde q_2' +\frac 2s \tilde{q}_2 |\leq \frac{C A^3}{s^2}.\]

\medskip

(iv) \textit{Estimates of $\hat q_1$, $\hat q_2$, $\hat q_j$ and $\tilde q_j$ for $3\leq j\leq M$:} Using  the definintion of the shrinking set $\Vg_A$ Definition \ref{defthess} and Corollary \ref{corprojPtls}, Lemma \ref{lembin}, Corollary \ref{projR*} from Part 1 and the fact that $q(s)\in \Vg_A(s)$, we see that for all $s\in [\tau,s_1]$, we have
\beqtn
\begin{array}{lll}
\dsp \left|\hat q_1^{'}+\frac 12 \hat q_1\right|\leq C\frac{A^5}{s^{\frac 32}},& \dsp\left|\hat q_{2}^{'}+\hat q_2\right|\leq \frac{C}{s},&  \\
\dsp\left|\hat q_{j}^{'}+\frac j2 \hat q_j\right|\leq C\frac{A^{j-1}}{s^{\frac{j+1}{4}}},&\mbox{if}&3\leq j\leq M,\\
 \dsp\left|\tilde q_{j}^{'} +\frac{j-2}{2}\tilde q_j\right|\leq C\frac{A^{j-1}}{s^{\frac{j+1}{4}}},&\mbox{if}&3\leq j\leq M.
\end{array}
\eeqtn
Integrating this inequalities between $\tau$ and $s_1$ gives the desired estimates.
\subsubsection{ The infinite dimensional part: $q_-$}

We proceed in 2 parts:
\begin{itemize}
\item In Part 1, we project equation \eqref{eqq} to get equations satisfied by $ q_-$. 
\item In Part 2, we prove the estimate on $q_-$.
\end{itemize}
\textbf{Part 1: Projection of equation \eqref{eqq} using the projector $P_-$}
In the following, we will project equation \eqref{eqq} term by term.

\medskip

\textbf{ First term: $\frac{\pa q}{\pa s}$}\\
From \eqref{decomp3}, its projection is
\beqtn
P_-(\frac{\pa q}{\pa s})=\frac{\pa q_-}{\pa s}
\eeqtn
\textbf{ Second term: $\tilde \Lg_0 q$}\\
From \eqref{eqq}, we have the following, 
\[P_-(\tilde \Lg q)= \Lg_0 q_-+P_-[(1+i\delta)\Re q_-].\]
\textbf{ Third term: $-i\left(\frac \mu s+\theta'(s)\right)q$}\\
Since $P_-$ commutes with the multiplication by $i$, we deduce that
\[P_-[-i\left(\frac \mu s+\theta'(s)\right)q]=-i\left(\frac \mu s+\theta'(s)\right)q_-.\]
\textbf{ Fourth term: $V_1 q$ and $V_2 \bar q$}\\
We have the following:
\begin{lemma}[Projection of $V_1 q$ and $V_2 \bar q$]
The projection of $V_1 q$ and $V_2\bar q$ satisfy for all $s\geq 1$,
\beqtn
\left\| \frac{P_-(V_1q)}{1+|y|^{M+1}}\right\|_{L^\infty}\leq \left(\|V_1\|_{L^\infty}+\frac{C}{s^{1/2}} \right)\left\|\frac{q_-}{1+|y|^{M+1}}\right\|_{L^\infty}
+\sum_{n=0}^{M}\frac{C}{s^{\frac{M+1-n}{4}}}(|\hat q_n|+|\tilde q_n|),
\eeqtn
and the same holds for $V_2 \bar q$.
\label{PVneg}
\end{lemma}
Using the fact that $q(s)\in \Vg_A(s)$, we get the following
\begin{corollary}
For all $A\geq 1$, there exists $s_{14}(A)$ such that for all $s\geq s_{14}$, if $q(s)\in \Vg_A(s)$ then,

\[\dsp   \left\|\frac{P_-(V_1 q)}{1+|y|^{M+1}}\right\|_{L^\infty}\leq \|V_1\|_{L^\infty}\left\| \frac{q_-}{1+|y|^{M+1}} \right \|_{L^\infty}+C\frac{A^M}{s^{\frac{M+2}{4}}},\]

and the same holds for $V_2 q$.
\end{corollary}

\textit{Proof of Lemma \ref{PVneg}:} We just give the proof for $V_1 q$ since the proof for $V_2\bar q$ is similar.

From Subsection \ref{sectidecompq}, we write $q=q_++q_-$ and
\[P_-(V_1 q)=V_1 q-P_+(V_1q_-)+P_-(V_1 q_+).\]
Moreover, we claim that the following estimates hold
\[
\begin{array}{lll}
\left\|\frac{V_1 q_-}{1+|y|^{M+1}}\right\|_{L^\infty} &\leq& \|V_1\|_{L^\infty}\left\|  \frac{q_-}{1+|y|^{M+1}}\right\|_{L^\infty}\\
\left\|\frac{P_+(V_1 q_-)}{1+|y|^{M+1}}\right\|_{L^\infty} &\leq& \frac{C}{s^{1/2}}\left\|  \frac{q_-}{1+|y|^{M+1}}\right\|_{L^\infty}.
\end{array}
\]
Indeed, the first one is obvious. To prove the second one, we use \eqref{bdVi1} to show that 

\[|\hat P_{n}(V_1 q_-)|+|\tilde P_{n}(V_1 q_-)|\leq \frac{C}{s^{1/2}}\left\| \frac{q_-}{1+|y|^{M+1}}\right\|_{L^\infty}.\]
To control $P_-(V_1 q_+)=\sum_{n\leq M}P_-(V_1(\hat q_n \hat h_n+\tilde q_n\tilde h_n))$, we argue as follows.\\
If $M-n$ is odd, we take $k=\frac{M-1-n}{2}$ in \eqref{dcVi}, hence
\[
\begin{array}{ll}
P_-(V_1(\hat q_n \hat h_n+\tilde q_n \tilde h_n))=\sum_{j=1}^{k}\frac{1}{s^{\frac j2}} P_-\left (W_{1,j}(\hat q_n\hat h_n+\tilde q_n \tilde h_n)\right )+P_-\left ((\hat q_n \hat h_n+\tilde q_n \tilde h_n)\tilde W_{1,k}\right)
\end{array}
\]
Since $2k+n\leq M$, we deduce that $P_-\left (W_{1,j}(\hat q_n\hat h_n+\tilde q_n \tilde h_n)\right )=0$ for all $0\leq j\leq k$. Moreover, using that
\[|\tilde W_{1,k}|\leq C\frac{(1+|y|^{2k+2})}{s^{\frac{k+1}{2}}}\]
 and applying Lemma  \ref{lemA3}, we deduce that

\beqtn
\left\|\frac{P_-(V_1(\hat q_n \hat h_n+\tilde q_n \tilde h_n))}{1+|y|^{M+1}}\right\|_{L^\infty}\leq C \frac{|\hat q_n|+|\tilde q_n|}{s^{\frac{M+1-n}{4}}}.
\label{estiodd}
\eeqtn
If $M-n$ is even, we take $k=\frac{M-n}{2}$ in \eqref{dcVi} and use that
\[|\tilde W_{1,k}|\leq C\frac{1+|y|^{2k+1}}{s^{\frac k2+\frac 14}},\]
to deduce that \eqref{estiodd} holds. This ends the proof of Lemma \ref{PVneg}.$\blacksquare$

\medskip

\textbf{Fifth term: $B(q,y,s)$.}\\
Using \eqref{estiquadinn}, we have the following estimate from Lemmas \ref{lemA3} and \ref{lemfifB}.
\begin{lemma} For all $K\geq 1$ and $A\geq 1$, there exists $s_{15}(K,A)$ such that for all $s\geq s_{15}$, if $q(s)\in \Vg_A(s)$, then
\beqtn
\left\|\frac{P_-(B(q,y,s))}{1+|y|^{M+1}} \right\|_{L^\infty}\leq C(M)    \left [\left(\frac{A^{M+2}}{s^{\frac 14}}\right)^{\bar p}+\frac{A^{5}}{s^{\frac 12}} \right]\frac{1}{s^{\frac{M+1}{4}}},
\label{Bmoins}
\eeqtn
where $\bar p=\min(p,2)$.
\end{lemma}
\medskip
\textit{Proof: }The proof is very similar to the proof of the previous lemma. From Lemma \ref{lemfifB}, we deduce that for all $s$ there exists a polynomial $B_M$ of degree $M$ in $y$ such that for all $y$ and $s$, we have
\beqtn
|B-B_M(y)|\leq C\left[\left(\frac{A^{M+2}}{s^{\frac 14}}\right)^{\bar p}+ \frac{A^{[5+(M+1)^2]}}{s^{\frac 12}}\right]\frac{(1+|y|^{M+1})}{s^{\frac{M+1}{4}}}.
\label{Bmoin}
\eeqtn
Indeed, we can take $B_M$ to be the polynomial
\[\dsp B_M=P_{+,M}\left [ \sum_{l=0}^{M}\sum_{\begin{array}{l}  0\leq j,k\leq M+1\\2\leq j+k\leq M+1   \end{array}  }\frac{1}{s^{l/2}} \left [ B_{j,k}^{l}(\frac{y}{s^{\frac 14}}) q_{+}^{j}\bar q_{+}^{k}\right ] \right ]  .\]
Then the fact that $B-B_M(y)$ is controlled by the right hand side of \eqref{Bmoin} is a consequence of the following estimates in the outer region and in the inner region.\\
First, in the region $|y|\geq s^{\frac 14}$, we have from Lemma \ref{lemfifB},
\[|B|\leq C|q|^{\bar p}\leq C\left(\frac{A^{M+2}}{s^{\frac 14}}\right)^{\bar p}\]
and from the proof of Lemma \ref{lembin}, we know that for $0\leq n\leq M$,
\[|\tilde P_{n}(B_M(q,y,s))|+|\hat P_{n}(B_M(q,y,s))|\leq C\frac{A^n}{s^{\frac{n+2}{4}}}.\]
Beside, in the region $|y|\leq s^{\frac 14}$, we can use the same argument as in the proof of Lemma \ref{lemfifB} to deduce that the coefficients of degree $k\geq M+1$ of the polynomial
\[\sum_{l=0}^{M}\sum_{\begin{array}{l}  0\leq j,k\leq M+1\\2\leq j+k\leq M+1   \end{array}  }\frac{1}{s^{l/2}} \left [ B_{j,k}^{l}(\frac{y}{s^{\frac 14}}) q_{+}^{j}\bar q_{+}^{k}\right ] -B_M,\]
is controlled by $C\frac{A^k}{s^{\frac k 4+\frac 12}}$
 and hence
 \[\left | \sum_{l=0}^{M}\sum_{\begin{array}{l}  0\leq j,k\leq M+1\\2\leq j+k\leq M+1   \end{array}  }\frac{1}{s^{l/2}} \left [ B_{j,k}^{l}(\frac{y}{s^{\frac 14}}) q_{+}^{j}\bar q_{+}^{k}\right ] -B_M\right|\leq C\frac{A^{2M+2}}{s^{\frac{M+3}{4}}}(1+|y|^{M+1}),\]
in the region $|y|\leq s^{\frac 14}$.\\
Moreover, using that $ |q|\leq C\frac{A^{M+1}}{s^{\frac 14}}$ in the region $|y|\leq s^{\frac 14}$, we deduce that for all $s\geq 2 A^2$, we have 
\[\left|\sum_{l=0}^{M}\sum_{\begin{array}{l}  0\leq j,k\leq M+1\\2\leq j+k\leq M+1   \end{array}  }\frac{1}{s^{l/2}} \left [ B_{j,k}^{l}(\frac{y}{s^{\frac 14}}) q_{+}^{j}\bar q_{+}^{k}\right ] \right|\leq  C\frac{A^{2M+2}}{s^{\frac{M+3}{4}}}(1+|y|^{M+1}).\]
Finally, to control the term $|q|^{M+2}$, we use the fact that in the region $|y|\leq s^{\frac 14}$, we have the following two estimates $|q|\leq C\frac{A^{M+1}}{s^{\frac 14}}$
 and $|q|\leq A^5 \frac{1}{s^{\frac 34}}(1+|y|^{M+1})$ if $\sqrt {s}\geq 2 A^2$. Hence
 \[|q|^{M+2}\leq C\frac{A^5}{s^{\frac 34}}\left(\frac{A^{M+1}}{s^{\frac 14}}\right)^{M+1}(1+|y|^{M+1}).\]
This ends the proof of estimate \eqref{Bmoin} and conclude the proof of \eqref{Bmoins} by applying Lemma \ref{lemA3}.$\blacksquare$

\medskip

\textbf{Sixth term: $R^*(\theta',y,s)$}.\\
We claim the following:
\begin{lemma}
If $|\theta'(s)|\leq \frac{CA^5}{s^{3/2}}$, then the following holds
\[\dsp \left\|\frac{P_-(R^*(\theta',y,s))}{1+|y|^{M+1}}\right\|\leq C\dsp\frac{1}{s^{\frac{M+3}{4}}}\]
\end{lemma}
\textit{Proof: }Taking $n=\frac M2+1$ (remember $M$ is even), we write from Lemma \ref{decompR}
$R^*(\theta',y,s)=\Pi_n(\theta',y,s)+\tilde \Pi_n(\theta',y,s)$. Since $2n-2=M$, we see from subsection \ref{sectidecompq} that
\beqtn
\dsp |\tilde\Pi_n(\theta',y,s)|\leq\dsp C\frac{1+|y|^{2n-2}}{s^{\frac{n+1}{2}}}\leq\dsp C\frac{1+|y|^{M+1}}{s^{\frac{M+3}{4}}}
\label{RPIneq}
\eeqtn
in the region $|y|< s^{1/4}$. It is easy to see using \ref{Rexpect} and the definition of $\Pi_n$ that \eqref{RPIneq} holds for all $y\in\R$ and $s\geq 1$. Then applying Lemma \ref{lemA3}, we conclude easily.$\blacksquare$

\medskip

\textbf{Part 2: Proof of the last but one identity in (iv) of Proposition \ref{propode} (estimate on $q_-$)}\\ 
If we apply the projection $P_-$ to the equation \eqref{eqq} satisfied by $q$, we see that $q_-$ satisfies the following equation:
\[\dsp\frac{\pa q_-}{\pa s}=\Lg_0 q_-+ P_-[(1+i\delta) \Re q_-]+P_-[-i(\frac\mu s+\theta'(s))q+ V_1 q+V_2\bar q+ B(q,y,s)+R^*(\theta',y,s)]. 
\]
Here, we have used the important fact that $P_-[(1+i\delta)\Re q_+]=0$. The fact that $M$ is large as fixed in \eqref{boundM} is crucial in the proof.
Using the kernel of the semigroup generated by $\Lg_0$, we get for all $s\in [\tau,s_1]$,
\[\begin{array}{lll}
q_-(s)&=&\dsp e^{(s-\tau)\Lg_0}q_-(\tau)\\
&&+\dsp \int_{\tau}^{s}e^{(s-s')\Lg_0}P_-[(1+i\delta) \Re  q_-]ds'\\
&&+\dsp \int_{\tau}^{s}e^{(s-s')\Lg_0}P_-\left[-i(\frac\mu s+\theta'(s'))q+ V_1 q+V_2\bar q+ B(q,y,s')+R^*(\theta',y,s')\right]ds'.
\end{array}
\]
Using Lemma \ref{lemA2}, we get
\[\begin{array}{l}
\left\| \frac{q_-(s)}{1+|y|^{M+1}}\right\|_{L^\infty}\leq \dsp e^{-\frac{M+1}{2}(s-\tau)}\left\|\frac{q_-(\tau)}{1+|y|^{M+1}}\right\|_{L^\infty}\\
+\dsp \int_{\tau}^{s}e^{-\frac{M+1}{2}(s-s')} \sqrt{1+\delta^2} \left\|\frac{q_-}{1+|y|^{M+1}}\right\|_{L^\infty}ds'\\
+\dsp \int_{\tau}^{s}e^{-\frac{M+1}{2}(s-s')} \left\|\frac{P_-\left[-i(\frac\mu s+\theta'(s'))q+ V_1 q+V_2\bar q+ B(q,y,s')+R^*(\theta',y,s')\right]}{1+|y|^{M+1}}\right\|_{L^\infty}ds'
\end{array}
\]
Assuming that $q(s')\in \Vg_A(s')$, the results from Part 1 yields (use (i) of Proposition \ref{propode} to bound $\theta'(s)$)
\[\begin{array}{l}
\left\| \frac{q_-(s)}{1+|y|^{M+1}}\right\|_{L^\infty}\leq \dsp e^{-\frac{M+1}{2}(s-\tau)}\left\|\frac{q_-(\tau)}{1+|y|^{M+1}}\right\|_{L^\infty}\\
+\dsp \int_{\tau}^{s}e^{-\frac{M+1}{2}(s-s')} \left(\sqrt{1+\delta^2}+\||V_1 |+|V_2|\|_{L^\infty} \right)\left\|\frac{q_-}{1+|y|^{M+1}}\right\|_{L^\infty}ds'\\
+\dsp C(M)\int_{\tau}^{s}e^{-\frac{M+1}{2}(s-s')} \left [    \frac{A^{(M+1)^2+5}}{(s')^{\frac{M+3}{4}}}   +\frac{A^{(M+2)\bar p}}{(s')^{\frac{\bar p-1}{2}}}\frac{1}{(s')^{\frac{M+2}{2}}}    +\frac{A^M}{(s')^{\frac{M+2}{2}}}              \right]ds'.
\end{array}
\]
Since we have already fixed $M$ in \eqref{boundM} such that
\[\dsp M\geq 4\left( \sqrt{1+\delta^2}+1+2 \max_{i=1,2,y\in \R,s\geq 1}|V_i(y,s)|\right),\]
using Gronwall's lemma or Maximum principle and  \eqref{taus}, we deduce that
\[\begin{array}{ll}
e^{\frac{M+1}{2}s}\left\| \frac{q_-(s)}{1+|y|^{M+1}}\right\|_{L^\infty}&\leq\dsp  e^{\frac{M+1}{4}(s-\tau)}e^{\frac{M+1}{2}\tau}\left\| \frac{q_-(\tau)}{1+|y|^{M+1}}\right\|_{L^\infty}\\
&+\dsp e^{\frac{M+1}{2}s}  2^{\frac{M+3}{4}}   \left [    \frac{A^{(M+1)^2+5}}{s^{\frac{M+3}{4}}}   +\frac{A^{(M+2)\bar p}}{s^{\frac{\bar p-1}{2}}}\frac{1}{(s')^{\frac{M+2}{2}}}    +\frac{A^M}{s^{\frac{M+2}{2}}}              \right]
\end{array}
\]
which concludes the proof of the last but one identity in (iv) of Proposition \ref{propode}.

\subsubsection{The outer region: $q_e$}\label{outreg}
Here, we finish the proof of Proposition \ref{propode} by proving the last inequality in (iv). Since $q(s)\in \Vg_A(s)$ for all $s\in [\tau,s_1]$, it holds from Claim \ref{propshrinset} and Proposition \ref{propode} that
\beqtn
\|q(s)\|_{L^\infty(|y|<2K s^{1/4})}\leq C\frac{A^{M+1}}{s^{1/4}}\mbox{ and }|\theta'(s)|\leq \frac{CA^5}{s^{3/2}}.
\label{estheq}
\eeqtn
Then, we derive from \eqref{eqq1} an equation satisfied by $q_e$, where $q_e$ is defined by \eqref{defiqe}:
\beqtn
\begin{array}{ll}
\dsp \frac{\pa q_e}{\pa s}&= \Lg_0 q_e-\frac{1}{p-1} q_e +(1-\chi)e^{\frac{i\delta}{p-1}s}\left\{L(q,\theta',y,s)+ R^*(\theta',y,s)\right\}\\
&\dsp -e^{\frac{i\delta}{p-1}s}q(s)\left( \pa_s \chi+\Delta \chi+\frac 12y\cdot\nabla \chi\right)+2 e^{\frac{i\delta}{p-1}s}div(q(s)\nabla \chi).
\end{array}
\eeqtn
Writing this equation in its integral form and using the maximum principle satisfied by $e^{\tau \Lg_0}$ (see Lemma \ref{lemA1}, see Apendix below), we write
\[
\begin{array}{ll}
\|q_e(s)\|_{L^\infty}& \leq \dsp e^{-\frac{s-\tau}{p-1}}\|q_e(\tau)\|_{L^\infty},\\
&+\dsp \int_{\tau}^{s}e^{-\frac{s-s'}{p-1}}\left(\|(1-\chi)L(q,\theta',y,s')\|_{L^\infty}+\|(1-\chi)R^*(\theta',y,s')\|_{L^\infty}\right)ds'\\
&+\dsp \int_{\tau}^{s} e^{-\frac{s-s'}{p-1}}\left\| q(s')\left(\pa_s \chi+\Delta \chi +\frac 12 y\cdot \nabla\chi\right)\right\|_{L^\infty}ds'\\
&+\dsp \int_{\tau}^{s} e^{-\frac{s-s'}{p-1}} \frac{1}{\sqrt{1-e^{-(s-s')}}}\| q(s')\nabla \chi\|_{L^\infty}ds'.
\end{array}
\]
Let us bounds the norms in the three last lines of this inequality.\\
First from \eqref{defchi14} and \eqref{estheq}
\beqtn
\begin{array}{ll}
\left\|q(s')\left(\pa_s \chi+\Delta \chi +\frac 12 y\cdot \nabla\chi\right)\right\|_{L^\infty}&\leq C(1+\frac{1}{K^2 s'})\|q(s')\|_{L^\infty(|y|<2 K {s'}^{1/4})}\\
&\dsp\leq C \frac{A^{M+1}}{{(s')}^{1/4}},
\label{estiqe1}
\end{array}
\eeqtn
\beqtn
\|q(s')\nabla \chi\|_{L^\infty}\leq \frac{C}{K {(s')}^{1/4}}\|q(s')\|_{L^\infty(|y|<2K {(s')}^{1/4})}\leq C\frac{A^{M+1}}{\sqrt{s'}},
\label{estiqe2}
\eeqtn
for $s'$ large enough.\\
Second note that the residual term $(1-\chi)R^*$ is small as well. Indeed, recalling the bound \eqref{estR*} on $R$, we write from the definition of $R^{*}$ \eqref{eqqd1} and \eqref{estheq}:
\beqtn
\|(1-\chi)R^*(\theta',y,s')\|_{L^\infty}\leq \frac{C}{{(s')}^{1/4}}+|\theta'(s')|\leq\frac{C}{{(s')}^{1/4}}
\label{estiqe3}
\eeqtn
for $s'$ large enough.\\
Third, the term $(1-\chi)L(q,\theta',y,s')$ given in \eqref{eqqd1} is less than $\epsilon |q_e| $ with $\epsilon =\frac{1}{2(p-1)}$. Indeed, it holds from \eqref{estheq} that: 
\beqtn
\begin{array}{l}
\dsp \|(1-\chi)L(q,\theta',y,s')\|_{L^\infty}\\
\leq \dsp C\|q_e(s')\|_{L^\infty}\left( \|\varphi(s')\|_{L^\infty(|y|\geq K{s'}^{1/4})}^{p-1} +\|q(s')\|_{L^\infty(|y|\geq K{s'}^{1/4})}^{p-1}+\frac{1}{s'}+|\theta'(s)|\right),\\
\leq \frac{1}{2(p-1)}\|q_e(s')\|_{L^\infty},
\label{estiqe4}
\end{array}
\eeqtn
whenever $K$ and $s'$ are large (in order to ensure that $\|\varphi(s')\|_{L^{\infty}(|y|\geq K{s'}^{1/4})}$ is small).\\
Notice that it is only here that we need the fact that $K$ is big enough. Using estimates \eqref{estheq}, \eqref{estiqe1}, \eqref{estiqe2}, \eqref{estiqe3} and \eqref{estiqe4}, we write
\[\begin{array}{ll}
\|q_e(s)\|_{L^\infty}&\leq \dsp e^{-\frac{s-\tau}{p-1}}\|q_e(\tau)\|_{L^\infty}\\
&\dsp + \int_{\tau}^{s}e^{-\frac{s-s'}{p-1}}\left(\frac{1}{2(p-1)}\|q_e(s')\|_{L^\infty}+C\frac{A^{M+1}}{(s')^{\frac 14}}+C\frac{A^{M+1}}{(s')^{\frac 12}}\frac{1}{\sqrt{1-e^{-(s-s')}}}\right)ds'.

\end{array}
\] 
Using Gronwall's inequality or Maximum principle, we end-up with
\[\|q_e(s)\|_{L^\infty}\leq e^{-\frac{(s-\tau)}{2(p-1)}}\|q_e(\tau)\|_{L^\infty}+\frac{CA^{M+1}}{ \tau^{\frac 14}}(s-\tau+\sqrt{s-\tau}),\]
which concludes the proof of Proposition \ref{propode}.

%%%%%%%%%%%%%%%

\section{Single point blow-up and final profile}
In this section, we prove Theorem \ref{thm1}\label{pth1}.
Here, we use the solution of problem \eqref{eqq}-\eqref{eqmod} constructed in the last section to exhibit a blow-up solution of equation \eqref{GL} and prove Theorem \ref{thm1}.\\
(i) Consider $(q(s), \theta(s))$ constructed in Section \ref{existence} such that  \eqref{limitev} holds. From \eqref{limitev} and the properties of the shrinking set given in Claim \ref{propshrinset}, we see that $\theta(s)\to \theta_0$ as $s\to \infty$ such that
\beqtn
\label{estitetas}
|\theta(s)-\theta_0|\leq CA^5\int_{s}^{\infty}\frac{1}{\tau^{\frac 32}}d\tau\leq \frac{CA^5}{\sqrt s} \mbox{ and }\|q(s)\|_{L^\infty(\R)}\leq \frac{C_0(K,A)}{\sqrt s}.
\eeqtn
Introducing $w(y,s)=e^{i(\mu\log s+\theta(s))}\left( \varphi(y,s)+q(y,s)\right)$, we see that $w$ is a solution of equation \eqref{GLauto} that satisfies for all $s\geq \log T$ and $y\in\R$,
\[|w(y,s)-e^{i\theta_0+i\mu\log s}\varphi(y,s)|\leq C\|q(s)\|_{L^\infty}+C|\theta(s)-\theta_0|\leq \frac{C_0}{s^{\frac 14}}.\] 
Introducing
\[u(x,t)=e^{-i\theta_0}\kappa^{i\delta}(T-t)^{\frac{1+i\delta}{p-1}}w\left(\frac{y}{\sqrt{T-t}},-\log(T-t)\right),\]
we see from \eqref{chauto} and the definition of $\varphi$ \eqref{defifi} that $u$ is a solution of equation \eqref{GL} defined for all $(x,T)\in\R\times [0,T)$ which satisfies \eqref{Preslt}.\\
If $x_0=0$. It remains to prove that when $x_0\not =0$, $x_0$ is not a blow-up point. The following result from Giga and Kohn \cite{GKCPAM89} allows us to conclude:
\begin{prop}[Giga and Kohn - No blow-up under the ODE threshold] For all $C_0>0$, there is $\eta_0>0$ such that if $v(\xi,\tau)$ solves
\[\left |  v_t -\Delta v\right |\leq C_0 (1+|v|^p)\]
and satisfies
\[|v(\xi,\tau)|\leq \eta_0(T-t)^{-1/(p-1)}\]
for all $(\xi,\tau)\in B(a,r)\times[T-r^2,T)$ for some $a\in \R$ and $r>0$, then $v$ does not blow up at ($a,T$).
\label{propositionGK}
\end{prop}
\textit{Proof: } See Theorem 2.1 page 850 in \cite{GKCPAM89}. $\blacksquare$\\

Indeed, we see from \eqref{Preslt} and \eqref{defifi} that
\[\sup_{|x-x_0|\leq\frac{|x_0|}{2}}(T-t)^{\frac{1}{p-1}}|u(x,t)|\leq \left |  \varphi_0\left( \frac{|x_0|/2}{\sqrt{(T-t)}|\log(T-t)|}\right)\right |+\frac{C}{|\log (T-t)|^{\frac 14}}\to 0\]
as $t\to T$, $x_0$ is not a blow-up point of $u$ from Proposition \ref{propositionGK}. This concludes the proof of (i) of Theorem \ref{thm1}.\\
(ii)  Arguing as Merle did in \cite{FMCPAM92}, we derive the existence of a blow-up profile $u^*\in C^2(\R^*)$ such that $u(x,t)\to u^*(x)$ as $t\to T$, uniformly on compact sets of $\R^*$. The profile $u^*(x)$ is not defined at the origin. In the following, we would like to find its equivalent as $x\to 0$ and show that it is in fact singular at the origin. We argue as in Masmoudi and Zaag \cite{MZ07}. Consider $K_0>0$ to be fixed large enough later. If $x_0\neq 0$ is small enough, we introduce for all $(\xi,\tau)\in \R\times [-\frac{t_0(x_0)}{T-t_0(x_0)},1)$, 

\begin{align}
v(x_0,\xi,\tau)&=(T-t_0(x_0))^{\frac{1+i\delta}{p-1}} v(x,t),\\
\mbox{where,   }x&=x_0+\xi\sqrt{T-t_0(x_0)},\; t=t_0(x_0)+\tau(T-t_0(x_0)),
\label{defVN}
\end{align}

and $t_0(x_0)$ is uniquely determined by
\beqtn
|x_0|=K_0\sqrt{(T-t_0(x_0))|\log(T-t_0(x_0))|^{\frac 12}}.
\label{xt0N}
\eeqtn
From the invariance of problem (\ref{GL}) under dilation, $v(x_0,\xi,\tau)$ is also a solution of (\ref{GL}) on its domain. From (\ref{defVN}), \eqref{xt0N}, \eqref{defifi}, we have
\[\sup_{|\xi|<2|\log(T-t_0(x_0))|^{1/8}}\left |v(x_0,\xi,0)-\varphi_0(K_0) \right|\leq \frac{C}{|\log(T-t_0(x_0))|^{\frac 18}}\to 0\mbox{ as }x_0\to 0. \]
Using the continuity with respect to initial data for problem (\ref{GL}) associated to a space-localization in the ball $B(0,|\xi|<|\log(T-t_0(x_0))|^{1/8})$, we show as in Section 4 of \cite{ZAIHPANL98} that
\[
\begin{array}{l}
\sup_{|\xi|\leq |\log(T-t_0(x_0))|^{1/8},\;0\leq\tau<1}|v(x_0,\xi,\tau)-U_{K_0}(\tau)|\leq \epsilon(x_0)\mbox{ as }x_0\to 0, \\
\end{array}
\]
where $U_{K_0}(\tau)=((p-1)(1-\tau)+b K_{0}^{2})^{-\frac{1+i\delta}{p-1}}$ is the solution of the PDE (\ref{GL}) with constant initial data $\varphi_0(K_0)$. Making $\tau \to 1$ and using (\ref{defVN}), we see that 
\[
\begin{array}{lll}
u^*(x_0)=\lim_{t\to T} v(x,t)&=&(T-t_0(x_0))^{-\frac{1+i\delta}{p-1}} |\log (T-t_0(x_0))|^{i\mu} \lim_{\tau \to 1} v(x_0,0,\tau)\\
&\sim&(T-t_0(x_0))^{-\frac{1+i\delta}{p-1}} |\log (T-t_0(x_0))|^{i\mu}U_{K_0}(1)
\end{array}
\]
as $x_0\to 0$. Since we have from (\ref{xt0N}) 
\[\log(T-t_0(x_0)) \sim 2 \log |x_0|\mbox{ and } T-t_0(x_0)\sim \frac{|x_0|^2}{\sqrt 2K_{0}^{2}\sqrt{|\log|x_0||}},\]
as $x_0\to 0$, this yields (ii) of Theorem \ref{thm1} and concludes the proof of Theorem \ref{thm1}. $\blacksquare$

\appendix
\section{Spectral properties of $\Lg_0$}
In this Appendix, we give some properties associated to the operator $\Lg_0$, defined in \eqref{eqqd1};
\[
 \Lg_0 v=\Delta v-\frac 12 y\cdot \nabla v.
\]
Indeed as stated by  \cite{BKN94}, we have the following Mehler's kernel for the semigroup  $\Lg_0$ l:
\beqtn
e^{s\Lg_0}(y,x)=\frac{1}{[4\pi(1-e^{-s})]^{N/2}}\exp\left [ -\frac{|x-ye^{-\frac s2}|}{4(1-e^{-s})}\right ].
\label{defsemgr}
\eeqtn
In the following, we give some properties associated to the kernel.
\begin{lemma}
a) The semigroup associated to $ \Lg_0$ satisfies the maximum principle:
\[
\|e^{s \Lg_0}\varphi\|_{L^\infty}\leq \|\varphi\|_{L^\infty}.
\]
b) Moreover, we have
\[\|e^{s \Lg_0} \div (\varphi)\|_{L^\infty}\leq \frac{C}{\sqrt{1-e^{-s}}}\|\varphi\|_{L^\infty},  \]
where $C$ is a constant.
\label{lemA1}
\end{lemma}
\textit{Proof: } a) It follows directly by part, this also follows from the definition of the semigroup \eqref{defsemgr}.\\
b) Using an integration by part, this also follows from the definition of the semigroup \eqref{defsemgr}.
\begin{lemma}
There exists a constant $C$ such that if $\phi$ satisfies
\[\forall x\in\R\;\; |\phi(x)|\leq (1+|x|^{M+1})\]
then for all $y\in \R$, we have
\[|e^{s \Lg_0}P_-(\phi(y))|\leq Ce^{-\frac{M+1}{2}s}(1+|y|^{M+1})\]
\label{lemA2}
\end{lemma}
\textit{Proof: } This also follows directly from the semigroup's definition, through an integration by part, for a similar case see page 556-558 from \cite{BKN94}.

Moreover, we have the following useful lemma about $P_-$.
\begin{lemma}
For all $k\geq 0$, we have
\[\left\| \frac{P_-(\phi)}{1+|y|^{M+k}}    \right\|_{L^\infty}\leq C\left\|\frac{\phi}{1+|y|^{M+k}}     \right\|.\]
\label{lemA3}
\end{lemma}
\textit{Proof: } Using \eqref{eqQn}, we have
\[|\phi_n| \leq C   \left\|\frac{\phi}{1+|y|^{M+k}}\right\|_{L^\infty}. \]  
Since for all $m\leq M$, $|h_m(y)|\leq C(1+|y|^{m+k})$ and 
\[|\phi|\leq   C\left\|\frac{\phi}{1+|y|^{M+k}} \right\|_{L^\infty}(1+|y|^{m+k}),\]
the result follows from definition \eqref{decomp1} of $\phi$.
\section{Details of expansions of the fourth term of equation \eqref{eqqd1} $V_1 q+V_2 \bar q$ } \label{detailsV1V2}
In the following, we will try to expand each term of $V_1(y,s)$ and $V_2(y,s)$ as a power series of $\dsp\frac 1s$ as $s\to \infty$, uniformly for $\dsp |y|\leq C\dsp s^{\frac 14}$, with $C$ a positive constant.

\beqtn
\begin{array}{l}
\left |\varphi_0(z)+\frac{a}{s^{1/2}}(1+i\delta)\right|^{p-1}\\
=\kappa^{p-1}\Big(1+\frac{(p-1)a}{\kappa s^{1/2}}+\frac{a^2(p-1)^2}{\kappa^2s}-\frac{b}{(p-1)}\frac{y^2}{s^{1/2}}+\frac{a^3}{3\kappa^3 s^{3/2}}(p-1)(p-3)(2p-1)\\

-\frac{2ab}{\kappa}\frac{y^2}{s}+\frac{a^2b(p-3)(1-2p)}{\kappa^2(p-1)}\frac{y^2}{s^{3/2}}\\

+\frac{b^2}{(p-1)^2}\frac{y^4}{s}

+\frac{ab^2}{2\kappa (p-1)^3}\left[ (p+1)(p-2)+2(p-1)(p-3)+2(p-1)^2\right ]\frac{y^4}{s^{3/2}}\\ 
-\frac{b^3}{(p-1)^3}\frac{y^6}{s^{3/2}}
+O(\frac{1}{s^{2}})+O(\frac{y^8}{s^{2}})\Big).
%\kappa^{p-1}\left(1+\frac{(p-1)a}{\kappa s^{1/2}}+\frac{a^2(p-1)^2}{\kappa^2s}-\frac{b}{(p-1)}\frac{y^2}{s^{1/2}}-\frac{2ab}{\kappa}\frac{y^2}{s}+\frac{b^2}%{(p-1)^2}\frac{y^4}{s}+O(\frac{y^2}{s^{3/2}})+O(\frac{y^4}{s^{3/2}})+O(\frac{y^6}{s^{3/2}})\right).
\end{array}
\eeqtn

\[
\begin{array}{l}
V_1=(1+i\delta)\frac{p+1}{2}\left(\left |\varphi_0(z)+\frac{a}{s^{1/2}}(1+i\delta)\right|^{p-1}-\frac{1}{p-1}\right)\\
=\left(\frac{1}{s^{1/2}}W_{1,1}+\frac 1 s W_{1,2}+O(\frac{1}{s^{3/2}})+O(\frac{y^6}{s^{3/2}}) \right)
\end{array}
\]
Recalling that $ b=\frac{(p-1)^2}{8\sqrt{p(p+1)}},\;\;\; a=\frac{\kappa}{4\sqrt{p(p+1)}},\;\;$, we have

\[
\begin{array}{ll}
W_{1,1}&=-(1+i\delta)\frac{b(p+1)}{2(p-1)^2}h_2(y)\\
F_2(W_{1,1}\tilde q_2\tilde h_2)&=-32(1+i\delta)^2\frac{b(p+1)}{(p-1)^2}\tilde q_2\\
\tilde P_2(W_{1,1}\tilde q_2\tilde h_2)&=32\frac{b(p+1)}{(p-1)}\tilde q_2,\\
F_0(W_{1,1}\tilde q_2\tilde h_2)&=-4(1-p+2i\delta)\frac{b(p+1)}{(p-1)^2}\tilde q_2\\
\hat P_0(W_{1,1}\tilde q_2\tilde h_2)&=-4\delta\frac{b(p+1)^2}{(p-1)^2}\tilde q_2
\end{array}\]

\[
\begin{array}{ll}
W_{1,2}&=(1+i\delta)\frac{b^2(p+1)}{2(p-1)^3}\left (y^4-4y^2+4\right)=(1+i\delta)\frac{b^2(p+1)}{2(p-1)^3}h_{2}^{2}(y),\\
F_2(W_{1,2}\tilde q_2\tilde h_2)&=60(1+i\delta)^2\frac{b^2(p+1)}{(p-1)^3}\tilde q_2\\
\tilde P_2(W_{1,2}\tilde q_2\tilde h_2)&=-60\frac{b^2(p+1)}{(p-1)^2}\tilde q_2
\end{array}
\]

%%and
%\[W_{1,3}=(1+i\delta)\frac{b^3(p+1)}{2(p-1)^6}
%\left(y^6+(p-1)(p-2)(5p-3)y^4-(p-1)(p-3)(2p-1)y^2+\frac 83 (p-1)(p-3)(2p-1)\right).
%\]

Let us now expand the term $V_2$;

\beqtn
\begin{array}{l}
\left |\varphi_0(z)+\frac{a}{s^{1/2}}(1+i\delta)\right|^{p-3}\\

=\kappa^{p-3}\left(1+\frac{2 a}{\kappa s^{1/2}}+\frac{a^2(p+1)}{\kappa^2s}\right)^{\frac{p-3}{2}}\Big (1-2\frac{b}{(p-1)^2}\frac{y^2}{s^{1/2}}-2\frac{ab}{\kappa(p-1)}\frac{y^2}{s}+\frac{b^2(p+1)}{(p-1)^4}\frac{y^4}{s}\\

+2\frac{a^2b (3p-1)}{\kappa^2(p-1)^2}\frac{y^2}{s^{3/2}}+\frac{ab^2(p+1)( p-2)}{\kappa(p-1)^4}\frac{y^4}{s^{3/2}}
-\frac 23\frac{b^3p (p+1)}{(p-1)^6}\frac{y^6}{s^{3/2}}
+O(\frac{y^6}{s^{2}})+O(\frac{y^8}{s^{2}})\Big)^{\frac{p-3}{2}}\\
=\kappa^{p-3}\left(1+\frac{a(p-3)}{\kappa}\frac{1}{s^{1/2}}+\frac{a^2(p-3)(p-2)}{\kappa^2}\frac 1s+\frac{2}{3}\frac{a^3(p-1)(p-3)(p-5)}{\kappa^3}\frac{1}{s^{3/2}}+O(\frac{1}{s^2})\right)\\
*\Big(1-\frac{b(p-3)}{(p-1)^2}\frac{y^2}{s^{1/2}}-\frac{ab(p-3)}{\kappa (p-1)}\frac{y^2}{s}
+\frac{b^2(p-3)(p-2)}{(p-1)^4}\frac{y^4}{s}+
\frac{ab^2(p-3)^2(2p-1)}{2\kappa (p-1)^4}\frac{y^4}{s^{3/2}}
-\frac{b^3(p-2)(p-3)(3p-5)}{3(p-1)^6}\frac{y^6}{s^{3/2}}
+O(\frac{y^6}{s^{2}})+O(\frac{y^8}{s^{2}})\Big)\\

=\kappa^{p-3}\Big(1+\frac{a(p-3)}{\kappa}\frac{1}{s^{1/2}}+\frac{a^2(p-3)(p-2)}{\kappa^2}\frac 1s+\frac{2}{3}\frac{a^3(p-1)(p-3)(p-5)}{\kappa^3}\frac{1}{s^{3/2}}
-\frac{b(p-3)}{(p-1)^2}\frac{y^2}{s^{1/2}}-\frac{2ab(p-3)(p-2)}{\kappa(p-1)^2}\frac{y^2}{s}\\
+\frac{b^2(p-3)(p-2)}{(p-1)^4}\frac{y^4}{s}+\frac{ab^2(p-3)^2(4p-5)}{2\kappa(p-1)^4}\frac{y^4}{s^{3/2}}
-\frac{b^3(p-2)(p-3)(3p-5)}{3(p-1)^6}\frac{y^6}{s^{3/2}}+O(\frac{1}{s^2})+O(\frac{y^8}{s^{2}})

\Big)

\end{array}
\eeqtn

\[
\begin{array}{l}
\left (\varphi_0(z)+\frac{a}{s^{1/2}}(1+i\delta)\right)^{2}\\
=\kappa^2\left(1+\frac{a}{\kappa\sqrt{s}}(1+i\delta)-\frac{ b}{(p-1)^2}(1+i\delta)\frac{y^2}{s^{1/2}}+i\frac{b^2\delta(p+1)}{2(p-1)^4}\frac{y^4}{s}
-\frac{b^3\delta (p+1)}{6(p-1)^6}\left(\delta+i(1-2p)\right)\frac{y^6}{s^{3/2}}
+O(\frac{y^8}{s^{2}})  \right)^2\\

=\kappa^2\Big(1+\frac{2a}{\kappa\sqrt{s}}(1+i\delta)-\frac{ 2b}{(p-1)^2}(1+i\delta)\frac{y^2}{s^{1/2}}+i\frac{b^2\delta(p+1)}{(p-1)^4}\frac{y^4}{s}
+\frac{ a^2}{\kappa^2 }(1+i\delta)^2\frac1s 
-\frac{2ab}{\kappa (p-1)^2}(1+i\delta)^2\frac{y^2}{s}\\
+\frac{b^2(1+i\delta)^2}{(p-1)^4}\frac{y^4}{s} +\frac{ab^2\delta(p+1)(-\delta+i)}{\kappa(p-1)^4}\frac{y^4}{s^{3/2}}
+\frac{(\delta-i)b^3\delta(p+1)}{(p-1)^6}\frac{y^6}{s^{3/2}}-\frac{b^3\delta (p+1)}{3(p-1)^6}\left(\delta+i(1-2p)\right)\frac{y^6}{s^{3/2}}
+O(\frac{y^8}{s^{2}})\Big)\\

=\kappa^2\Big(1+(1+i\delta)\frac{2a}{\kappa}\frac{1}{\sqrt{s}}+\left((1-p)+i2\delta\right)\frac{a^2}{\kappa^2}\frac 1 s
-(1+i\delta)\frac{ 2b}{(p-1)^2}\frac{y^2}{s^{1/2}}
-\left((1-p)+i2\delta\right)\frac{2ab}{\kappa (p-1)^2}\frac{y^2}{s}\\
+\left(1-p+i\delta(p+3)\right)\frac{b^2}{(p-1)^4}\frac{y^4}{s}
+\frac{ab^2(p+1)(-p+i\delta )}{\kappa(p-1)^4}\frac{y^4}{s^{3/2}}
+\frac 23\frac{b^3(p+1)}{(p-1)^6}(p+i\delta (p-2))\frac{y^6}{s^{3/2}}
+O(\frac{y^6}{s^{4}})+O(\frac{y^8}{s^{4}})
\Big)
\end{array}
\]

\[
\begin{array}{ll}
W_{2,1}&=-(1+i\delta)\frac{p-1}{2}\frac{b}{(p-1)^3}\left(p-1+2i\delta\right)h_2(y)\\
F_2(W_{2,1}\tilde q_2\bar{\tilde{h_2}})&=-32(\left(p-1+2i\delta\right)(p+1)\frac{b}{(p-1)^2}\tilde q_2\\
\tilde P_2(W_{2,1}\tilde q_2\bar{\tilde{h_2}})&=-32\frac{b(p+1)}{p-1}\tilde q_2,\\
F_0(W_{2,1}\tilde q_2\bar{\tilde{h_2}})&=-4(\left(p-1+2i\delta\right)(p+1)\frac{b}{(p-1)^2}\tilde q_2\\
\hat P_0(W_{2,1}\tilde q_2\bar{\tilde{h_2}})&=4\delta\frac{b(p+1)(p-3)}{(p-1)^2}\tilde q_2,
\end{array}
\]

\[
\begin{array}{lll}
W_{2,2}&=&\frac{p-1}{2}(1+i\delta)\Big \{\frac{a^2}{\kappa^2 s}\left(p^2-4p+1+i2\delta(p-2)\right)\\
&&-y^2\frac{2ab}{\kappa(p-1)^2}\left(p^2-4p+1+i2\delta(p-2)\right)\\
&&+y^4\frac{b^2}{(p-1)^4}\left(p^2-4p+1+i3\delta(p-1)\right)\Big\}\\
&=&(1+i\delta)\frac{b^2}{2(p-1)^3}\Big \{4\left(p^2-4p+1+i2\delta(p-2)\right)\\
&&-4y^2\left(p^2-4p+1+i2\delta(p-2)\right)\\
&&+y^4\left(p^2-4p+1+i3\delta(p-1)\right)\Big\}.\\
&=&(1+i\delta)\frac{b^2}{2(p-1)^3}\Big[(p^2-4p+1)h_{2}^{2}(y)\\
&&+i\delta\left (8(p-2)(1-y^2)+y^43(p-1)\right)\Big],

\end{array}
\]
Then, we can write
\[
\begin{array}{lll}
W_{2,2}\tilde q_2\bar{\tilde{h_2}}&=&(p+1)\frac{b^2}{2(p-1)^3}\Big[(p^2-4p+1)h_{2}^{3}(y)+i\delta\left (8(p-2)(1-y^2)+y^43(p-1)\right)h_2(y)\Big]
\end{array}
\]
and we obtain
\[\tilde P_2(W_{2,2}\tilde q_2\bar{\tilde{h_2}})=60(p+1)\frac{b^2}{2(p-1)^3}(p^2-4p+1)\tilde q_2.\]

%-\frac{2ab}{\kappa}\frac{y^2}{s}+\frac{a^2b(p-3)(1-2p)}{\kappa^2(p-1)}\frac{y^2}{s^{3/2}}\\

%+\frac{b^2}{(p-1)^2}\frac{y^4}{s}

%+\frac{ab^2}{2\kappa (p-1)^3}\left[ (p+1)(p-2)+2(p-1)(p-3)+2(p-1)^2\right ]\frac{y^4}{s^{3/2}}\\ 
%-\frac{b^3}{(p-1)^3}\frac{y^6}{s^{3/2}}
%+O(\frac{1}{s^{2}})+O(\frac{y^8}{s^{2}})\Big).\\

%\end{array}
%\]
\section{Details of expansions of the sixth term of equation \eqref{eqqd1} $R^*(\theta',y,s)$ }\label{Detailscal}
Using the definition of $\varphi$, the fact that $\varphi_0$ satisfies \eqref{eqfi0} and \eqref{eqq}, we see that $R^*$ is in fact a function of $\theta'$, $z=\frac{y}{s^{1/4}}$ and $s$ that can be written as 
\beqtn
\begin{array}{l}
\displaystyle R^*(\theta',y,s)=\frac{1}{4}\frac{z}{s}\nabla_z \varphi_0(z)+\frac{1}{2}\frac{a}{s^{3/2}}(1+i\delta)+ \frac{1}{s^{1/2}}\Delta_z \varphi_0(z)\\
\displaystyle -\frac{a(1+i\delta)^2}{(p-1)s^{1/2}}+\left( F\left(\varphi_0(z)+\frac{a}{s^{1/2}}(1+i\delta)\right)-F(\varphi_0)\right)\\
\displaystyle -i\frac{\mu}{s}\left( \varphi_0(z)+\frac{a}{s^{1/2}} (1+i\delta)\right)-i\theta'(s)\left( \varphi_0(z)+\frac{a}{s^{1/2}} (1+i\delta)\right)\mbox{ with, }F(u)=(1+i\delta)|u|^{p-1}u.
\end{array}
\label{eqreste}
\eeqtn
%Since $|z|<1$, there exists positive $C_0$ and $s_0$ such that $|\varphi_0(z)|$ 
In the following, we will try to expand each term of \eqref{eqreste} as a power series of $\dsp\frac 1s$ as $s\to \infty$, uniformly for $\dsp |y|\leq C\dsp s^{\frac 14}$, with $C$ a positive constant.
\begin{itemize}
\item $\frac{1}{4}\frac{z}{s}\nabla_z \varphi_0(z)$:\\
\[
\begin{array}{lll}
\frac{1}{4}\frac{z}{s}\nabla_z \varphi_0(z)&=&-\frac 12\kappa\frac{1+i\delta}{p-1}\frac{b}{p-1}\frac{z^2}{s}\left(1+\frac{b}{p-1}z^2\right)^{-\frac{p+i\delta}{p-1}}\\
&=& -\frac 12\kappa\frac{1+i\delta}{p-1}\frac{b}{p-1}\frac{y^2}{s^{3/2}}\left \{1-\frac{p+i\delta}{p-1}\frac{b}{p-1}\frac{y^2}{s^{1/2}} +O(\frac{y^4}{s})\right\}\\
&=&\dsp -\frac 12\kappa(1+i\delta)\frac{b}{(p-1)^2}\frac{y^2}{s^{3/2}}+O(\frac{y^4}{s^2})+O(\frac{y^6}{s^2}).
\end{array}
\]

\item $ \frac{1}{s^{1/2}}\Delta_z \varphi_0(z)$:

\[
\begin{array}{l}
\frac{1}{s^{1/2}}\Delta_z \varphi_0(z)=\\
\dsp -2\kappa(1+i\delta)\frac{b}{(p-1)^2}\frac{1}{s^{1/2}}\left(1+\frac{b}{p-1}z^2\right)^{-\frac{p+i\delta}{p-1}}+4
\kappa\frac{1+i\delta}{p-1}\frac{p+i\delta}{p-1}\frac{b^2}{(p-1)^2}\frac{z^2}{s^{1/2}}\left(1+\frac{b}{p-1}z^2\right)^{-\frac{2p-1+i\delta}{p-1}}\\

=\dsp 2\kappa(1+i\delta)\frac{b}{(p-1)^2}\frac{1}{s^{1/2}}\Big\{-1+\frac{b(p+i\delta)}{(p-1)^2} \frac{y^2}{s^{1/2}}-\frac{b^2(p+i\delta)(2p-1+i\delta)}{2(p-1)^4} \frac{y^4}{s}
\\
\dsp +2\frac{p+i\delta}{p-1}\frac{b}{p-1}\frac{y^2}{s^{1/2}}\left(1-\frac{b(2p-1+i\delta)}{(p-1)^2}\frac{y^2}{s^{1/2}} \right) +O(\frac{y^6}{s^{3/2}})\Big\}\\
=\dsp -2\kappa(1+i\delta)\frac{b}{(p-1)^2}\frac{1}{s^{1/2}}+i 6\kappa \frac{b^2\delta(p+1)}{(p-1)^2}\frac{y^2}{s}+5\kappa\frac{b^3\delta(p+1)}{(p-1)^6}\left (\delta-i(2p-1)\right)\frac{y^4}{s^{3/2}}+O(\frac{y^6}{s^2}).
\end{array}
\]

\item $F\left(\varphi_0(z)+\frac{a}{s^{1/2}}(1+i\delta)\right)$
We start by the following term and recall that we work in the critical case $\delta^2=p$ and  we expand this term as a power series of $\dsp\frac 1s$ as $s\to \infty$, uniformly for $\dsp |y|\leq C \dsp s^{\frac 14}$, with $C$ a positive constant.

\medskip

\[
\begin{array}{l}
\left |\varphi_0(z)+\frac{a}{s^{1/2}}(1+i\delta)\right|^{p-1}\\
=\left |\kappa+ \frac{a}{s^{1/2}}(1+i\delta)-\kappa\frac{b(1+i\delta)}{(p-1)^2} \frac{y^2}{s^{1/2}}+i\kappa \frac{b^2\delta (p+1)}{2(p-1)^4}\frac{y^4}{s}  
+\kappa \frac{b^3}{6(p-1)^6}\delta(1+p) (\delta-i(2p-1))\frac{y^6}{s^{3/2}}
+O(\frac{y^8}{s^{2}})\right|^{p-1}\\

= [\left (\kappa+ \frac{a}{s^{1/2}}-\frac{\kappa b}{(p-1)^2}\frac{y^2}{s^{1/2}}+\kappa \frac{b^3\delta^2 (1+p)}{6(p-1)^6} \frac{y^6}{s^{3/2}}
+O(\frac{y^8}{s^{2}})\right)^2+\delta^2\left(\frac{a}{s^{1/2}}-\frac{\kappa b}{(p-1)^2}\frac{y^2}{s^{1/2}} +\kappa \frac{b^2 (p+1)}{2(p-1)^4} \frac{y^4}{s}+O(\frac{y^6}{s^{3/2}} \right)^2         ]^{\frac{p-1}{2}}\\

=\Big ((\kappa+ \frac{a}{s^{1/2}})^2+\frac{\delta^2 a^2}{s}-2(\kappa+ \frac{a}{s^{1/2}})\frac{\kappa b}{(p-1)^2}\frac{y^2}{s^{1/2}}-2\frac{\delta^2 \kappa ab}{(p-1)^2}\frac{y^2}{s}+\frac{\kappa^2b^2(p+1)}{(p-1)^4}\frac{y^4}{s}\\
+\kappa\frac{ab^2\delta^2(p+1)}{(p-1)^4}\frac{y^4}{s^{3/2}}+\kappa^2\frac{b^3\delta^2(p+1)}{(p-1)^6}\frac{1-3}{3}\frac{y^6}{s^{3/2}}+O(\frac{y^6}{s^{2}})+O(\frac{y^8}{s^{2}})\Big)^{\frac{p-1}{2}}\\

=\Big (\kappa^2+\frac{2\kappa a}{s^{1/2}}+\frac{a^2(p+1)}{s}-2\frac{b\kappa^2 }{(p-1)^2}\frac{y^2}{s^{1/2}}-2\frac{\kappa ab(p+1)}{(p-1)^2}\frac{y^2}{s}+\frac{\kappa^2b^2(p+1)}{(p-1)^4}\frac{y^4}{s} 
+\kappa\frac{ab^2 p(p+1)}{(p-1)^4}\frac{y^4}{s^{3/2}}-\frac 2 3\kappa^2\frac{b^3p(p+1)}{(p-1)^6}\frac{y^6}{s^{3/2}}
\\+O(\frac{y^6}{s^{2}})+O(\frac{y^8}{s^{2}})\Big)^{\frac{p-1}{2}}.\\

=\kappa^{p-1}\left(1+\frac{2 a}{\kappa s^{1/2}}+\frac{a^2(p+1)}{\kappa^2s}\right)^{\frac{p-1}{2}}\Big(1-2\frac{b}{(p-1)^2}\frac{1}{1+\frac{2 a}{\kappa s^{1/2}}+\frac{a^2(p+1)}{\kappa^2s}}\Big[\frac{y^2}{s^{1/2}}+\frac{a(p+1)}{\kappa}\frac{y^2}{s}-\frac{b(p+1)}{2(p-1)^2}\frac{y^4}{s}\\
-\frac{ab p(p+1)}{2\kappa(p-1)^2}\frac{y^4}{s^{3/2}}+\frac 13\frac{b^2 p(p+1)}{(p-1)^4} \frac{y^6}{s^{3/2}}
+O(\frac{y^6}{s^{2}})+O(\frac{y^8}{s^{2}}) \Big ]\Big)^{\frac{p-1}{2}}\\

=\kappa^{p-1}\left(1+\frac{2 a}{\kappa s^{1/2}}+\frac{a^2(p+1)}{\kappa^2s}\right)^{\frac{p-1}{2}}\Big (1-2\frac{b}{(p-1)^2}\left[ 1-\frac{2 a}{\kappa s^{1/2}}-\frac{a^2(p+1)}{\kappa^2s}   +\frac{4 a^2}{\kappa^2 s}+O(\frac{1}{s^{3/2}})\right ]\\
*\Big[\frac{y^2}{s^{1/2}}+\frac{a(p+1)}{\kappa}\frac{y^2}{s}-\frac{b(p+1)}{2(p-1)^2}\frac{y^4}{s}
-\frac{ab p(p+1)}{2\kappa (p-1)^2}\frac{y^4}{s^{3/2}}+\frac 13\frac{b^2 p(p+1)}{(p-1)^4}\frac{y^6}{s^{3/2}}
+O(\frac{y^6}{s^{2}})+O(\frac{y^8}{s^{2}}) \Big ]\Big)^{\frac{p-1}{2}}\\

=\kappa^{p-1}\left(1+\frac{2 a}{\kappa s^{1/2}}+\frac{a^2(p+1)}{\kappa^2s}\right)^{\frac{p-1}{2}}\Big (1-2\frac{b}{(p-1)^2}\left[ 1-\frac{2 a}{\kappa s^{1/2}}-\frac{a^2(p-3)}{\kappa^2s} +O(\frac{1}{s^{3/2}})\right ]\\
*\left[\frac{y^2}{s^{1/2}}+\frac{a(p+1)}{\kappa}\frac{y^2}{s} -\frac{b(p+1)}{2(p-1)^2}\frac{y^4}{s}-\frac{ab p(p+1)}{2\kappa(p-1)^2}\frac{y^4}{s^{3/2}}+\frac 13\frac{b^2 p(p+1)}{(p-1)^4} \frac{y^6}{s^{3/2}}
+O(\frac{y^6}{s^{2}})+O(\frac{y^8}{s^{2}})\right]\Big )^{\frac{p-1}{2}}\\

=\kappa^{p-1}\left(1+\frac{2 a}{\kappa s^{1/2}}+\frac{a^2(p+1)}{\kappa^2s}\right)^{\frac{p-1}{2}}\Big (1-2\frac{b}{(p-1)^2}\frac{y^2}{s^{1/2}}-2\frac{ab}{\kappa(p-1)}\frac{y^2}{s}+\frac{b^2(p+1)}{(p-1)^4}\frac{y^4}{s}\\

+2\frac{a^2b (3p-1)}{\kappa^2(p-1)^2}\frac{y^2}{s^{3/2}}+\frac{ab^2(p+1)( p-2)}{\kappa(p-1)^4}\frac{y^4}{s^{3/2}}
-\frac 23\frac{b^3p (p+1)}{(p-1)^6}\frac{y^6}{s^{3/2}}
+O(\frac{y^6}{s^{2}})+O(\frac{y^8}{s^{2}})\Big)^{\frac{p-1}{2}}\\

=\kappa^{p-1}\left(1+\frac{(p-1)a}{\kappa s^{1/2}}+\frac{a^2(p^2-1)}{2\kappa^2s}+\frac{(p-1)(p-3)a^2}{2\kappa^2 s}+
+\frac{a^3}{2 \kappa^3}\frac{(p^2-1)(p-3)}{s^{3/2}}+\frac{a^3}{6\kappa^3}\frac{(p-1)(p-3)(p-5)}{s^{3/2}}+

O(\frac{1}{s^{2}})\right)\\

*\Big(1-\frac{b}{(p-1)}\frac{y^2}{s^{1/2}}-\frac{ab}{\kappa}\frac{y^2}{s}+\frac{b^2(p+1)}{2(p-1)^3}\frac{y^4}{s}
+\frac{(p-3)b^2}{2(p-1)^3}\frac{y^4}{s}\\
+\frac{a^2b (3p-1)}{\kappa^2(p-1)}\frac{y^2}{s^{3/2}}
+\frac{ab^2(p+1)(p-2)}{2\kappa(p-1)^3}\frac{y^4}{s^{3/2}}+\frac{ab^2(p-1)(p-3)}{\kappa (p-1)^3}\frac{y^4}{s^{3/2}}\\
-\frac{b^3}{(p-1)^5}\left [\frac 1 3p(p+1)+\frac 1 2(p+1)(p-3)+\frac 16 (p-3)(p-5)

\right ]\frac{y^6}{s^{3/2}}

+O(\frac{y^6}{s^{2}})+O(\frac{y^8}{s^{2}})\Big)\\

\\
=\kappa^{p-1}\left(1+\frac{(p-1)a}{\kappa s^{1/2}}+\frac{a^2(p-1)^2}{\kappa^2s}+\frac{a^3}{3\kappa^3 s^{3/2}}(p-1)(p-3)(2p-1)+O(\frac{1}{s^{2}})\right)\\

*\Big(1-\frac{b}{(p-1)}\frac{y^2}{s^{1/2}}-\frac{ab}{\kappa}\frac{y^2}{s}+\frac{b^2}{(p-1)^2}\frac{y^4}{s}\\

+\frac{a^2b (3p-1)}{\kappa^2(p-1)}\frac{y^2}{s^{3/2}}\\
+\frac{ab^2}{2\kappa(p-1)^3}\left[(p+1)(p-2)+2(p-1)(p-3)\right]\frac{y^4}{s^{3/2}}\\
-\frac{b^3}{(p-1)^3}\frac{y^6}{s^{3/2}}
+O(\frac{y^6}{s^{2}})+O(\frac{y^8}{s^{2}})\Big ),\\

\\
\\

=\kappa^{p-1}\Big(1+\frac{(p-1)a}{\kappa s^{1/2}}+\frac{a^2(p-1)^2}{\kappa^2s}-\frac{b}{(p-1)}\frac{y^2}{s^{1/2}}+\frac{a^3}{3\kappa^3 s^{3/2}}(p-1)(p-3)(2p-1)
\end{array}
\]
Now we can write

\[
\begin{array}{l}
F\left(\varphi_0(z)+\frac{a}{s^{1/2}}(1+i\delta)\right)\\
=  (1+i\delta)\kappa^{p-1}
\Big(1+\frac{(p-1)a}{\kappa s^{1/2}}+\frac{a^2(p-1)^2}{\kappa^2s}-\frac{b}{(p-1)}\frac{y^2}{s^{1/2}}+\frac{a^3}{3\kappa^3 s^{3/2}}(p-1)(p-3)(2p-1)\\

-\frac{2ab}{\kappa}\frac{y^2}{s}+\frac{a^2b(p-3)(1-2p)}{\kappa^2(p-1)}\frac{y^2}{s^{3/2}}\\

+\frac{b^2}{(p-1)^2}\frac{y^4}{s}

+\frac{ab^2}{2\kappa (p-1)^3}(p-2)(5p-3)\frac{y^4}{s^{3/2}}\\
-\frac{b^3}{6(p-1)^3}\frac{y^6}{s^{3/2}}
+O(\frac{1}{s^{2}})+O(\frac{y^8}{s^{2}})\Big).\\

*\left[\kappa +\frac{a(1+i\delta)}{s^{1/2}}-\frac{\kappa b(1+i\delta)}{(p-1)^2}\frac{y^2}{s^{1/2}} +i\kappa\frac{b^2\delta(p+1)}{2(p-1)^4}\frac{y^4}{s}
-\kappa\frac{b^3}{6(p-1)^6}\delta(1+p)\left(\delta+i(1-2p)\right)\frac{y^6}{s^{3/2}}
+O(\frac{y^8}{s^{2}})\right].
\end{array}
\]

\item $F(\varphi_0(z))$

\[
\begin{array}{l}
F(\varphi_0(z))=\\
 (1+i\delta)\kappa^p \left | 1+\frac{b}{p-1}\frac{y^2}{s^{1/2}}\right |^{-1}\left(1+\frac{b}{p-1}\frac{y^2}{s^{1/2}}\right)^{-\frac{1+i\delta}{p-1}}\\
=(1+i\delta) \kappa^p \left(1-\frac{b}{p-1}\frac{y^2}{s^{1/2}}+\frac{b^2}{(p-1)^2}\frac{y^4}{s}-\frac{b^3}{(p-1)^3}\frac{y^6}{s^{3/2}}+O(\frac{y^8}{s^{2}})\right)\\

*\left\{1-\frac{b(1+i\delta)}{(p-1)^2}\frac{y^2}{s^{1/2}}+i\frac{b^2\delta(p+1)}{2(p-1)^4}\frac{y^4}{s}-\frac{b^3}{6(p-1)^6}\delta(1+p)\left(\delta+i(1-2p)\right)\frac{y^6}{s^{3/2}}
+O(\frac{y^8}{s^{2}})\right\}\\
= (1+i\delta)\kappa^p\Big(1-\frac{b}{p-1}\frac{p+i\delta}{p-1}\frac{y^2}{s^{1/2}}+\left \{ i\frac{b^2\delta(p+1)}{2(p-1)^4}+\frac{b^2}{(p-1)^2}+\frac{b^2(1+i\delta)}{(p-1)^3} \right \}\frac{y^4}{s}\\
-\left \{\frac{b^3}{(p-1)^3}+\frac{b^3}{6(p-1)^6}\delta(1+p)\left(\delta+i(1-2p)\right)+i\frac{b^3\delta(p+1)}{2(p-1)^5}+\frac{b^3}{(p-1)^4}(1+i\delta)
\right \}\frac{y^6}{s^{3/2}}

+O(\frac{y^8}{s^{2}})\Big).
\end{array}
\]
\item Finally the term$\varphi_0(z)+\frac{a}{s^{1/2}} (1+i\delta)$:

\[
%\begin{array}{lll}
\varphi_0(z)+\frac{a}{s^{1/2}} (1+i\delta)=\kappa+\frac{a}{\sqrt{s}}(1+i\delta)-\frac{\kappa b}{(p-1)^2}(1+i\delta)\frac{y^2}{s^{1/2}}+i\kappa\frac{b^2\delta(p+1)}{2(p-1)^4}\frac{y^4}{s}+O(\frac{y^6}{s^{3/2}}).
%\end{array}
\]
\end{itemize}

Now, let us see the term of order $\frac 1 s$ in the function $R^*$, by expansions above, we have
\[\begin{array}{l}
(1+i\delta)\frac{a^2}{\kappa}\left (p-1+1+i\delta\right)-i\mu\kappa,\\
=\frac{a^2}{\kappa}(1+i\delta)(p+i\delta )-i\mu\kappa,\\
=\frac{a^2}{\kappa}(p-\delta^2+i\delta(p+1))-i\mu\kappa\mbox{  we recall that we are in the critical case ,}p=\delta^2,\\
=i \frac{a^2}{\kappa}\delta (1+p)-i\mu\kappa%\mbox{  we recall that,}\mu=\frac{\delta}{8p},\;\;a=\frac{\kappa}{4\sqrt{p(p+1)}},\\
%=i \kappa\frac{1}{16}\frac{\delta}{p}-i\kappa\frac{1}{8}\frac{\delta}{p}.
%\mbox{\textcolor{red}{\danger we must change $\mu=\frac{\delta}{16p}$.}}
\end{array}
\]

Now, let us see the term of order $\frac{1}{s^{1/2}}$ in the function $R^*$
\[\begin{array}{l}
(1+i\delta)\left\{-2\frac{\kappa b}{(p-1)^2}-(1+i\delta)\frac{a}{(p-1)}+a+(1+i\delta)\frac{a}{(p-1)}\right\}\\
=(1+i\delta)\left\{-2\frac{\kappa b}{(p-1)^2}+a\right\}=0\mbox{, (we recall }b=\frac{(p-1)^2}{8\sqrt{p(p+1)}},\;\;\; a=\frac{\kappa}{4\sqrt{p(p+1)}}).
\end{array}
\]

Now, let us see the term of  order $\frac{y^2}{s^{1/2}}$ in the function $R^*$, by expansions above, we have
\[-2i\frac{b\kappa}{(p-1)^2}(p+1)\delta+2i\frac{b\kappa}{(p-1)^2}(p+1)\delta=0\]
%+(1+i\delta)\theta'(s)\frac{\kappa b}{(p-1)^2} =(1+i\delta)\theta'(s)\frac{\kappa b}{(p-1)^2}\]
Now, let us see the term of  order $\frac{y^2}{s}$ in the function $R^*$, by expansions above, we have
%and recalling that $b=\frac{(p-1)^2}{8\sqrt{p(p+1)}},\;\;\; a=\frac{\kappa}{4\sqrt{p(p+1)}}$, we obtain
\[i6\kappa\frac{\delta b^2(p+1)}{(p-1)^4}-i\frac{2ab\delta (p+1)}{(p-1)^2}\]%=i6\kappa\frac{\delta}{64p}-i2\kappa\frac{\delta}{32p}=i\kappa\frac{\delta}{32p}\]
The term of order $\frac{y^4}{s}$ is equal to zero.

\medskip

Now, let us see the term of order $\frac{1}{s^{3/2}}$ in the function $R^*$, by expansions above, we have
\[ a(\frac 12 +\mu\delta)+\frac{a^3}{3\kappa^2}(p-3)(2p-1)+i\left[ \frac{\delta a^3}{3\kappa^2}(p-3)(2p-1)-\mu a\right ].\]
The term of order $\frac{y^2}{s^{3/2}}$ in the function $R^*$, 
\[\kappa\frac{b}{(p-1)^2}\left(-\frac 12 -\mu\delta\right)
+\frac{a^2b}{\kappa (p-1)^2}p(p+1)
+i\left[\kappa\frac{b}{(p-1)^2}\left(-\frac \delta2 +\mu\right)
+\frac{a^2b}{\kappa (p-1)^2}\delta(-2p^2+p+3)
\right]
\]
The term of order $\frac{y^4}{s^{3/2}}$ in the function $R^*$, 
\[\kappa\frac{b^3}{(p-1)^6}5p(p+1)-\frac{ab^2}{(p-1)^4}p(p+1)+i\delta \left[\kappa\frac{b^3}{(p-1)^6}5(1-2p)+\frac{ab^2}{(p-1)^4}(5p^2-5)\right]
\]
The term of order $\frac{y^6}{s^{3/2}}$ in the function $R^*$, is equal to zero.
\section{Taylor expansion of $B$}
\label{TEB}
\medskip

Let us recall from \eqref{eqqd} that:
\[B(q,y,s)=(1+i\delta)\left ( |\varphi+q|^{p-1}(\varphi+q)-|\varphi|^{p-1}\varphi-|\varphi|^{p-1} q-\frac{p-1}{2}|\varphi|^{p-3}\varphi(\varphi \bar q+\bar\varphi q)\right).\]

Next, we compute the taylor expansion of $F(q)= |\varphi+q|^{p-1}= \left((\varphi_1+q_1)^2+(\varphi_2+q_2)^2\right)^{(p-1)/2}$, then an easy calculation shows
\[\nabla F(q)=(p-1) |\varphi+q|^{p-3}\left(\varphi_1+q_1,\;\;\varphi_2+q_2 \right)
\]
where $\varphi=\varphi_1+i\varphi_2$ and $q=q_1+iq_2$, and
\[H(F)(q)
=\dsp (p-1)\left(
\begin{array}{lll}
 \dsp|\varphi+q|^{p-3}+(p-3)  |\varphi+q|^{p-5}(\varphi_1+q_1)^2  && \dsp(p-3) |\varphi+q|^{p-5}(\varphi_1+q_1)(\varphi_2+q_2)\\
 \dsp (p-3) |\varphi+q|^{p-5}(\varphi_1+q_1)(\varphi_2+q_2) && \dsp |\varphi+q|^{p-3}+(p-3)  |\varphi+q|^{p-5}(\varphi_2+q_2)^2.  
 \end{array}
\right)
\]
Then 
\[\nabla F(0)=(p-1) |\varphi|^{p-3}\left(\varphi_1,\;\;\varphi_2 \right)
\]
\[H(F)(0)
=\dsp  (p-1)\left(
\begin{array}{lll}
 \dsp|\varphi|^{p-3}+ (p-3)   |\varphi|^{p-5}\varphi_{1}^{2}  && \dsp (p-3)  |\varphi|^{p-5} \varphi_1\varphi_2\\
 \dsp  \frac{(p-3)}{2}  |\varphi|^{p-5}\varphi_1 \varphi_2 && \dsp |\varphi|^{p-3}+ (p-3)   |\varphi|^{p-5}\varphi_{2}^{2}.  
 \end{array}
\right)
\]
And we obtain

\[
\begin{array}{lll}
B(q,y,s)&=&(1+i\delta)\Big \{(p-1)|\varphi |^{p-3}  (\varphi_1 q_1+\varphi_2 q_2 )q\\
&&+\frac{(p-1)}{2}\Big [(|\varphi|^{p-3}+ (p-3)   |\varphi|^{p-5}\varphi_{1}^{2}) \varphi q_{1}^{2}+ (|\varphi|^{p-3}+ (p-3)  |\varphi|^{p-5}\varphi_{2}^{2})\varphi q_{2}^{2}\\
&&+2(p-3)  |\varphi|^{p-5}\varphi \varphi_1 \varphi_2 q_1q_2 )
\Big ]\Big \}+O(q^3)
\end{array}
\]
Using the fact that $\varphi_1=\kappa+ O(\frac{y^2}{\sqrt s})$, $\varphi_2= O(\frac{y^2}{\sqrt s})$ and by the decomposition of $q$ given by \eqref{decompq}, we can deduce that the contribution of  $\tilde q_{2}^{2}$ is 
given by the following

\[
\begin{array}{l}
 (1+i\delta)\left\{(p-1)  \kappa^{p-2}(1+i\delta)+\frac{p-1}{2}\left(1+(p-3)+\delta^2 \right )\right \} \tilde q_{2}^{2} h_{2}^{2} ,\\

 %\left(1+p-3+\delta^2  \right)\right\}\tilde q_{2}^{2}\\

=\dsp(p-1)\kappa^{p-2}\Big[(1+i\delta)^2 +(1+i\delta) (p-1)\Big ]\tilde q_{2}^{2} h_{2}^{2}
\end{array}
\]
We note that 
\[ \tilde P_2\left((1+i\delta)^2 h_{2}^{2}\right)= 8(1-\delta^2)=8(1-p).\]
\[ \tilde P_2\left((1+i\delta) h_{2}^{2}\right)= 8.\]

Then, we conclude that the contribution of $\tilde q_{2}^{2}$ in $\tilde P_2(B)$ is zero.
%\[\tilde P_2(B)=0*\tilde q_{2}^{2}+......\] 
\newpage
\textbf{Funding:}\\
Hatem Zaag is supported by the ERC Advanced Grant  no. 291214, BLOWDISOL, and by ANR project no. ANR-13-BS01-0010-03, ANA\'E.

\bigskip

\textbf{Conflict of Interest:}\\
The authors declare that they have no conflict of interest.

\newpage
%\bibliography{mabiblio}

\begin{thebibliography}{PSKK98}

\bibitem[AK02]{AKRMP02}
I.S. Aranson and L.~Kramer.
\newblock The world of the complex {G}inzburg-{L}andau equation.
\newblock {\em Rev. Modern Phys.}, 74(1):99--143, 2002.

\bibitem[BK94]{BKN94}
J.~Bricmont and A.~Kupiainen.
\newblock Universality in blow-up for nonlinear heat equations.
\newblock {\em Nonlinearity}, 7(2):539--575, 1994.

\bibitem[Caz03]{CNYUCIM03}
T.~Cazenave.
\newblock {\em Semilinear {S}chr\"odinger equations}, volume~10 of {\em Courant
  Lecture Notes in Mathematics}.
\newblock New York University, Courant Institute of Mathematical Sciences, New
  York; American Mathematical Society, Providence, RI, 2003.

\bibitem[CZ13]{CZCPAM13}
R.~C{\^o}te and H.~Zaag.
\newblock Construction of a multi-soliton blow-up solution to the semilinear
  wave equation in one space dimension.
\newblock {\em Comm. Pure Appl. Math.}, 66(10):1541--1581, 2013.

\bibitem[EZ11]{EZSMJ11}
M.~A. Ebde and H.~Zaag.
\newblock Construction and stability of a blow up solution for a nonlinear heat
  equation with a gradient term.
\newblock {\em S$\vec{ \rm e}$MA J.}, 55:5--21, 2011.

\bibitem[GK89]{GKCPAM89}
Y.~Giga and R.~V. Kohn.
\newblock Nondegeneracy of blowup for semilinear heat equations.
\newblock {\em Comm. Pure Appl. Math.}, 42(6):845--884, 1989.

\bibitem[GV96]{GV96}
J.~Ginibre and G.~Velo.
\newblock The {C}auchy problem in local spaces for the complex
  {G}inzburg-{L}andau equation.
\newblock {\em Differential equations, asymptotic analysis, and mathematical
  physics}, 100:138--152, 1996.

\bibitem[GV97]{GVCMP97}
J.~Ginibre and G.~Velo.
\newblock The {C}auchy problem in local spaces for the complex
  {G}inzburg-{L}andau equation. {II}. {C}ontraction methods.
\newblock {\em Comm. Math. Phys.}, 187(1):45--79, 1997.

\bibitem[HS72]{HSPRSLS72}
L.~M. Hocking and K.~Stewartson.
\newblock On the nonlinear response of a marginally unstable plane parallel
  flow to a two-dimensional disturbance.
\newblock {\em Proc. Roy. Soc. London Ser. A}, 326:289--313, 1972.

\bibitem[HSSB72]{HSSBJFM72}
L.~M. Hocking, K.~Stewartson, Stuart.J.T., and Brown.S.N.
\newblock A nonlinear instability in plane parallel flow.
\newblock {\em J. Fluid. Mech}, 51:705--735, 1972.

\bibitem[KBS88]{KBSPRL88}
P.~Kolodner, D.~Bensimon, and M.~Surko.
\newblock Traveling wave convection in an annulus.
\newblock {\em Phys.Rev.Lett.}, 60:1723--, 1988.

\bibitem[KSAL95]{KSALPNP95}
P.~Kolodner, S.~Slimani, N.~Aubry, and R.~Lima.
\newblock Characterization of dispersive chaos and related states of
  binary-fluid convection.
\newblock {\em Phys. D}, 85(1-2):165--224, 1995.

\bibitem[Mer92]{FMCPAM92}
F.~Merle.
\newblock Solution of a nonlinear heat equation with arbitrary given blow-up
  points.
\newblock {\em Comm. Pure Appl. Math.}, 45(3):263--300, 1992.

\bibitem[MZ97]{MZDuke97}
F.~Merle and H.~Zaag.
\newblock Stability of the blow-up profile for equations of the type
  $u_t={\Delta} u+\vert u\vert ^{p-1}u$.
\newblock {\em Duke Math. J.}, 86(1):143--195, 1997.

\bibitem[MZ08]{MZ07}
N.~Masmoudi and H.~Zaag.
\newblock Blow-up profile for the complex {G}inzburg-{L}andau equation.
\newblock {\em J. Funct. Anal.}, 225:1613--1666, 2008.

\bibitem[NZ10]{NZTAMS08}
N.~Nouaili and H.~Zaag.
\newblock A liouville theorem for vector valued semilinear heat equations with
  no gradient structure and applications to blow-up.
\newblock {\em Trans. Amer. Math. Soc.}, 362(7):3391--3434, 2010.

\bibitem[NZ15a]{NVTZ15}
V.T. Nguyen and H.~Zaag.
\newblock Construction of a stable blow-up solution for a class of strongly
  perturbed semilinear heat equations.
\newblock 2015.
\newblock Submitted.

\bibitem[NZ15b]{NZCPDE15}
N.~Nouaili and H.~Zaag.
\newblock Profile for a simultaneously blowing up solution to a complex valued
  semilinear heat equation.
\newblock {\em Comm. Partial Differential Equations}, (7):1197--1217, 2015.

\bibitem[P{\v{S}}01]{PSCPAM01}
P.~Plech{\'a}{\v{c}} and V.~{\v{S}}ver{\'a}k.
\newblock On self-similar singular solutions of the complex {G}inzburg-{L}andau
  equation.
\newblock {\em Comm. Pure Appl. Math.}, 54(10):1215--1242, 2001.

\bibitem[PSKK98]{PSKKPD98}
S.~Popp, O.~Stiller, E.~Kuznetsov, and L.~Kramer.
\newblock The cubic complex {G}inzburg-{L}andau equation for a backward
  bifurcation.
\newblock {\em Phys. D}, 114(1-2):81--107, 1998.

\bibitem[Rot08]{RPD08}
V.~Rottsch\"afer.
\newblock Multi-bump, self-similar, blow-up solutions of the
  {G}inzburg-{L}andau equation.
\newblock {\em Phys. D}, 237(4):510--539, 2008.

\bibitem[Rot13]{REJAM13}
V.~Rottsch\"afer.
\newblock Asymptotic analysis of a new type of multi-bump, self-similar, blowup
  solutions of the {G}inzburg-{L}andau equation.
\newblock {\em European J. Appl. Math.}, 24(1):103--129, 2013.

\bibitem[RS13]{RSCPAM13}
P.~Rapha{\"e}l and R.~Schweyer.
\newblock Stable blowup dynamics for the 1-corotational energy critical
  harmonic heat flow.
\newblock {\em Comm. Pure Appl. Math.}, 66(3):414--480, 2013.

\bibitem[RS14]{RSMA14}
P.~Rapha{\"e}l and R.~Schweyer.
\newblock On the stability of critical chemotactic aggregation.
\newblock {\em Math. Ann.}, 359(1-2):267--377, 2014.

\bibitem[SS71]{SSJFM71}
K.~Stewartson and J.~T. Stuart.
\newblock A non-linear instability theory for a wave system in plane
  {P}oiseuille flow.
\newblock {\em J. Fluid Mech.}, 48:529--545, 1971.

\bibitem[TZ15]{TZ17}
S.~Tayachi and H.~Zaag.
\newblock Existence of a stable blow-up profile for the nonlinear heat equation
  with a critical power nonlinear gradient term.
\newblock 2015.
\newblock arXiv preprint:1506.08306.

\bibitem[Vel92]{VCPDE92}
J.~J.~L. Vel{\'a}zquez.
\newblock Higher-dimensional blow up for semilinear parabolic equations.
\newblock {\em Comm. Partial Differential Equations}, 17(9-10):1567--1596,
  1992.

\bibitem[Vel93a]{VTAMS93}
J.~J.~L. Vel{\'a}zquez.
\newblock Classification of singularities for blowing up solutions in higher
  dimensions.
\newblock {\em Trans. Amer. Math. Soc.}, 338(1):441--464, 1993.

\bibitem[Vel93b]{VINDIANA93}
J.~J.~L. Vel{\'a}zquez.
\newblock Estimates on the $(n-1)$-dimensional {H}ausdorff measure of the
  blow-up set for a semilinear heat equation.
\newblock {\em Indiana Univ. Math. J.}, 42(2):445--476, 1993.

\bibitem[VGH91]{VGH91}
J.~J.~L. Vel\'azquez, V.~A. Galaktionov, and M.~A. Herrero.
\newblock The space structure near a blow-up point for semilinear heat
  equations: a formal approach.
\newblock {\em Zh. Vychisl. Mat. i Mat. Fiz.}, 31(no.~3), 1991.

\bibitem[Zaa98]{ZAIHPANL98}
H.~Zaag.
\newblock Blow-up results for vector-valued nonlinear heat equations with no
  gradient structure.
\newblock {\em Ann. Inst. H. Poincar\'e Anal. Non Lin\'eaire}, 15(5):581--622,
  1998.

\bibitem[Zaa01]{ZCPAM01}
H.~Zaag.
\newblock A {L}iouville theorem and blowup behavior for a vector-valued
  nonlinear heat equation with no gradient structure.
\newblock {\em Comm. Pure Appl. Math.}, 54(1):107--133, 2001.

\bibitem[Zaa02a]{ZIHP02}
H.~Zaag.
\newblock On the regularity of the blow-up set for semilinear heat equations.
\newblock {\em Ann. Inst. H. Poincar\'e Anal. Non Lin\'eaire}, 19(5):505--542,
  2002.

\bibitem[Zaa02b]{ZCMP02}
H.~Zaag.
\newblock One-dimensional behavior of singular {$N$}-dimensional solutions of
  semilinear heat equations.
\newblock {\em Comm. Math. Phys.}, 225(3):523--549, 2002.

\bibitem[Zaa02c]{ZMME02}
H.~Zaag.
\newblock Regularity of the blow-up set and singular behavior for semilinear
  heat equations.
\newblock In {\em Mathematics \& mathematics education (Bethlehem, 2000)},
  pages 337--347. World Sci. Publishing, River Edge, NJ, 2002.

\end{thebibliography}
\bibliographystyle{alpha}

\end{document}